%% file: drinfeld.tex
\documentstyle{amsppt}
\nonstopmode
\voffset=+.0in
\hoffset=+.1in
\vsize=7.5in
\magnification=\magstep1
\parindent=12pt
\baselineskip=14pt
\hfuzz=5pt
\mathsurround 2pt
\parskip 6pt
\overfullrule=0pt
\def \today {\ifcase \month \or January\or February\or
  March\or April\or May\or June\or July\or August\or
  September\or October\or November\or December\fi\
  \number \day, \number \year}
\def\now{
\countdef\hours=1
\countdef\minutes=2
\count0=\time
\divide\count0 by 60
\hours=\count0
\count3=\count0
\multiply\count3 by 60
\minutes=\time
\advance\minutes by -\count3
\number\hours:\number\minutes}

\define\A{{\Bbb A}}
\define\C{{\Bbb C}}
\define\PP{{\Bbb P}}    
\define\Q{{\Bbb Q}}
\define\R{{\Bbb R}}
\define\Z{{\Bbb Z}}
\define\a{\alpha}
\redefine\b{\beta}

\define\e{\varepsilon}
\predefine\l{\ll}
\redefine\l{\lambda}
\redefine\o{\omega}
\define\ph{\varphi}

\redefine\P{\Phi}
\predefine\Sec{\S}


\define\back{\backslash}

\define\lra{\longrightarrow}
\redefine\tt{\otimes}
\define\scr{\scriptstyle}
\define\liminv#1{\underset{\underset{#1}\to\leftarrow}\to\lim\,}
\define\limdir#1{\underset{\underset{#1}\to\rightarrow}\to\lim\,}

\define\isoarrow{\ {\overset{\sim}\to{\longrightarrow}}\ }


\parindent=12pt
\baselineskip=14pt
\define\F{\Bbb F}

\define\fm{\frak m}

\define\nass{\noalign{\smallskip}}

\redefine\j{\text{\bf j}}

\define\und#1{\underline{#1}}

\parskip=12pt
\define\End{\text{\rm End}}
\define\Spec{\text{\rm Spec}}


\predefine\oldvol{\vol}
\redefine\vol{\text{\rm vol}}

\define\Sym{\text{\rm Sym}}

\define\mult{\text{\rm mult}}    

\define\diag{\text{\rm diag}}    


\define\hak{\widehat{\Cal A_K}}
\define\hakw{\big(\hak\big)_W}
\define\ha{\hat{\Cal A}}
\define\hc{\hat{\Cal C}}

\define\tha{\ha^{{}\thicksim{}}}
\define\thc{\hc^{{}\thicksim{}}}

\define\uuj{\underline{j}}
\define\bull{\bullet}

\define\pr{\text{\rm pr}}
\define\Spf{\text{\rm Spf}}
\define\ord{\text{\rm ord}}   


\define\AM{\bf1} 
\define\boutotcarayol{\bf2}
\define\deligne{\bf3}
\define\delignetwo{\bf4}
\define\drinfeld{\bf5}
\define\genestierbook{\bf6}   
\define\genestier{\bf7}
\define\gross{\bf8}
\define\grosskeating{\bf9}
\define\katsurada{\bf10}
\define\kitaoka{\bf11}
\define\kitaokaII{\bf12}
\define\kottwitz{\bf13}  
\define\kottwitztwo{\bf14}
\define\annals{\bf15}
\define\krsiegel{\bf16}
\define\krHB{\bf17}
\define\myers{\bf18}
\define\rapoportzink{\bf19}
\define\yang{\bf20}
\define\zink{\bf21}  

\centerline{\bf Height pairings on Shimura curves and $p$-adic uniformization}
\smallskip
\centerline{by}
\smallskip
\centerline{\bf Stephen S. Kudla\footnote{Partially supported by NSF Grant
DMS-9622987}}
\smallskip
\centerline{and}
\smallskip
\centerline{\bf Michael Rapoport}

\vskip 25pt

\noindent
{\bf Introduction.}
\par\noindent
In a recent paper \cite{\annals} of one of us it was shown that
there is a close connection between the value of the height pairing
of certain arithmetic 0-cycles on Shimura curves and the values at
the center of their symmetry of the derivatives of certain
metaplectic Eisenstein series of genus 2. On the one hand, the
height pairing can be written as a sum of local height pairings.
For example, if the 0-cycles have disjoint support on the generic
fiber, then their height pairing is a sum of an archimedean
contribution and a contribution from each of the (finitely many)
finite primes $p$ for which the cycles meet in the fiber at $p$. On
the other hand, it turns out that the non-singular part of the
Fourier expansion of the central derivative of the metaplectic
Eisenstein series also has a decomposition into a sum of
contributions indexed by the places of $\Q$. Then, one would like
to compare the height pairing and the Fourier coefficients by
proving an identity of local contributions place by place. In loc.\
cit.\ the identity for the archimedean place was proved, and it was
shown that the identity at a non-archimedean place of good
reduction is a consequence of results of Gross and Keating,
\cite{\grosskeating}, (for the algebraic-geometric side) and of
Kitaoka, \cite{\kitaoka}, (for the analytic side). It then remains
to consider the finite primes $p$ where the Shimura curve has bad
reduction. These are of two sorts: (i) the primes $p$ at which the
quaternion algebra defining the Shimura curve is split, but which
divide the level, and (ii) the primes $p$ at which the quaternion
algebra remains division.

In the present paper we
consider the case of a non-archimedean place $p$ of bad reduction of the second type,
and hence where $p$-adic uniformization in the sense of Cherednik-Drinfeld holds.
It turns out that the identity to be proved in this case can be
reduced to a purely local statement concerning the Drinfeld
$p$-adic upper half plane $\hat\Omega$ (the formal scheme version).
Therefore the bulk of this paper (sections 1-7) is concerned with
the local situation, and in this introduction we will concentrate
on the local aspects of the problem.

Let $W=W(\overline{\Bbb F}_p)$ be the ring of Witt vectors of
$\overline{\F}_p$. Also let $B$ be the division quaternion algebra over $\Q_p$
and let $O_B$ be its maximal order.

Recall that $\hat\Omega\times_{\roman{Spf}\,
\Z_p}{\roman{Spf}}\, W$ parametrizes pairs $(X,\varrho)$, where $X$
is a special formal (s.f.) $O_B$-module (of dimension 1 and height 2) over
a $W$-scheme $S$ on which $p$ is locally nilpotent and
$$\varrho:{\Bbb X}\times_{\roman{Spec}\, \overline{\F}_p}\overline
S\longrightarrow X\times_S\overline S\tag0.1$$
is an $O_B$-linear
quasi-isogeny of height 0. Here  ${\Bbb X}$ is a
fixed s.f.\ $O_B$-module over ${\roman{Spec}}\, \overline{\F}_p$, and
$\overline S = S\times_{\Spec\,\F_p}\Spec\,\overline{\F}_p$ is the special fibre of $S$. We note that
${\roman{End}}^0_{O_B}({\Bbb X})\simeq M_2(\Q_p)$.

We define cycles in $\hat\Omega\times_{\roman{Spf}\,\Z_p}{\roman{Spf}}\, W$ by
imposing additional endomorphisms as follows. Let
$$V=\{ j\in {\roman{End}}^0_{O_B}({\Bbb X});\ {\roman{tr}}^0(j)=0\},\tag0.2$$
equipped with the quadratic form $q=-{\roman{det}}$,
also given by $j^2=q(j)\cdot{\roman{id}}$. We call the elements of
$V$ {\it special endomorphisms.}\/ For any $j\in V$ with
$q(j)\in\Z_p\setminus \{ 0\}$ we define a {\it special cycle}\/
$Z(j)$ which is a closed formal subscheme of
$\hat\Omega\times_{\roman{Spf}\, \Z_p}{\roman{Spf}}\, W$. It is the
locus of pairs $(X, \varrho)$ such that $\varrho\circ j\circ
\varrho^{-1}$ extends to an isogeny of $X$.

We assume, from now on, that $p\ne 2$.

Our first task is to investigate the structure of a single special
cycle. For this we use two methods. The first is completely
elementary and uses the Bruhat-Tits building of $PGL_2(\Q_p)$. This
method is sufficient to give a fairly accurate picture of the point
set  of the special fibre of $Z(j)$. More precisely, let $\Cal B$
be the Bruhat-Tits building of $PGL_2(\Q_p)$, and recall that the
irreducible components of the special fiber of
$\hat\Omega\times_{\roman{Spf}\, \Z_p}{\roman{Spf}}\, W$ are
projective lines indexed by the vertices of $\Cal B$, i.e., by
homothety classes $[\Lambda]$ of $\Z_p$-lattices $\Lambda$ in
$\Q_p^2$. If $\Bbb P_{[\Lambda]}$ is the component corresponding to
$[\Lambda]$, then $Z(j)\cap \Bbb P_{[\Lambda]}\ne \emptyset$ if and
only if $j(\Lambda)\subset \Lambda$. It follows that the support of
$Z(j)$ is contained in the union of the $\Bbb P_{[\Lambda]}$'s for
$[\Lambda]$'s lying in the tube $\Cal T(j)$ of radius $\frac12
\roman{ord}\,q(j)$ around the fixed point set $\Cal B^j$ of $j$ in
$\Cal B$.
The second method is due to
Genestier,  \cite{\genestier}. His crucial observation is that $Z(j)$ may be
identified with the fixed point locus on
$\hat\Omega\times_{\roman{Spf}\,
\Z_p}{\roman{Spf}}\, W$ of the action of $j$ if
${\roman{ord}}_p(q(j))=0$ (resp.\ of ${\roman{id}} +j$ if
${\roman{ord}}_p(q(j))>0$). This observation allows one to write explicit
local equations for the cycle $Z(j)$ and thus yields a good
understanding of its local structure.  The combination of the two methods determines $Z(j)$ completely and
shows that, apart from degenerate cases, $Z(j)$ is purely
one-dimensional. Moreover, $Z(j)$ can contain (multiples of) irreducible components of the
special fibre and can even have embedded components! These latter
phenomena are in contrast to the case of good reduction. At the end
of the introduction there is a schematic picture of the various
possibilities of $Z(j)$.

We next turn to the calculation of the intersection product $(Z(j),Z(j'))$
of two special cycles, assuming that $j$ and $j'$ span a
2-dimensional non-degenerate quadratic $\Z_p$-submodule $\j$ of
$V$. In contrast with the case of good reduction, we have to deal here
with cases of excess intersection. We proceed in two steps. We
first prove that $(Z(j), Z(j'))$ only depends on the $\Z_p$-span
$\j$ of $j$ and $j'$. This is achieved by showing that the
Genestier equations globalize to give a resolution of the structure
sheaves of $Z(j)$ and $Z(j')$. This part of our paper is in the
spirit of the venerable theory of M\"obius transformations. In the
case of good reduction the analogue of this independence statement
is trivial, whereas at an archimedean place it was one of the
main and most difficult steps in the proof of the local identity,
\cite{\annals}. Since $p\ne2$, we may then assume that $j,j'$ diagonalize the
quadratic form on $\j$, i.e.\ that $jj'=-j' j$. In this case, the
calculation of the intersection number becomes a piece of recreational mathematics, involving
the various facts about the structure of special cycles mentioned above
and combinatorial arguments involving the tubes $\Cal T(j)$ in the building $\Cal B$. The
end result then is:
\proclaim{Theorem A} Let $j$ and $j'$ be special endomorphisms with
$q(j), q(j')\in\Z_p\setminus \{ 0\}$ such that their $\Z_p$-span
$\j
=
\Z_p j +
\Z_p j'$ is of rank $2$ and is nondegenerate for the quadratic form.
Let
$$T= \pmatrix q(j)&\frac12(j,j')\\\frac12(j',j)&q(j')\endpmatrix,$$
and suppose that $T$ is $GL_2(\Z_p)$-equivalent to
$\text{diag}(\e_1 p^\a,\e_2 p^\b)$, where $\varepsilon_1,
\varepsilon_2\in\Z_p^{\times}$, and $\a$ and $\b$ are integers with $0\le
\a\le\b$. Then $(Z(j),Z(j')) = e_p(T)$ depends only on the $GL_2(\Z_p)$-equivalence
class of $T$, and is given explicitly by:
$$e_p(T) = \a+\b+1 - \cases p^{\a/2}+2\, \frac{\scr p^{\a/2}-1}{\scr p-1}
&\text{if $\a$ is even and $\chi(\e_1)=-1$,}\\
\nass
(\b-\a+1)p^{\a/2} + 2\,\frac{p^{\a/2}-1}{p-1}&
\text{if $\a$ is even and $\chi(\e_1)=1$,}\\
\nass
 2\,\frac{p^{(\a+1)/2}-1}{p-1}&\text{if $\a$ is odd.}
\endcases
$$
\comment
$$\align
& \\
\nass
{}& =\a+\b+1 - \cases p^{\a/2}+2\, \frac{\scr p^{\a/2}-1}{\scr p-1}
&\text{if $\a$ is even and $\chi(\e_1)=-1$,}\\
\nass
(\b-\a+1)p^{\a/2} + 2\,\frac{p^{\a/2}-1}{p-1}&
\text{if $\a$ is even and $\chi(\e_1)=1$,}\\
\nass
 2\,\frac{p^{(\a+1)/2}-1}{p-1}&\text{if $\a$ is odd.}
\endcases
\endalign
$$
\endcomment
\endproclaim

The entity $e_p(T)$ appearing in this theorem may be related to
local representation densities of quadratic forms and their
derivatives. For simplicity, we continue to assume that $p\ne 2$.
Recall that for nonsingular symmetric matrices $S\in \Sym_m(\Z_p)$ and
$T\in \Sym_n(\Z_p)$, the classical representation density is defined by
\medskip\noindent
$$\a_p(S,T) = \lim_{t\rightarrow\infty} p^{-t n(2m-n-1)/2}\,
|\{\, x\in M_{m,n}(\Z/p^t\Z)\ ;\ S[x] - T \in p^t \Sym_m(\Z_p)\,
\}|.\tag0.3$$
Let
$$S=-\pmatrix 1&{}&{}\\{}&1&{}\\{}&{}&-1\endpmatrix\tag0.4$$
be the matrix for the determinant quadratic form on the space
$V(\Z_p)$ of special endomorphisms, i.e. on the lattice $\{x\in
M_2(\Z_p)\ ;\ tr(x)=0\}$, and let
$$S'=-\pmatrix \eta&{}&{}\\{}&p&{}\\{}&{}&-\eta p\endpmatrix,\tag0.5$$
be the matrix for the reduced norm quadratic form on the space
$$V'(\Z_p)=\{x\in O_B\ ;\  tr(x)=0\}.\tag0.6$$
Here $\eta\in \Z_p^\times$ with $\chi(\eta):=(\eta,p)_p=-1$, and $(a,b)_p$ is the
quadratic Hilbert symbol for $\Q_p$.
Also let
$$S''=-\pmatrix 1&{}&{}\\{}&p&{}\\{}&{}&-p\endpmatrix.\tag0.7$$
Then for a pair of special endomorphisms $j$ and $j'$ with
associated matrix $T$, as in Theorem~A above, we have
$\a_p(S',T)=0$. In this situation, it is possible, as in
\cite{\annals}, to define the {\it derivative} $\a_p'(S',T)$ of the
representation density (see (7.4)). For the unimodular quadratic
form $S$ and for any binary form $T$, the representation densities
and their derivatives can be calculated using the results of
Kitaoka, \cite{\kitaoka}, \cite{\annals}, and Proposition~7.1
below. For $S'$ and $S''$, the analogous information is provided by
the work of B. Myers, \cite{\myers}, and T. Yang, \cite{\yang}.

These representation densities are then related to the intersection number $e_p(T)$
as follows:
\proclaim{Theorem B} For a pair of special endomorphisms $j$ and $j'$ with associated matrix $T$, as in
Theorem~A,
$$e_p(T) =
-\frac{1}{p+1}\, \a_p'(S',-T) + \frac{2 p^2}{p+1}\, \a_p(S,-T)
+\frac{1}{2(p-1)} \a_p(S'',-T).$$
\endproclaim

This result illustrates again the rather remarkable connection
between arithmetic intersection numbers on certain moduli spaces on
the one hand, and the arithmetic theory of quadratic forms on the
other. For additional examples, see (at least)
\cite{\grosskeating}, \cite{\annals}, \cite{\krsiegel}, and
\cite{\krHB}.

It is instructive to compare the statements of Theorems~A and B
with the following reformulation of the result of Gross and
Keating, \cite{\grosskeating}. Changing notation slightly, we fix a
formal p-divisible group $\Bbb X$ of dimension $1$ and height $2$
over ${\roman{Spec}}\, \overline{\F}_p$, and let $\Cal M$ be the
moduli scheme where $\Cal M(S)$ is the set of pairs $(X,\rho)$ for
$\rho$ a quasi-isogeny of height $0$ as in (0.1). Then, there is a
(non-canonical) analogue of the Drinfeld isomorphism:
$$\Cal M\simeq \roman{Spf}\, W[[t]].$$
Fix an isomorphism $\End^0(\Bbb X)=B$, and let
$$V'=\{j\in B\ ;\ tr^0(j)=0\},$$
with (anisotropic) quadratic form defined by $j^2=q(j)\,\roman{id}$.
To a special endomorphism $j\in V'$ with $q(j)\ne0$, define a cycle $Z(j)\subset \Cal M$, as before.
\proclaim{Theorem} {\rm(Gross-Keating)}\  (i) If $q(j)\notin \Z_p$,
then $Z(j)=\emptyset$. Otherwise, $Z(j)$ is a divisor on $\Cal M$,
and is flat over $\roman{Spf}\,W$. \hfill\break (ii) Let $j$ and
$j'\in V'$ be special endomorphisms such that $j$ and $j'$ span a
2-dimensional non-degenerate subspace of $V'$. Let $T = q(j,j')$ be
the matrix of the quadratic form with respect to the basis $j$,
$j'$, as in Theorem~A. Then $Z(j)\cap Z(j')$ has support at the
origin in $\Cal M$ and the intersection multiplicity is
$$(Z(j),Z(j')) = \cases \sum_{i=0}^{(\a-2)/2} (\a+\b-4i)p^i +\frac12(\b-\a+1)p^{\a/2}&\text{ if $\a$ is even,}\\
\nass
\sum_{i=0}^{(\a-1)/2} (\a+\b-4i) p^i&\text{ if $\a$ is odd.}
\endcases
$$
In particular, this multiplicity depends only on the
$GL_2(\Z_p)$-equivalence class of $T$.
\endproclaim

By Kitaoka's formula, Corollary~8.5 of \cite{\annals}, we then have
the (simpler) analogue of Theorem~B in this case:
$$(Z(j),Z(j'))  = - \frac{p^2}{p^2-1} \a_p'(S,T)\ \ .$$

Two features of the theory developed here are worth pointing out.
Recall that the proofs of Gross and Keating in the case of good
reduction make heavy use of the theory of (formal) complex
multiplication which connects this case, via Gross's theory of
quasi-canonical liftings, \cite{\gross}, with Kronecker's
Jugendtraum. The first remark is that in our case this connection
does not appear (explicitly! - it is of course hidden to some
degree in Drinfeld's representability theorem). The second remark
is that it is the global nature of the Cherednik-Drinfeld
uniformization which allowed us here to prove the independence
statement on the intersection numbers. In other, higher-dimensional
cases \cite{\krsiegel}, \cite{\krHB} when the analogous
independence property is problematical, {\it global} uniformization
of the special fiber is not available. In these cases these
problems remain a challenge.

In the last two sections, we draw the global consequences of the
local results of sections 1--7  just described and obtain an extension of the
results of \cite{\annals}.

For an indefinite quaternion algebra $B$ over $\Q$, let
$H=B^\times$, and let $V=\{x\in B\ ;\  tr^0(x)=0\}$. Fix a prime
$p$ ( $p\ne 2$) which ramifies in $B$ and let $K= K_pK^p$ be a
compact open subgroup of $H(\A_f)$ such that $K_p =
O_{B_p}^\times$, where $O_{B_p}$ denotes the maximal order in
$B_p$. Associated to this data is a model $\Cal A_K$ over
$\Z_{(p)}$ of the Shimura curve $A_K$ over $\Q$ attached to $B$ and
$K$. It is the moduli space of certain abelian surfaces with
$O_B$-action and $K^p$-level structure,   For each pair $(t,\o)$
with $t\in
\Z_{(p)}$, $t<0$ and $\o\subset V(\A_f^p)$ a $K^p$-invariant
compact open subset, there is a special cycle
$$\Cal C(t,\o)\lra \Cal A_K\tag0.8$$
defined, as in \cite{\annals}, by imposing an additional special
endomorphism. The generic fiber of $\Cal C(t,\o)$ is a
$\Q$-rational $0$-cycle on the Shimura curve $A_K$. Given a pair of
special cycles $\Cal C_1 =\Cal C(t_1,\o_1)$ and $\Cal C_1=\Cal
C(t_2,\o_2)$, we form their intersection
$$\Cal C= \Cal C_1\times_{\Cal A_K}\Cal C_2.\tag0.9$$
This scheme has a decomposition
$$\Cal C=\coprod_T\Cal C_T,\tag0.10$$
where each $\Cal C_T$ is a union of connected components of $\Cal
C$, and where $T\in \roman{Sym}_2(\Z_{(p)})$ runs over negative
semi-definite matrices with diagonal entries $t_1$ and $t_2$. For
$\det(T)\ne0$, the image of $\Cal C_T$ lies in the special fiber of
$\Cal A_K$. Utilizing the p-adic uniformization of $\Cal A_K$, and
the results on intersections of special cycles in Drinfeld space,
we obtain the following statement.
\proclaim{Theorem C} Assume that $t_1t_2$ is not a square in $\Q^\times$. Then the
special cycles $\Cal C(t_1,\o_1)$ and $\Cal C(t_2,\o_2)$ do not meet in the
generic fiber and their intersection number is given by
$$(\Cal C(t_1,\o_1),\Cal C(t_2,\o_2)) = 2 \sum_{T} e_p(T)\, \vol(K^p)^{-1}\, I_T(\ph_1^p\tt\ph_2^p),$$
where $e_p(T)$ is as in Theorem~A and $I_T(\ph_1^p\tt\ph_2^p)$ is a certain orbital integral
associated to the data $\o_1$, $\o_2$, and $T$. The summation runs over the
same range as in (0.10).
\endproclaim

Finally, this result, combined with Theorem~B above yields the
connection between the $p$ part of the height pairing
of the cycles $\Cal C_1$ and $\Cal C_2$ and certain Fourier coefficients
of the derivative of a metaplectic Eisenstein series (Theorem~9.1).
This is analogous to Theorem~14.11 of \cite{\annals} and Theorem~9.2 of \cite{\krsiegel}.
We content ourselves here with
pointing out one essential difference with the case of good
reduction treated in \cite{\annals}. In that case, the choice of
the local component $\Phi_p$
of the function occurring in the Eisenstein series is {\it
canonical}\/ and it is in fact {\it standard}\/ in the sense that
its restriction to the maximal compact subgroup of $Mp_{2,\Q_p}$ is
independent of the complex parameter $s$. In our case, the choice of
$\Phi_p$, given in Corollary~7.4, is no longer canonical; rather, we are able to single
out a whole class of functions for which the main identity holds,
and these functions are definitely not standard. It seems quite
likely that among those functions there are preferred choices,
namely those that match up the Eisenstein series with $L$-functions
via the doubling method, but this will not concern us here.

Let us give a brief description of the contents of the various
sections. Section 1 contains recollections about the Drinfeld
moduli space and introduces our terminology (esp.\ {\it ordinary
special}\/ and {\it superspecial}\/ points). In section 2 we use
the building to determine the point set of the special fibre of a
special cycle (resp.\ the intersection of two of them). In section 3
we use the Genestier equations to determine the local structure of
a special cycle (multiplicities of vertical components, occurrence of
embedded components, etc.). In section 4 we explain the kind of intersection
theory we are using. In section 5 we construct a global resolution of the
structure sheaves of the special cycles and prove the above-mentioned
invariance property of the intersection numbers. Section 6
contains the calculation of the intersection number in the diagonalized case.
In section 7 we
review the results of Myers, \cite{\myers}, and Yang, \cite{\yang}, and establish the connection
with representation densities. In section 8 we pose the global
intersection problem for cycles on models of Shimura curves over $\Z_p$
and relate it to the local theory. The final section 9
gives the relation, extending that of \cite{\annals}, between the intersection numbers and
special values of the derivatives of Fourier coefficients of Eisenstein series.

In conclusion we wish to point out again that one of the main
ingredients of this paper is due to Genestier \cite{\genestier}; regrettably,
he did not pursue the further implications of his idea. We also
thank T.\ Yang for communicating to us his results on
representation densities at an early stage of this project. This
work was begun at the University of Cologne in August 1997 and
continued at the University of Maryland in March 1998. We thank
both institutions for their hospitality and the NSF and the DFG for
their support.
\vskip .5in

\centerline{\bf Contents}
\roster
\item"{1.}" The Drinfeld moduli space and the p-adic upper half plane
\item"{2.}" Special cycles and their support in the special fiber
\item"{3.}" Local equations for special cycles
\item"{4.}" Intersection calculus of special cycles
\item"{5.}" An invariance property of intersection numbers
\item"{6.}" Computation of intersection numbers
\item"{7.}" Intersection numbers and representation densities
\item"{8.}" Intersection numbers on Shimura curves
\item"{9.}" Intersection numbers and Fourier coefficients
\endroster

\vfill
\eject



\input drinfigmittext.tex

\vfill\eject

\subheading{Notation}
\par\noindent
The following notation will be used in the local part of this
paper (sections 1-7):
$$\align
k &\ \
\text{algebraically closed field of characteristic $p>2$}\\
W=W(k) &\ \
\text{the ring of Witt vectors of $k$ with its Frobenius automorphism
$\sigma$.}\\ W_{\Q} &
\ \ \text{the fraction field of $W$.}\\
\Z_{p^2}=W(\F_{p^2})&\ \
\text{where $\F_{p^2}= \{ x\in k;\ x^{p^2}=1\}$ .}\\
B &\ \
\text{a quaternion division algebra over $\Q_p$.}\\
O_B &\ \
\text{the ring of integers which we identify with}\\
&\ \ O_B= \Z_{p^2}[\Pi] /\Pi^2 -p,\ \Pi a= a^{\sigma}\Pi\ \ (a\in
\Z_{p^2}).\\ {\Cal B} &\ \ \text{the Bruhat-Tits building of
$PGL_2(\Q_p)$.}
\\ {\chi} &\ \ \text{the quadratic residue character on
$\Z_p^{\times}$ resp.\ ${\Bbb F}_p^{\times}$.}
\endalign$$

\subheading{\Sec 1. The Drinfeld moduli space and
the $p$-adic upper half plane}
\par\noindent
In this section we recall some facts from Drinfeld's paper
\cite{\drinfeld}, cf.\ also \cite{\boutotcarayol}, \cite{\genestierbook}, \cite{\rapoportzink}. A {\it
special formal} (s.f.) $O_B$-{\it module} over a $W$-scheme $S$ is
a $p$-divisible formal group $X$ of height 4, with an $O_B$-action
$\iota: O_B\to {\roman{End}}_S(X)$ such that the induced
$\Z_{p^2}\otimes {\Cal O}_S$-module ${\roman{Lie}}\, X$ is, locally
on $S$, free of rank 1. Let us fix once and for all a s.f.\
$O_B$-module ${\Bbb X}$ over ${\roman{Spec}}\, k$. It is unique up
to $O_B$-linear isogeny, and ${\roman{End}}^0_{O_B}({\Bbb X})\simeq
M_2(\Q_p)$. Let us consider the following functor ${\Cal M}$ on the
category Nilp of $W$-schemes $S$ such that $p$ is locally nilpotent
in ${\Cal O}_S$. The $S$-valued points of ${\Cal M}$ are the
isomorphism classes of pairs $(X, \varrho)$ consisting of a s.f.\
$O_B$-module $X$ over $S$ and a quasi-isogeny of height zero,
$$\varrho: {\Bbb X}\times_{\roman{Spec}\, k}
\overline S \to X\times_S \overline S\ \ .$$
Here $\overline S =S\times_{\roman{Spec}\, W} {\roman{Spec}}\, k$.
We consider
$$G(\Q_p)^0 =\{ g\in GL_2(\Q_p);\ {\roman{ord}}\, {\roman{det}}\, g=0\}
\ .\tag1.1$$
Then, after choosing an identification ${\roman{End}}^0_{O_B}({\Bbb
X})=M_2(\Q_p)$, the group $G(\Q_p)^0$ acts to the left on ${\Cal
M}$, via
$$g: (X,\varrho)\mapsto (X, \varrho\circ g^{-1})\ \ .$$
According to Drinfeld, this functor ${\Cal M}$ is representable by
$\hat\Omega\times_{\roman{Spf}\, \Z_p} \roman{Spf}\, W$. Here
$\hat\Omega =\hat\Omega^2_{\Q_p}$ is the formal model of the
$p$-adic upper half plane associated to the local field $\Q_p$ that
was introduced by Deligne (comp.\ \cite{\boutotcarayol}, chap.\ 1). The
isomorphism
$${\Cal M}\to \hat\Omega\times_{\roman{Spf}\,\Z_p}\roman{Spf}\, W\tag1.2$$
is equivariant for the action of $G(\Q_p)^0$, for a suitable
identification $\roman{End}^0_{O_B}({\Bbb X})=M_2(\Q_p)$. The group
$G(\Q_p)^0$ acts on the RHS via the natural action of $PGL_2(\Q_p)$
on $\hat\Omega$.

We need to describe some features of $\hat\Omega$, comp.\
\cite{\boutotcarayol}, chap.\ 1. We denote by ${\Cal B}={\Cal B}(PGL_2(\Q_p))$
the Bruhat-Tits building of $PGL_2(\Q_p)$. The formal scheme
$\hat\Omega$ is obtained by glueing open formal subschemes
$\hat\Omega_{\Delta}$ where $\Delta$ runs over the simplices of
${\Cal B}$, and
$$\hat\Omega_{\Delta}\cap \hat\Omega_{\Delta'} =
\cases \hat\Omega_{\Delta\cap\Delta'}
&
\text{ if $\Delta\cap\Delta'$ is a simplex}\\
\emptyset&
\text{ if $\Delta\cap\Delta'=\emptyset$.}
\endcases\tag1.3$$
For the action of $PGL_2(\Q_p)$ on $\hat\Omega$ we have
$$g\hat\Omega_{\Delta} =\hat\Omega_{g\Delta}\ \ .\tag1.4$$

We now describe the open charts in detail. Because of (1.4) it will
suffice to do this for $\Delta =$ the standard vertex and for
$\Delta =$ the standard edge.

Let $\Delta= [\Lambda_0]$ be the homothety class of the standard
lattice
$$\Lambda_0= [e_1, e_2]\ \ .\tag1.5$$
Here $[e_1, e_2]$ denotes the $\Z_p$-span of the standard basis of
$\Q^2_p$. Then
$$\hat\Omega_{[\Lambda_0]} =({\Bbb P}(\Lambda_0)-
{\Bbb P}(\Lambda_0)({\Bbb F}_p))^{\wedge}\ \ .\tag1.6$$ Here ${\Bbb
P}(\Lambda_0)\simeq {\Bbb P}^1_{\Z_p}$ denotes the relative
projective line over $\Z_p$ associated to $\Lambda_0$ and the
``hat'' indicates the completion along the special fibre. If we use
$e_1, e_2$ to identify ${\Bbb P}(\Lambda_0)$ with ${\Bbb
P}^1_{\Z_p}$, we have
$$\hat\Omega_{[\Lambda_0]} ={\roman{Spf}}\,
\Z_p[T, (T^p-T)^{-1}]^{\wedge}\ \ ,\tag1.7$$
where $T=X_0/X_1$ in terms of the canonical coordinates on ${\Bbb
P}^1_{\Z_p}$. The subgroup $GL_2(\Z_p)$ preserves $\Lambda_0$ and
hence acts on $\hat\Omega_{[\Lambda_0]}$. The action of $g$ on
${\Bbb P}(\Lambda_0)$ is induced by the automorphism
$g:\Lambda_0\to \Lambda_0$. For the left action $f\mapsto
g_*^{-1}(f)$ of $GL_2(\Z_p)$ on the ring of holomorphic functions
on $\hat{\Omega}_{[\Lambda_0]}$ we therefore have
$$g_*^{-1}:(X_0, X_1)\mapsto ({\roman{det}} (g)^{-1}
(dX_0-cX_1),\ {\roman{det}} (g)^{-1}(-bX_0+ aX_1)),\ g=\pmatrix
a&b\\ c&d\\\endpmatrix\ \ .$$ In terms of the coordinate $T$ the
action $g_*^{-1}$ is therefore given by
$$g_*^{-1}:T\mapsto \frac{dT-c}{-bT+a}\ \ ,\ \
g=\pmatrix a&b\\ c&d\\\endpmatrix\tag1.8$$ (homography associated
to ${}^tg^{-1}$).

Next, let $\Delta_0=([\Lambda_0], [\Lambda_1])$ be the edge
corresponding to
$$\Lambda_0=[e_1, e_2]\ \ ,\ \ \Lambda_1= [pe_1, e_2]\ \ .\tag1.9$$
In this case we have an identification
$$\hat\Omega_{\Delta_0}= {\roman{Spf}}\, \Z_p[T_0, T_1,
(1- T_0^{p-1})^{-1}, (1- T_1^{p-1})^{-1}]^{\wedge} / (T_0T_1-p)\
\ .\tag1.10$$ The action of the Iwahori subgroup
$$\pmatrix \Z_p^{\times} & p\Z_p\\ \Z_p & \Z_p^{\times}\\
\endpmatrix\tag1.11$$
is given by
$$g_*^{-1}:T_0\mapsto \frac{dT_0 -pc}{-b_0T_0+a}\ \ ,\ \ T_1\mapsto
\frac{aT_1-pb_0}{-cT_1+d}\ \ ,\tag1.12$$
for $$g=\pmatrix a&b\\ c&d\\\endpmatrix =\pmatrix a&pb_0\\
c&d\\\endpmatrix\ \ .\tag1.13$$ The open immersion
$\hat\Omega_{[\Lambda_0]} \to \hat\Omega_{\Delta_0}$ is induced by
the open immersion
$${\roman{Spf}}\, \Z_p[T, T^{-1}]^{\wedge} \to {\roman{Spf}}
(\Z_p[T_0, T_1] / T_0 T_1-p)^{\wedge}\tag1.14$$ induced by
$$T_1\mapsto T^{-1}\ \ ,\ \ T_0\mapsto p\cdot T\ \ .\tag1.15$$
It is easy to check that this morphism is indeed equivariant for
the action of the Iwahori subgroup (1.11).

The special fibre of $\hat\Omega$ is a union of projective lines
parametrized by the vertices in ${\Cal B}$. More precisely, ${\Cal
B}$ can be identified with the dual graph of the special fibre,
compatibly with the action of $PGL_2(\Q_p)$. Let $[\Lambda]$ be a
vertex with corresponding projective line ${\Bbb P}_{[\Lambda]}$.
Then ${\Bbb P}_{[\Lambda]}$ may be identified with
$${\Bbb P}_{[\Lambda]} ={\Bbb P}(\Lambda /p\Lambda)\ \ ,\tag1.16$$
where $\Lambda$ is any lattice in the homothety class $[\Lambda]$,
and
$$\hat\Omega_{[\Lambda]} \times_{\roman{Spf}\, \Z_p}
\roman{Spec}\, {\Bbb F}_p = {\Bbb P}_{[\Lambda]}
 -{\Bbb P}_{[\Lambda]} ({\Bbb F}_p)\ \ .\tag1.17$$ If $\Delta=
([\Lambda], [\Lambda'])$ is an edge, we denote the corresponding
point in the special fibre by $pt_{\Delta}$,
$$pt_{([\Lambda],[\Lambda'])} = {\Bbb P}_{[\Lambda]}\cap
{\Bbb P}_{[\Lambda']} =\hat\Omega_{\Delta} ({\Bbb F}_p)\ \
.\tag1.18$$

We shall also need to know how the Drinfeld isomorphism (1.2) looks
on the set of closed points. Let $(X, \varrho)\in {\Cal M}(k)$ and
let $(M, F,V)$ be the covariant Dieudonn\'e module of $X$. It is a
free $W$-module of rank 4. From the action of $\Z_{p^2}$ on $X$ we
obtain a $\Z /2$-grading,
$$M=M_0\oplus M_1\ \ ,\tag1.19$$
and $F,V$ and $\iota(\Pi)$ are all homogeneous of degree 1. Since
${\roman{Lie}}\, X =M/VM$ and by the condition on $X$ to be
special, we have inclusions with index 1 of free $W$-modules of
rank 2,
$$pM_0\mathop{\subset}\limits_{\ne} VM_1\mathop{\subset}\limits_{\ne} M_0
\ \ ,\ \ pM_1\mathop{\subset}\limits_{\ne} VM_0
\mathop{\subset}\limits_{\ne} M_1\ \ .\tag1.20$$
An index $i\in \Z /2$ is called {\it critical,} if
$$VM_i=\Pi M_i\ \ .\tag1.21$$
Since ${\roman{Lie}}(\iota(\Pi))^2 =0$ and $\dim_k M_0/VM_1=\dim_k
M_1/VM_0 =1$, there always exists $i$ with $\Pi M_i\subset VM_i$.
Since both modules have index 1 in $M_{i+1}$, it follows that $i$
is critical. Hence there always exists a critical index.
\hfill\break
If $i$ is a critical index, $V^{-1}\Pi$ is a $\sigma$-linear
automorphism of $M_i$. If we set $\eta_i= M_i^{V^{-1}\Pi}$, then
$\eta_i$ is a $\Z_p$-module with
$$M_i= \eta_i\otimes_{\Z_p}W\ \ .$$
Recall our fixed s.f.\ $O_B$-module ${\Bbb X}$. We may choose
${\Bbb X}$ in its isogeny class so that 0 and 1 are critical. Let
$$M^0=M_0^0\oplus M_1^0\ \ ,\ \ \eta^0_i =M_i^{0^{V^{-1}\Pi}}$$
be the Dieudonn\'e module of ${\Bbb X}$ and let us fix an
isomorphism
$$U:= \eta_0^0\otimes_{\Z_p}\Q_p\simeq \Q_p^2\ \ .\tag1.22$$
We can now describe the Drinfeld isomorphism on ${\Cal M}(k)$. Let
$(X,\varrho)\in {\Cal M}(k)$ and let $M$ be the Dieudonn\'e module
of $X$. If 0 is critical we have the $\Z_p$-module $\eta_0$ and
$\varrho$ defines an isomorphism
$$\eta_0\otimes \Q_p\to U =\Q_p^2\ \ .\tag1.23$$
Let $\Lambda_0$ be the image of $\eta_0$ under (1.23). Then
$\Lambda_0$ is a lattice in $\Q_p^2$ which (since the height of
$\varrho$ is zero) is of the same volume as $\Z_p^2$. The point
$(X, \varrho)$ is then mapped to the point of ${\Bbb
P}_{[\Lambda_0]} (k)$ corresponding to the line
$$\ell = VM_1/pM_0\subset M_0/pM_0 =\Lambda_0\otimes_{\Z_p}k\ \ .\tag1.24$$
If 1 is critical, $\varrho$ defines an isomorphism
$$\eta_1\otimes \Q_p\to \eta_1^0 \otimes \Q_p\ \ .\tag1.25$$
The action of $\Pi$ identifies $\eta_1^0\otimes \Q_p$ with $U$.
Hence we obtain a lattice $\Lambda_1$ in $\Q_p^2$, of the same
volume as $p\Z_p\oplus \Z_p$. The point $(X,\varrho)$ is mapped to
the point of ${\Bbb P}_{[\Lambda_1]}(k)$ corresponding to the line
$$\ell = VM_0/pM_1\subset M_1/pM_1 =\Lambda_1\otimes_{\Z_p}k\ \ .\tag1.26$$
If both 0 and 1 are critical, both procedures above are applicable
and the lattices $\Lambda_1\subset \Lambda_0$ define an edge
$\Delta =([\Lambda_0], [\Lambda_1])$ in ${\Cal B}$. In this case
the line $VM_1/pM_0$ in ${\Bbb P}_{[\Lambda_0]}$ and $VM_0/pM_1$ in
${\Bbb P}_{[\Lambda_1]}$ define the same point of $\hat\Omega(k)$,
namely $pt_{\Delta}$. We thus see that the set of irreducible
components of $\hat\Omega\otimes k$ can be partitioned into two
subsets: on the irreducible components corresponding to even
lattice classes the index 0 is critical and on the irreducible
components corresponding to odd lattice classes the index 1 is
critical.
\hfill\break
We shall use the following terminology.
\medskip\noindent
{\bf Definition 1.1.} A point $(X,\varrho)\in {\Cal M}(k)$ is
called
\hfill\break
{\it superspecial}\/ if both indices 0,1 are critical
\hfill\break
{\it ordinary}\/ if only one index is critical
\hfill\break
{\it special}\/ if $V^2M=pM$.

Obviously, a superspecial point is special since in this case
$VM=\Pi M$. Assume that $(X,\varrho)$ is ordinary and let e.g.\ 0
be the unique critical index,
$$VM_0=\Pi M_0\ \ .\tag1.27$$
Since $V^2M_0= pM_0$, we see that $(X, \varrho)$ is special if and
only if
$$V^2M_1= pM_1,\ \ i.e.\ \ (V^{-1}\Pi)^2VM_1=VM_1\ \ .\tag1.28$$
Recalling that $VM_1/pM_0\subset M_0/pM_0$ is the line in ${\Bbb
P}_{[\Lambda_0]}$ associated to $M$, we obtain the following
characterization.
\proclaim{Proposition 1.2}
Let ${\Bbb P}_{[\Lambda]}\subset\hat\Omega \otimes_{\Z_p}k$ be the
irreducible component corresponding to the vertex $[\Lambda]$ in
${\Cal B}$. Then a point $x\in {\Bbb P}_{[\Lambda]}(k)$ is
\hfill\break
superspecial iff $x\in {\Bbb P}_{[\Lambda]}({\Bbb F}_p)$
\hfill\break
special iff $x\in {\Bbb P}_{[\Lambda]} ({\Bbb F}_{p^2})$.
\endproclaim
\par\noindent
{\bf Remark 1.3.} Let $(X, \varrho)\in {\Cal M}(k)$ be special.
Then we may find a $W$-basis of its Dieudonn\'e module with
$$e_3= Ve_1\ ,\ e_4=Ve_2\ ,\ pe_1=Ve_3\ ,\ pe_2=Ve_4\ \ .\tag1.29$$
Hence $X$ is isomorphic to the product of the $p$-divisible group
of a supersingular elliptic curve with itself. The converse also
holds.

\subheading{\Sec 2. Special cycles and their support in the special fibre}
\par\noindent
Recall our fixed s.f.\ $O_B$-module ${\Bbb X}$ over $k$. We
introduce
$$V=\{ j\in {\roman{End}}^0_{O_B}({\Bbb X});
\ {\roman{tr}}^0(j)=0\}\ \ .\tag2.1$$
We will always identify $j\in {\roman{End}}^0_{O_B}({\Bbb X})$ with
its image under the identifications (cf.\ (1.22)),
$${\roman{End}}^0_{O_B}({\Bbb X})= {\roman{End}}(U)= M_2(\Q_p)\ \ .
\tag2.2$$
Then $V$ is equipped with a quadratic form given by squaring,
$$j^2=q(j)\cdot {\roman{id}},\ \ q(j)\in \Q_p\ \ .\tag2.3$$
We note that
$$q(j)=-{\roman{det}}(j)\ \ .\tag2.4$$
The elements of $V$ will be called {\it special endomorphisms.}
\medskip\noindent
{\bf Definition 2.1.} Let $j\in V$ be a special endomorphism with
$q(j)\ne 0$. The {\it special cycle}\/ $Z(j)$ {\it associated to}\/
$j$ is the closed formal subscheme of ${\Cal M}$ consisting of all
points $(X, \varrho)$ such that $\varrho\circ j\circ
\varrho^{-1}$ lifts to an endomorphism of $X$.

Let $(X, \varrho)\in {\Cal M}(S)$. Then $\lambda = \varrho\circ
j\circ \varrho^{-1}$ is a quasi-isogeny of $X\times_S\overline S$.
By the rigidity property of quasi-isogenies (\cite{\drinfeld}) we may also
consider $\lambda$ as a quasi-isogeny of $X$. The locus in $S$
where $\lambda$ is an isogeny is a closed subscheme $S'$ of $S$
(\cite{\rapoportzink}, prop.\ 2.9). The closed formal subscheme $Z(j)$ of
${\Cal M}$ is characterized by
$$S'=Z(j)\times_{\Cal M}S\ \ .\tag2.5$$

In this section we will study the point set $Z(j)(k)$ of special
cycles and their intersection properties.
\proclaim{Proposition 2.1}
Let $[\Lambda]$ be a vertex in ${\Cal B}$ and let $\Lambda\subset
U=\Q_p^2$ be a lattice in its homothety class. Then
$${\Bbb P}_{[\Lambda]} \cap Z(j)\ne\emptyset\ \ \hbox{iff}\ j(\Lambda)
\subset\Lambda\ \ .$$
In particular, if $Z(j)\ne\emptyset$ then $q(j)\in \Z_p$.
\endproclaim
\par\noindent
We will see in Corollary 2.5 that conversely, if $q(j)\in\Z_p$,
then $Z(j)\neq \emptyset$.
\demo{Proof}
Suppose e.g.\ that $[\Lambda]$ is even and let $x\in {\Bbb
P}_{[\Lambda]} (k)$. Then $x$ corresponds under the Drinfeld
isomorphism to $(X, \varrho)$ where $0$ is a critical index and
with Dieudonn\'e module $M=M_0\oplus M_1$, where
$$M_0= \Lambda \otimes_{\Z_p}W\ \ ,\ \ VM_1/pM_0=\ell_x\subset \Lambda
\otimes k\ \ .\tag2.6$$
Then $x$ lies in $Z(j)$ iff $j\otimes_{\Q_p}{\roman{id}}_{W_{\Q}}$
preserves the $W$-lattices $M_0$ and $VM_1$ in
$U\otimes_{\Q_p}W_{\Q}$. This proves one implication. Assume now
that $j(\Lambda)\subset
\Lambda$, i.e.\ $j(M_0)\subset M_0$. Since $k$ is algebraically closed
there exists a line $\ell\subset \Lambda\otimes_{\Z_p}k$ stable
under the endomorphism induced by $j$. The corresponding point
$x\in {\Bbb P}_{[\Lambda]}(k)$ is associated to $(X,\varrho)$ with
$M=M_0\oplus M_1$ where
$$j(VM_1/ pM_0)\subset VM_1/pM_0,\ \text{i.e.,}\ j(VM_1)\subset VM_1\ \ .\tag2.7$$
Hence $x\in Z(j)$.
\hfill\break
The last assertion follows by taking a $\Z_p$-basis of a lattice
$\Lambda$ with ${\Bbb P}_{[\Lambda]} \cap Z(j)\ne\emptyset$ and
using it to calculate ${\roman{det}}(j)=-q(j)$.
\qed
\enddemo
{\it From now on we will always assume that} $q(j)\in\Z_p
- \{ 0\}$. We next do a local analysis of $Z(j)$ along a line ${\Bbb
P}_{[\Lambda]}$ which intersects $Z(j)$. Let us fix a lattice
$\Lambda\subset U$ with $j(\Lambda)\subset \Lambda$. We use the
notation
$${\roman{red}}_{\Lambda}(j)\in {\roman{End}}(\Lambda /p\Lambda)
\tag2.8$$
for the induced endomorphism.
\proclaim{Proposition 2.2}
With the notation introduced above, there are the following
possibilities.
\hfill\break
(i) ${\roman{rk}}({\roman{red}}_{\Lambda}(j))=2$. Then
$q(j)=\varepsilon\in \Z_p^{\times}$ and ${\roman{red}}(j)$
preserves precisely two lines in ${\Bbb P}_{[\Lambda]}(k)$. The
corresponding points are both superspecial, if
$\chi(\varepsilon)=1$ and both ordinary special, if
$\chi(\varepsilon)=-1$.
\hfill\break
(ii) ${\roman{rk}}({\roman{red}}_{\Lambda}(j))=1$. Then
${\roman{ord}}\, q(j)\ge 1$ and ${\roman{red}}_{\Lambda}(j)$ is a
nilpotent endomorphism. The line
$$\ell = {\roman{Ker}}\, {\roman{red}}_{\Lambda}(j)= {\roman{Im}}
\, {\roman{red}}_{\Lambda}(j)$$
is the unique line stable under ${\roman{red}}_{\Lambda}(j)$. The
corresponding point of ${\Bbb P}_{\Lambda}$ is superspecial.
\hfill\break
(iii) ${\roman{red}}_{\Lambda}(j)=0$. Then ${\roman{ord}}\, q(j)\ge
2$ and all lines $\ell\in {\Bbb P}_{[\Lambda]} (k)$ are stable
under ${\roman{red}}_{\Lambda}(j)$, i.e.\ ${\Bbb
P}_{[\Lambda]}\subset Z(j)$.
\endproclaim
\demo{Proof}
Let us prove (i). Since ${\roman{red}}_{\Lambda}(j)$ is a traceless
automorphism it has two distinct eigen lines. The characteristic
polynomial of ${\roman{red}}_{\Lambda}(j)$ has the form
$$X^2-q(j)\ \ {\roman{mod}}\, p\ \ ,$$
since $q(j)= -{\roman{det}}(j)$, and the determinant may be
calculated using a basis of $\Lambda$. The assertion follows. The
other assertions are proved in a similar way.
\qed\enddemo
\proclaim{Corollary 2.3}
$Z(j)(k)$ consists only of isolated points if and only if
${\roman{ord}}\, q(j)\le 1$.
\endproclaim
\demo{Proof}
If ${\Bbb P}_{[\Lambda]}\subset Z(j)$ then
${\roman{red}}_{\Lambda}(j) = 0$ and ${\roman{ord}}\, q(j)\ge 2$.
Conversely, let ${\roman{ord}}\, q(j)\ge 2$. By Corollary 2.5 below
there exists $\Lambda$ with ${\Bbb P}_{[\Lambda]}\cap
Z(j)\ne\emptyset$. Then ${\roman{red}}_{\Lambda}(j)$ cannot have
full rank because otherwise ${\roman{ord}}\, q(j)=0$. Therefore, if
${\Bbb P}_{[\Lambda]}\not\subset Z(j)$ we are in case (ii) of the
previous proposition. Let $\Lambda'$ be the lattice neighbour of
$\Lambda$ corresponding to the intersection point of ${\Bbb
P}_{[\Lambda]}$ with $Z(j)$, i.e.\ to the line ${\roman{Im}}\,
{\roman{red}}_{\Lambda}(j)$ in $\Lambda /p\Lambda$. Then
$$\Lambda' =j(\Lambda)+ p\Lambda\ \ .\tag2.9$$
But then
$$\align
j(\Lambda') &
=j^2(\Lambda)+ pj(\Lambda)\\
&
=p^{\alpha}\Lambda+pj(\Lambda),\ \ \alpha\ge 2\\
&
=p(j(\Lambda)+ p^{\alpha-1}\Lambda)\\
&
\subset p\Lambda'\ \ .
\endalign$$
Hence ${\roman{red}}_{\Lambda'}(j)=0$ and ${\Bbb
P}_{[\Lambda']}\subset Z(j)$.
\qed\enddemo
We next will get a global overview of the lattices $\Lambda$ which
satisfy the criterion of Proposition 2.1, i.e., of the set
$${\Cal T}(j)= \{ [\Lambda]\in {\Cal B};\ \ {\Bbb P}_{[\Lambda]}
\cap Z(j)\ne\emptyset\}\ \ .\tag2.10$$
\proclaim{Lemma 2.4}
Let $j\in GL_2(\Q_p)$. Let $\Lambda\subset
\Q_p^2$ be a lattice. The following conditions are equivalent:
\hfill\break
(i) $j(\Lambda)\subset\Lambda$
\hfill\break
(ii) Let $[\Lambda]\in {\Cal B}$ be the vertex corresponding to
$\Lambda$. Then
$$d([\Lambda], [j(\Lambda)])\le {\roman{ord}}\, {\roman{det}}\, j\ \ .$$
Here $d$ denotes the distance in the building.
\hfill\break
If ${\roman{tr}}(j)=0$, these conditions are also equivalent to
\hfill\break
(iii) $d([\Lambda], {\Cal B}^j)\leq \frac{1}{2}\cdot
{\roman{ord}}\, {\roman{det}}\, j$.
\hfill\break
Here ${\Cal B}^j$ denotes the fixed point set of $j$ in ${\Cal B}$.
\endproclaim
\demo{Proof}
Put $\alpha= {\roman{ord}}\, {\roman{det}}\, j$. Let
$(e,f)\in\Z^2$, $e\ge f$, be the elementary divisors for the
lattice pair $\Lambda, j(\Lambda)$. Then $e+f=\alpha$ and
$d([\Lambda], [j(\Lambda)])=e-f$. On the other hand
$$j(\Lambda)\subset\Lambda \Longleftrightarrow f\ge 0\Longleftrightarrow
e-f\le \alpha\ \ .\tag2.11$$ This shows the equivalence of (i) and
(ii). For (iii) we note that if ${\roman{tr}}(j)=0$, then $j$
induces an involution of ${\Cal B}$. The unique geodesic from
$[\Lambda]$ to $[j(\Lambda)]$ consists of the geodesic from
$[\Lambda]$ to ${\Cal B}^j$ and its image under $j$ which is the
geodesic from $[j(\Lambda)]$ to ${\Cal B}^j$.
\qed\enddemo
\par\noindent
{\it Remark:} This observation appears already in
\cite{\kottwitztwo} where it is attributed to Tate.
\proclaim{Corollary 2.5}
Let $j\in V$ with $q(j)=\varepsilon\cdot p^{\alpha}$, where
$\varepsilon\in\Z_p^{\times}$, $\alpha\ge 0$. The set ${\Cal T}(j)$
(cf.\ (2.10)) is as follows:
\hfill\break
(i) If $\alpha$ is even and  $\chi(\varepsilon)=-1$, then ${\Cal
T}(j)$ is a ball of radius $\alpha/2$ around the unique vertex
$[\Lambda_0]$ fixed by $j$.
\hfill\break
(ii) If $\alpha$ is even and $\chi(\varepsilon)=1$, then ${\Cal
T}(j)$ is a tube of width $\alpha/2$ around the apartment fixed by
$j$.
\hfill\break
(iii) If $\alpha$ is odd and $\chi(\varepsilon)$ arbitrary, then
${\Cal T}(j)$ is a ball of radius $\alpha/2$ around the unique
fixed point of $j$, which is the midpoint of an edge.
\par\noindent
In particular, in all cases ${\Cal T}(j)\neq\emptyset$.
\endproclaim
\demo{Proof}
By the preceding lemma we only have to determine ${\Cal B}^j$. In
cases (i) and (ii) and after replacing $j$ by a scalar multiple we
may assume that $j^2=\varepsilon\cdot {\roman{id}}$. In case (i) we
may find a basis $e_1, e_2$ of $\Q_p^2$ such that $j$ has the
matrix
$$j=\pmatrix 0&\varepsilon\\ 1&0\\ \endpmatrix\ \ .$$
In this case $[\Lambda]$ with $\Lambda =[e_1, e_2]$ is the unique
vertex fixed by $j$ (${\roman{red}}_{\Lambda}(j)$ fixes no ${\Bbb
F}_p$-rational line in $\Lambda /p\Lambda$). In case (ii) we may
assume that $j^2={\roman{id}}$ and can find a basis $e_1, e_2$ of
$\Q_p^2$ such that $j$ has the matrix
$$j=\pmatrix 1&0\\ 0&-1\\\endpmatrix\ \ .\tag2.12$$
In this case the fixed point set is given by the apartment with
vertices $[\Lambda_i]$,
$$\Lambda_i= [p^ie_1, e_2]\ \ ,\ \ i\in\Z\ \ .\tag2.13$$
The case (iii) is also easy and is left to the reader.
\qed\enddemo
\proclaim{Corollary 2.6}
In the cases when $Z(j)(k)$ is a set of isolated points (cf.\
Corollary 2.3), this point set is of the following form.
\hfill\break
(i) $q(j)=\varepsilon\in\Z_p^{\times}$, $\chi(\varepsilon)=-1$. In
this case let $[\Lambda_0]$ be the unique fixed vertex of $j$. Then
$Z(j)(k)$ consists of two ordinary special points on ${\Bbb
P}_{[\Lambda_0]}$ (namely the two eigen lines of
${\roman{red}}_{[\Lambda_0]}(j)$).
\hfill\break
(ii) $q(j)=\varepsilon\in\Z_p^{\times}$, $\chi(\varepsilon)=1$. In
this case let (2.13) be the apartment fixed by $j$. Then
$$Z(j)(k)=\{ pt_{([\Lambda_i, \Lambda_{i-1}])};\ \ i\in\Z\}\ \ .$$
(iii) $q(j)=\varepsilon\cdot p$, $\chi(\varepsilon)$ arbitrary. Let
$\Delta$ be the edge containing the unique midpoint fixed by $j$.
Then
$$Z(j)(k)=\{ pt_{\Delta}\}\ \ .$$
\endproclaim
\demo{Proof}
This is simply a combination of the previous corollary with the
local analysis of Proposition 2.2.
\qed\enddemo
For the next corollary we need the following slight strengthening
of Lemma 2.4. The proof is identical.
\proclaim{Lemma 2.7}
Let $j\in GL_2(\Q_p)$ with ${\roman{tr}}(j)=0$, and let
$q(j)=\varepsilon\cdot p^{\alpha}$, $\varepsilon\in \Z_p^{\times}$.
Let $\Lambda\subset\Q_p^2$ be a lattice. Let
$$m_{[\Lambda]}(j):= {\roman{max}}\{ r;\ j(\Lambda)\subset p^r\Lambda\}\ \ .
\tag2.14$$
Then
$$m_{[\Lambda]}(j)= \alpha/2-d([\Lambda], {\Cal B}^j)\ \ .\tag2.15$$
\endproclaim
\proclaim{Corollary 2.8}
Let $j\in V$ with $q(j)=\varepsilon\cdot p^{\alpha}$, $\alpha\ge
1$. For all $[\Lambda]\in {\Cal T}(j)$ not on the boundary of
${\Cal T}(j)$ we have ${\Bbb P}_{[\Lambda]}\subset Z(j)$. If
$[\Lambda]$ is on the boundary of ${\Cal T}(j)$, i.e.\
$d([\Lambda], {\Cal B}^j)=\alpha/2$, then ${\Bbb P}_{[\Lambda]}\cap
Z(j)$ consists of a single superspecial point, namely the one
corresponding to the unique neighbouring vertex of $[\Lambda]$ in
${\Cal T}(j)$.
\endproclaim
We next turn to the intersection of two special cycles. Obviously
the intersection $Z(j, j')=Z(j)\cap Z(j')$ will only depend on the
$\Z_p$-span ${\bold j}$ of $j$ and $j'$, which we will assume to be
of rank 2. Let
$$(\ ,\ ): V\times V\to\Q_p\tag2.16$$
be the bilinear form corresponding to the quadratic form $q$,
$$(x,y)=q(x+y)-q(x)-q(y)\ \ .\tag2.17$$
Since $p\ne 2$, the restriction of $(\ ,\ )$ to ${\bold j}$ may be
diagonalized. {\it We will always assume that ${\bold j}$ is
non-degenerate.} Then we may choose a $\Z_p$-basis $j, j'$ of
${\bold j}$ such that the restriction of $(\ ,\ )$ to $\j$ has
matrix
$$T:=
\pmatrix
q(j)&\frac{1}{2}(j,j')\\
\frac{1}{2}(j,j')& q(j')
\endpmatrix
={\roman{diag}} (\varepsilon_1 p^{\alpha},
 \varepsilon_2p^{\beta})\tag2.18$$
with $\varepsilon_1, \varepsilon_2\in\Z_p^{\times}$, and $\alpha\ge
0$, $\beta\ge 0$. In particular $j$ and $j'$ anticommute,
$$jj'=-j'j\ \ .\tag2.19$$
We wish to determine the $k$-rational points of
$$Z(\j)=Z(j)\cap Z(j')\ \ .\tag2.20$$
We introduce
$${\Cal T}(\j)=\{ [\Lambda]\in {\Cal B};
\ {\Bbb P}_{[\Lambda]}\cap Z(\j)\ne \emptyset\}\ \ .\tag2.21$$
As in the case of a special single endomorphism, we start with a
local analysis. Let $\Lambda$ be a lattice such that the
corresponding vertex $[\Lambda]$ lies in ${\Cal T}(\j)$. Let
$${\frak m}={\roman{red}}_{\Lambda}(\j)\subset {\roman{End}}(\Lambda/p\Lambda)
\ \ .\tag2.22$$
Then ${\frak m}$ lies in the subspace of traceless matrices in
${\roman{End}}(\Lambda /p\Lambda)\simeq M_2({\Bbb F}_p)$. Let
$$(\ ,\ ): {\frak m}\times {\frak m}\to {\Bbb F}_p$$
be the bilinear form associated to the quadratic form $j\mapsto
-{\roman{det}}(j)$. Note that the matrix of this form with respect
to the basis ${\roman{red}}_{\Lambda}(j),
{\roman{red}}_{\Lambda}(j')$ is just the reduction modulo $p$ of
the matrix $T$ with respect to $j,j'$. Let ${\roman{rk}}\, {\frak
m}$ be the rank of this reduction.
\proclaim{Proposition 2.9}
With the previous notation, we have ${\roman{rk}}\, {\frak m}\le
1$. Furthermore,
\smallskip\noindent
(i) If ${\roman{rk}}\, {\frak m}=1$ and $\frak m$ represents 1,
there are two possibilities
\hfill\break
a) ${\roman{dim}}\, {\frak m}=2$. In this case ${\Bbb
P}_{[\Lambda]}\cap Z(\j)$ consists of a single superspecial point,
which is an isolated point of $Z(\j)\times_{\roman{Spf}\,
W}{\roman{Spec}}\, k.$
\hfill\break
b) ${\roman{dim}}\, {\frak m}=1$. In this case ${\Bbb
P}_{[\Lambda]}\cap Z(\j)$ consists of two superspecial points,
which are isolated points of $Z(\j)\times_{{\roman{Spf}}\,
W}{\roman{Spec}}\, k.$
\par\noindent
(ii) If ${\roman{rk}}\, {\frak m}=1$ and ${\frak m}$ does not
represent 1, then ${\roman{dim}}\, {\frak m}=1$ and ${\Bbb
P}_{[\Lambda]}\cap Z(\j)$ consists of two ordinary special points,
which are isolated points of $Z(\j)\times_{\roman{Spf}\,
W}{\roman{Spec}}\, k$.
\par\noindent
(iii) If ${\roman{rk}}\, {\frak m}=0$, then ${\roman{dim}}\, {\frak
m}\le 1$ and there are two possibilities.
\hfill\break
a) ${\roman{dim}}\, {\frak m}=1$. Then ${\Bbb P}_{[\Lambda]}\cap
Z(\j)$ consists of a single superspecial point. This is an isolated
point of $Z(\j)\times_{\roman{Spf}\, W}{\roman{Spec}}\, k$ iff
$p^2\not\hskip1.5pt\mid T$.
\hfill\break
b) ${\frak m}=0$. Then ${\Bbb P}_{[\Lambda]}\subset Z(\j)$.
\endproclaim
\demo{Proof}
Suppose by contradiction that ${\roman{rk}}\, {\frak m}=2$. Then
${\roman{red}}_{\Lambda}(j)$ and ${\roman{red}}_{\Lambda}(j')$
would be traceless invertible linear transformations which
anticommute. But then each has to interchange the two eigenspaces
of the other. But these eigen lines correspond precisely to
$${\Bbb P}_{[\Lambda]}\cap Z(j)\ \ \text{resp.}
\ \ {\Bbb P}_{[\Lambda]}\cap Z(j')\ \ .$$
It would follow that ${\Bbb P}_{[\Lambda]}\cap Z(\j)=\emptyset$,
contrary to our assumption.
\hfill\break
Let us prove (i). If ${\frak m}$ represents 1, we may choose the
basis $j,j'$ of $\j$ such that
$${\roman{red}}_{\Lambda}(j)^2={\roman{id}}\ \ ,
\ \ {\roman{red}}_{\Lambda}(j')^2=0\ \ .\tag2.23$$
If ${\roman{dim}}\, {\frak m}=2$, then
$${\roman{Im}}\, {\roman{red}}_{\Lambda}(j') =
{\roman{Ker}}\, {\roman{red}}_{\Lambda}(j')$$ is a one-dimensional
subspace of $\Lambda/p\Lambda$ preserved by
${\roman{red}}_{\Lambda}(j)$ since $jj'=-j'j$. The superspecial
point corresponding to this ${\Bbb F}_p$-rational point of ${\Bbb
P}_{[\Lambda]}$ is the unique point of ${\Bbb P}_{[\Lambda]}\cap
Z(\j)$. If ${\roman{dim}}\, {\frak m}=1$, then
${\roman{red}}_{\Lambda}(j')=0$ and ${\Bbb P}_{[\Lambda]}\cap
Z(\j)$ consists of the two superspecial points given by the two
${\Bbb F}_p$-rational eigenlines of ${\roman{red}}_{\Lambda}(j)$.
By corollary 2.3 in both cases we are dealing with isolated points
of $Z(j)\times_{\roman{Spf}\, W}{\roman{Spec}}\, k$.
\hfill\break
Let us now prove (ii). If ${\roman{rk}}\, {\frak m} =1$ and ${\frak
m}$ does not represent 1 we may assume that
$${\roman{red}}_{\Lambda}(j)^2\in {\Bbb F}_p^{\times}
- {\Bbb F}_p^{\times, 2}\ \ .$$
Then ${\roman{red}}_{\Lambda}(j)$ has two non-rational eigen lines.
Since ${\roman{red}}_{\Lambda}\, (j')$ anticommutes with
${\roman{red}}_{\Lambda}(j)$ it has to take any one of these lines
into the other one. Since ${\roman{red}}_{\Lambda}(j')$ is not
invertible, its restriction to at least one of the eigen lines has
to be zero. But then ${\roman{red}}_{\Lambda}(j')$ has to kill both
eigenlines since they are not ${\Bbb F}_p$-rational. Hence
${\roman{dim}}\, {\frak m}=1$ and ${\Bbb P}_{[\Lambda]}\cap Z(\j) =
{\Bbb P}_{[\Lambda]}\cap Z(j)$ consists of the two ordinary special
points corresponding to the two eigen lines. These again are
isolated points of $Z(j)\times_{\roman{Spf}\, W}{\roman{Spec}}\,
k$.
\hfill\break
Finally, let us prove (iii). If ${\roman{rk}}\, {\frak m}=0$, then
${\roman{dim}}\, {\frak m}\leq 1$. If ${\roman{dim}}\, {\frak m}=1$
we may assume that
$${\roman{red}}_{\Lambda}(j)\ne 0\ , \ {\roman{red}}_{\Lambda}(j)^2=0
\ ,\ {\roman{red}}_{\Lambda}(j')=0\ \ .$$
In this case ${\Bbb P}_{[\Lambda]}\cap Z(\j) = {\Bbb
P}_{[\Lambda]}\cap Z(j)$ consists of a single superspecial point
corresponding to the ${\Bbb F}_p$-rational line
$${\roman{Im}}\, {\roman{red}}_{\Lambda}(j)= {\roman{Ker}}\,
{\roman{red}}_{\Lambda}(j)\subset \Lambda /p\Lambda\ \ .$$ If
$p^2\not\hskip1.5pt\mid T$ then ${\roman{ord}}\, q(j)\leq 1$ and
$Z(j)\times_{\roman{Spf}\, W}{\roman{Spec}}\, k$ consists of
isolated points only. If $p^2\mid T$, then ${\roman{ord}}\,
q(j)\geq 2$ and we put $\Lambda'=j(\Lambda)+ p\Lambda$, cf.\ (2.9).
Then our superspecial point lies on ${\Bbb P}_{[\Lambda]}\cap {\Bbb
P}_{[\Lambda']}$. But
$$j'(\Lambda)\subset p\Lambda\ \ \text{and}
\ \ j(\Lambda')\subset p\Lambda'\ \ ,\tag2.24$$
cf.\ (2.9). Hence
$$j'(\Lambda') =j'(j(\Lambda)+p\Lambda)\subset p(j(\Lambda)+p\Lambda) =
p\Lambda'\ \ .$$ It follows that ${\Bbb P}_{[\Lambda']}\subset
Z(\j)$ and hence our intersection point is not isolated in $Z(\j)
\times_{\roman{Spf}\, W}{\roman{Spec}}\, k$.
\hfill\break
Finally, if ${\frak m}=(0)$, then obviously all lines in ${\Bbb
P}_{[\Lambda]}$ are preserved by ${\roman{red}}_{\Lambda}(j)$,
$\forall j\in\j$, i.e. ${\Bbb P}_{[\Lambda]}\subset Z(\j)$.
\qed\enddemo
We next turn to a global analysis of ${\Cal T}(\j)$. We obviously
have an inclusion
$${\Cal T}(\j)\subset {\Cal T}(j)\cap {\Cal T}(j')\ \ ,\tag2.25$$
for any set $j,j'$ of generators of the $\Z_p$-module $\j$. We
therefore start by determining ${\Cal T}(j)\cap {\Cal T}(j')$ in
the case when $j,j'$ diagonalize the bilinear form on $\j$, cf.\
(2.18). We shall do this according to the following table.

\input tabelle
$$
\vcenter{
\begintab{|c|c|c|c|}
\hline
$j\setminus j'$&$\beta$ even&$\beta$ even&$\beta$ odd
\cr
&$\chi(\varepsilon_2)=-1$ & $\chi(\varepsilon_2)=1$ &
\cr
\hline
$\alpha$ even &&&
\cr
$\chi(\varepsilon_1)=-1$&&&$\emptyset$
\cr
\hline
$\alpha$ even &&&
\cr
$\chi(\varepsilon_1)=1$ &&&
\cr
\hline
&&&
\cr
&&&
\cr
&&&
\cr
 $\alpha$ odd & $\emptyset$ &&
\cr
\hline
\endtab
}
\tag2.26
$$
\par\noindent
{\bf Remark 2.10.} We note that the existence of a basis $j,j'$ of
the type (2.18) imposes restrictions on the matrix $T$. Indeed, $j$
and $j'$ generate a non-commutative $\Q_p$-subalgebra of
$M_2(\Q_p)$ which therefore has to be all of $M_2(\Q_p)$. On the
other hand, the subalgebra generated by $j$ and $j'$ is just the
quaternion algebra over $\Q_p$ with invariant
$$(\varepsilon_1p^{\alpha}, \varepsilon_2p^{\beta})_p\in\Z/2\ \ .$$
We conclude that
$$(\varepsilon_1p^{\alpha}, \varepsilon_2p^{\beta})_p
=\chi(-1)^{\alpha\beta}\cdot\chi(\varepsilon_1)^{\beta}\cdot
\chi(\varepsilon_2)^{\alpha}
=1\ \ .\tag2.27$$
This excludes the following cases.
$$\align
&
\text{$\alpha$ odd and $\beta$ even and
$\chi(\varepsilon_2)=-1$, resp.}\tag2.28\\ &
\text{$\beta$ odd and $\alpha$ even and $\chi(\varepsilon_1)=-1$.}
\endalign
$$
$$\text{$\alpha$ and $\beta$ odd and $\chi(-\varepsilon_1
\varepsilon_2)=-1$.}\tag2.29$$
\par\noindent
Before stating the result, we point out that, since $j$ and $j'$
anticommute, each one of the induced automorphisms of ${\Cal B}$
will preserve the fixed point set of the other.
\hfill\break
(This remark incidentally gives another reason why the cases (2.28)
are excluded: in these cases one fixed point set is a vertex and
the other a midpoint.)
\proclaim{Proposition 2.11}
Let $j,j'$ be a basis of $\j$ such that the matrix of the bilinear
form is given by (2.18). The set ${\Cal T}(j)\cap {\Cal T}(j')$ is
of the following form.
\par\noindent
(i) $\alpha$ and $\beta$ even,
$\chi(\varepsilon_1)=\chi(\varepsilon_2)=-1$. In this case ${\Cal
T}(j)\cap {\Cal T}(j')$ is the ball with radius
${\roman{min}}(\alpha/2, \beta/2)$ around the unique common fixed
vertex of $j$ and $j'$.
\par\noindent
(ii) $\alpha$ and $\beta$ even, $\chi(\varepsilon_1)\ne
\chi(\varepsilon_2)$. Suppose e.g.\ that $\chi(\varepsilon_1)=-1$.
Then the unique vertex $[\Lambda_0]$ fixed by $j$ lies on the
apartment fixed by $j'$. In this case ${\Cal T}(j)\cap {\Cal
T}(j')$ is the intersection of the ball with radius $\alpha/2$
around $[\Lambda_0]$ with the tube of width $\beta/2$ around the
fixed apartment. The case when $\chi(\varepsilon_2)=-1$ is
analogous with the roles of $j$ and $j'$ interchanged.
\par\noindent
(iii) $\alpha$ and $\beta$ even, $\chi(\varepsilon_1)=
\chi(\varepsilon_2)=1$. In this case the fixed apartment
 of $j$ and the fixed apartment of $j'$ have a unique vertex in
 common. The set ${\Cal T}(j)\cap {\Cal T}(j')$ is the intersection
 of the tube of width $\alpha/2$ around the fixed apartment of $j$
 and the tube of width $\beta/2$ around the fixed apartment of
 $j'$. \par\noindent (iv) $\alpha$ odd, $\beta$ even and
 $\chi(\varepsilon_2)=1$ (resp.\ $\beta$ odd, $\alpha$ even and
 $\chi(\varepsilon_1)=1$). In this case the midpoint fixed by $j$
 lies on the fixed apartment of $j'$. The set ${\Cal T}(j)\cap
 {\Cal T}(j')$ is the intersection of the ball of radius $\alpha/2$
 around the midpoint fixed by $j$ with the tube of width $\beta/2$
 around the apartment fixed by $j'$.
\par\noindent
(v) $\alpha$ and $\beta$ odd. In this case $j$ and $j'$ fix the
same midpoint and ${\Cal T}(j)\cap {\Cal T}(j')$ is the ball with
radius ${\roman{min}}(\alpha/2, \beta/2)$ around it.
\endproclaim
\demo{Proof}
Let us prove the first statement in (iii). After correcting $j$ and
$j'$ by scalars we may assume that $j^2=j^{'2}={\roman{id}}$. We
may choose a basis $e_1, e_2$ of $\Q_p^2$ such that $j(e_1)=e_1,
j(e_2)=-e_2$, cf.\ (2.12). Since $j$ and $j'$ anticommute, $j'$ has
to permute the two eigenlines of $j$, hence $j'(e_1)=\pm e_2$,
$j(e_2)=\pm e_1$. If $\Lambda = [e_1, e_2]$ we see that $j'$
interchanges the two neighbouring vertices $[p^{\pm 1}e_1,e_2]$ of
$[\Lambda]$ in the fixed apartment of $j$, cf.\ (2.13). Since the
common fixed point set of $j$ and $j'$ is convex, it consists of
$[\Lambda]$ only.
\hfill\break
The other assertions are equally easy.
\qed\enddemo
Combining this now with our analysis of the special fibres of
$Z(j)$ and $Z(j')$ we obtain the following corollary.
\proclaim{Corollary 2.12}
Let $j,j'$ be a basis of $\j$ such that the matrix of the bilinear
form is given by (2.18). We have an equality ${\Cal T}(\j)={\Cal
T}(j)\cap {\Cal T}(j')$ except when $\alpha=\beta=0$. The latter
condition is equivalent to $Z(\j)=\emptyset$.
\endproclaim
\demo{Proof}
If $[\Lambda]$ does not lie on the boundary of either ${\Cal T}(j)$
or ${\Cal T}(j')$ we have by Corollary 2.8 that ${\Bbb
P}_{[\Lambda]}\subset Z(j)$ or ${\Bbb P}_{[\Lambda]}\subset Z(j')$
and hence $[\Lambda]\in {\Cal T}(\j)$ iff $[\Lambda]\in {\Cal
T}(j)\cap {\Cal T}(j')$. Therefore the only problematical vertices
are those which lie on the boundary of both ${\Cal T}(j)$ and
${\Cal T}(j')$. Now one goes through all combinations as to whether
${\Cal T}(j)$ resp.\ ${\Cal T}(j')$ is a ball (around a vertex or a
midpoint of an edge) resp.\ a tube, using the second assertion of
corollary 2.8. In the analysis one uses the fact (comp.\
\cite{\krHB}, Lemma 6.8) that if
$$d([\Lambda], {\Cal B}^j)\leq d([\Lambda], {\Cal B}^{j'})\ \ ,$$
then the unique geodesic from $[\Lambda]$ to ${\Cal B}^{j'}$ first
runs along the geodesic from $[\Lambda]$ to ${\Cal B}^j$ and then
stays inside ${\Cal B}^j$. The analysis of the extreme cases in
which both $Z(j)\times_{\roman{Spf}\, W}{\roman{Spec}}\, k$ and
$Z(j')\times_{\roman{Spf}\, W}{\roman{Spec}}\, k$ are sets of
isolated points is left to the reader (use Proposition 2.9).
\qed\enddemo
\par\noindent
\proclaim{Corollary 2.13}
Assume that $\j$ is non-degenerate of rank 2. Then
$Z(\j)\times_{\roman{Spf}\, W}{\roman{Spec}}\, k$ is a projective
scheme over ${\roman{Spec}}\, k$. Furthermore, it is either empty,
or a finite set of points, or a connected finite union of
projective lines.
\endproclaim
\demo{Proof}
The first statement is equivalent to the assertion that
$Z(\j)\times_{\roman{Spf}\, W}{\roman{Spec}}\, k$ is quasi-compact,
i.e.\ to the assertion that ${\Cal T}(\j)$ is finite. The
finiteness of ${\Cal T}(\j)$ follows from Proposition 2.11.
\hfill\break
The proposition also implies that, except in extreme cases,
$Z(\j)\times_{\roman{Spf}\, W}{\roman{Spec}}\, k$ is a connected
union of projective lines. The extreme cases are, case by case:
\par\noindent
(i) $\alpha=\beta=0$. Then $Z(\j)=\emptyset$, cf.\ Corollary 2.12.
\par\noindent
(ii) $\alpha=0$, $\beta >0$ (or $\beta =0$, $\alpha >0$). In this
case $Z(\j)\times_{\roman{Spf}\, W}{\roman{Spec}}\, k$ consists of
two ordinary special points if $\chi(\varepsilon_1)=-1$, or of a
finite set of superspecial points if $\chi(\varepsilon_1)=1$. The
case where $\beta=0$ and $\alpha
>0$ is symmetric.\enddemo

\subheading{\Sec 3. Local equations for the special cycles}
\par\noindent
In this section we will write down equations which describe a
special cycle $Z(j)$ in a neighbourhood of a point $x$ in the
special fibre. As mentioned in the introduction, the basic idea of
how to do this is due to Genestier \cite{\genestier}. In the beginning of
this section we fix a special endomorphism $j$ with
$q(j)=\varepsilon\cdot q^{\alpha}\in
\Z_p - \{ 0\}$.
\proclaim{Theorem 3.1} {\rm (Genestier)} The Drinfeld isomorphism (1.2)
induces isomorphisms
$$\align
Z(j)\simeq \hat\Omega^j
\mathop{\times}\nolimits_{\roman{Spf}\,\Z_p} {\roman{Spf}}\, W &
\quad\text{if}\ \alpha=0
\tag3.1\\
Z(j)\simeq \hat\Omega^{\roman{id} +j}
\mathop{\times}\nolimits_{\roman{Spf}\, \Z_p}{\roman{Spf}}\,
 W &
\quad\text{if}\ \alpha >0\tag3.2\endalign$$
(fixed point formal schemes for the action of elements of
$GL_2(\Q_p)$ on $\hat\Omega$).
\endproclaim
\demo{Proof}
Suppose first that $j$ is an arbitrary element of $G(\Q_p)^0$ (cf.\
(1.1)) and define a closed formal subscheme $Z(j)$ in the same way
as for a special endomorphism, cf.\ Definition 2.1. If
$(X,\varrho)\in Z(j)(S)$, then the unique isogeny $j_X$ lifting
$\varrho\circ j\circ \varrho^{-1}$ is of height 0, and hence an
automorphism of $X$. It follows that $(X,\varrho)$ is a fixed point
for the action of $j^{-1}$ on ${\Cal M}$. Conversely, a fixed point
defines an element of $Z(j)(S)$. Therefore
$Z(j)=\hat\Omega^j\mathop{\times}\nolimits_{\roman{Spf}\,\Z_p}
{\roman{Spf}}\, W$, since the Drinfeld isomorphism is equivariant
for the action of $G(\Q_p)^0$, cf.\ (1.2).
\hfill\break
Suppose now that $j$ is a special endomorphism with
${\roman{ord}}\, {\roman{det}}\, j>0$. But then
${\roman{det}}({\roman{id}} +j)= 1+{\roman{det}}(j)$, hence
${\roman{id}} +j\in G(\Q_p)^0$. The result therefore follows from
the previous case since $Z({\roman{id}}+j)= Z(j)$.
\qed\enddemo

Our next task will be to write down equations describing the fixed
point schemes. Put $\tilde j=j$ if $\alpha=0$ and $\tilde
j={\roman{id}} +j$ if $\alpha >0$. Let $x\in \hat\Omega^{\tilde
j}(k)$. We distinguish two cases.
\par\noindent
{\bf Case $x$ ordinary:} After replacing $x$ by $gx$ and $j$ by
$gjg^{-1}$ we may assume that $x$ lies in
$\hat\Omega_{[\Lambda_0]}$, where $\Lambda_0=[e_1, e_2]$ denotes
the standard lattice. Furthermore, ${\roman{ord}}\, {\roman{det}}\,
\tilde j=0$. But then, since $\tilde j$ fixes $x$ and hence $[\Lambda_0]$, it
follows that $j\in M_2(\Z_p)$ and $\tilde j\in GL_2(\Z_p)$. We
write
$$j=\pmatrix a&b\\ c&d\\\endpmatrix\ \ .\tag3.3$$
Then, using the description (1.8) of the action of $GL_2(\Z_p)$ on
$\hat\Omega_{[\Lambda_0]}$ we obtain the following equation for
$\hat\Omega^{\tilde j}_{[\Lambda_0]}$,
$$-\tilde bT^2+(\tilde a-\tilde d)\cdot T+\tilde c=0\ \ ,\tag3.4$$
where $\tilde a, \tilde b=b, \tilde c=c, \tilde d$ are the
coefficients of $\tilde j$. The equation for $Z(j)$ therefore is,
regardless of whether ${\roman{ord}}\, {\roman{det}}\, j=0$ or
$>0$,
$$bT^2-2aT-c=0\ \ .\tag3.5$$
\par\noindent
{\bf Case $x$ superspecial:} In this case, after replacing $x$ by
$gx$ and $\tilde j$ by $g\tilde jg^{-1}$ we may assume that
$x=pt_{\Delta_0}$, where $\Delta_0= ([\Lambda_0], [\Lambda_1])$ is
the standard edge (1.9). In this case, $\tilde j$ will lie in the
Iwahori subgroup (1.11). We write
$$j=\pmatrix a&b\\ c&d\\\endpmatrix =
\pmatrix a&pb_0\\ c&d\\\endpmatrix\ \ .\tag3.6$$
In this case we obtain from (1.12) the following pair of equations
describing $\hat\Omega^{\tilde j}_{\Delta_0}$,
$$\align
T_0(\tilde b_0T_0-(\tilde a-\tilde d)-\tilde cT_1)&=0\tag3.7\\
T_1(\tilde b_0T_0-(\tilde a-\tilde d)-\tilde cT_1) &=0,\\
\endalign$$
where $\tilde a, p\tilde b_0= pb_0, \tilde c=c, \tilde d$ are the
coefficients of $\tilde j$. We have used the fact that $T_0T_1=p$.
The equations for $Z(j)$ are, regardless of whether
${\roman{ord}}\, {\roman{det}}\, j=0$ or $>0$,
$$\align
T_0(b_0T_0-2a-cT_1)&=0\tag3.8\\
T_1(b_0T_0-2a-cT_1)&=0.\\\endalign$$ The reader may reassure
himself that the locus defined by the equations (3.8) induces on
the open formal subscheme $\hat\Omega_{[\Lambda_0]}$ of
$\hat\Omega_{\Delta_0}$, cf.\ (1.14), the locus defined by the
equation (3.5).

In the following statement we denote by an upper index
${}^{\roman{ord}}$ the intersection with the ordinary locus of
${\Cal M}$ resp.\ $\hat\Omega$, i.e.\ the open formal subscheme
formed by the complement of the superspecial points. Define
$${\roman{mult}}_{[\Lambda]}(j) ={\roman{max}}(m_{[\Lambda]}(j), 0)\ \ ,\tag3.9$$
with $m_{[\Lambda]}(j)$ as defined in Lemma 2.7.
\proclaim{Proposition 3.2}
The closed formal subscheme $Z(j)^{\roman{ord}}$ of
$(\hat\Omega\times_{\roman{Spf}\,\Z_p}{\roman{Spf}}\,
W)^{\roman{ord}}$ is a divisor. We have the following equality of
divisors,
$$Z(j)^{\roman{ord}}=\sum_{[\Lambda]\in
{\Cal T}(j)} {\roman{mult}}_{[\Lambda]}(j)\cdot {\Bbb
P}_{[\Lambda]}^{\roman{ord}}\ \ ,\tag3.10$$ unless
$q(j)=\varepsilon\cdot p^{\alpha}$ where
$\varepsilon\in\Z_p^{\times}$ and with $\alpha$ even and
$\chi(\varepsilon)=-1$, in which case
$$Z(j)^{\roman{ord}}=Z(j)^{\roman{ord}, h} +\sum_{[\Lambda]\in{\Cal T}(j)}
{\roman{mult}}_{[\Lambda]}(j)\cdot {\Bbb
P}^{\roman{ord}}_{[\Lambda]}\ \ ,\tag3.11$$ where
$Z(j)^{\roman{ord},h}$ (the ``horizontal divisor'') is isomorphic
to the disjoint union of two copies of ${\roman{Spf}}\, W$ and
meets the special fibre in two ordinary special points of ${\Bbb
P}_{[\Lambda(j)]}$. Here $[\Lambda(j)]$ denotes the unique vertex
fixed by $j$.
\endproclaim
\demo{Proof}
We wish to determine $Z(j)\cap
(\hat\Omega_{[\Lambda]}\times_{\roman{Spf}\,\Z_p}{\roman{Spf}}\,
W)$ for a vertex $[\Lambda]$ where this intersection is non-empty.
Note that this implies that $m_{[\Lambda]}(j)\ge 0$, cf.\
Proposition 2.1. Replacing $[\Lambda]$ by $[g\Lambda]$ and $j$ by
$gjg^{-1}$ we may assume that $[\Lambda]=[\Lambda_0]$ is the
standard lattice. Let us write
$$j=\pmatrix a&b\\ c&-a\\\endpmatrix =p^m\cdot\pmatrix \bar a&\bar b\\
\bar c& -\bar a\endpmatrix =p^m\cdot\bar j\ \ ,\ \
m=m_{[\Lambda_0]}(j)\ \ ,\tag3.12$$ where $\bar a, \bar b, \bar
c\in\Z_p$ are not simultaneously divisible by $p$. The equation
(3.5) for $Z(j)\cap
(\hat\Omega_{[\Lambda_0]}\times_{\roman{Spf}\,\Z_p}{\roman{Spf}}\,
W)$ may be written as
$$p^m\cdot (\bar bT^2-2\bar a T-\bar c)=0\ \ .\tag3.13$$
It follows already that $Z(j)\cap (\hat\Omega_{[\Lambda]}
\times_{\roman{Spf}\,\Z_p}{\roman{Spf}}\, W)$ is a divisor. The
second factor in the LHS of (3.13) is not divisible by $p$, i.e.\
is a unit after localizing at the ideal $(p)$. Since ${\Bbb
P}_{[\Lambda_0]}^{\roman{ord}}$ is defined by $p=0$ in
$\hat\Omega_{[\Lambda_0]}$, the multiplicity of ${\Bbb
P}_{[\Lambda_0]}$ in the divisor is equal to $m$, as asserted. We
still have to determine the zero set $Z(j)^{\roman{ord},h}\cap
(\hat\Omega_{[\Lambda_0]}\times_{\roman{Spf}\,\Z_p}{\roman{Spf}}\,
W)$ of the second factor of (3.13). Now
$$Z(j)^{\roman{ord},h}\cap (\hat\Omega_{[\Lambda_0]}\times_{\roman{Spf}\,\Z_p}{\roman{Spf}}\, W)
=Z(\bar j)\cap (\hat\Omega_{[\Lambda_0]} \times_{\roman{Spf}\,\Z_p}{\roman{Spf}}\, W)\ \ .\tag3.14$$
Furthermore, ${\Bbb P}_{[\Lambda_0]}^{\roman{ord}}\cap Z(\bar j)$
is precisely the set of lines in ${\Bbb
P}_{[\Lambda_0]}^{\roman{ord}}$ preserved by
${\roman{red}}_{\Lambda_0}(\bar j)$. Since
${\roman{red}}_{\Lambda_0}(\bar j)\ne 0$, Proposition 2.2 implies
that this last intersection is empty, unless
${\roman{red}}_{\Lambda}(\bar j)$ is an invertible transformation
of $\Lambda/p\Lambda$ and this intersection consists of two
ordinary special points. In the latter case we have $j(\Lambda_0)=
p^m\Lambda_0$, hence $[\Lambda_0]$ is a fixed vertex under $j$. It
follows that $\alpha$ is even and $\chi(\varepsilon)=-1$, cf.\
Corollary 2.5, hence $[\Lambda_0]$ is the unique vertex fixed by
$j$. Furthermore, the extension
$$\Z_p[T] /\bar bT^2-2\bar a T-\bar c\tag3.15$$
is unramified, which finishes the proof.
\qed\enddemo

We next turn to the local equations of $Z(j)$ at a superspecial
point.
\proclaim{Proposition 3.3}
Let $q(j)=\varepsilon\cdot p^{\alpha}$, $\varepsilon\in
\Z_p^{\times}$, $\alpha\ge 0$
\par\noindent
(i) If $\alpha=0$ and $\chi(\varepsilon)=-1$, then $Z(j)$ contains
no superspecial points, cf.\ Corollary 2.6, (i)
\par\noindent
(ii) If $\alpha=0$ and $\chi(\varepsilon)=1$, then $Z(j)$ is equal
to the disjoint union of the {\rm reduced} superspecial points
corresponding to the edges in the fixed apartment ${\Cal B}^j$,
cf.\ Corollary 2.6, (ii).
\par\noindent
(iii) If $\alpha\ge 1$, then $Z(j)$ is purely one-dimensional and
contains superspecial points. Let $x=pt_{\Delta}$ be one of them,
where $\Delta$ is an edge contained in ${\Cal T}(j)$. Then, locally
around $pt_{\Delta}$, $Z(j)$ is the union of a divisor with support
in the special fibre and an embedded component at $pt_{\Delta}$,
except in the case where $\alpha$ is odd and $\Delta=\Delta(j)$ is
the edge containing the midpoint fixed by $j$. In the latter case,
if $\alpha=1$, then $Z(j)$ is the union of $Z(j)^h$ and an embedded
component at $pt_{\Delta(j)}$, where $Z(j)^h$ is a divisor
isomorphic to ${\roman{Spf}}\, W'$ where $W'$ is the ring of
integers in the ramified quadratic extension of $W_{\Q}$. Finally,
in the latter case, if $\alpha >1$, then $Z(j)$ is locally at
$pt_{\Delta(j)}$ the union of a divisor $Z(j)^h= {\roman{Spf}}\,
W'$, an embedded component at $pt_{\Delta}$ and a divisor with
support in the special fibre.
\endproclaim
\demo{Proof}
As before we may assume that $\Delta_0= ([\Lambda_0], [\Lambda_1])$
is the standard simplex. Let us first assume that $[\Lambda_0]$ is
strictly closer than $[\Lambda_1]$ to the fixed point set ${\Cal
B}^j$. We therefore have
$$m:= m_{[\Lambda_0]}(j)\ \ ,\ \ m_{[\Lambda_1]}(j)= m-1\ge 0
\ \ .\tag3.16$$
In terms of the canonical basis $e_1,e_2$ of $\Lambda_0$ we may
write
$$j=\pmatrix a&b\\ c&-a\\\endpmatrix =p^m\cdot
\pmatrix \bar a&\bar b\\ \bar c&-\bar a\\\endpmatrix\tag3.17$$
where $\bar a,\bar b,\bar c$ are not simultaneously divisible by
$p$. On the other hand, $b=p\cdot b_0$, cf.\ (3.6), and in terms of
the canonical basis $pe_1, e_2$ of $\Lambda_1$ the matrix of $j$ is
$$\pmatrix a&b_0\\ pc&-a\\\endpmatrix =p^{m-1}
\pmatrix p\bar a&\bar b_0\\ p^2\bar c&-p\bar a\\\endpmatrix
\ \ ,$$
with $b_0=p^{m-1}\cdot \bar b_0$. Since $m_{[\Lambda_1]}(j)= m-1$
we conclude that $\bar b_0$ is a unit. Therefore the system of
equations (3.8) for $Z(j)\cap
(\hat\Omega_{\Delta_0}\times_{\roman{Spf}\,\Z_p}{\roman{Spf}}\, W)$
is given by
$$\align
p^{m-1}\cdot T_0\cdot (\bar b_0\cdot T_0-2p\bar a-p\bar c\cdot T_1)
=T_0^{m+1}\cdot T_1^{m-1}\cdot u
&
=0\tag3.18\\
p^{m-1}\cdot T_1(\bar b_0T_0- 2p\bar a-p\bar c\cdot T_1)
=T_0^m\cdot T_1^m\cdot u
&
=0\ \ ,\\
\endalign$$
where
$$u=\bar b_0-2\bar a\cdot T_1-\bar c\cdot T_1^2$$
is a unit in the local ring of $pt_{\Delta_0}$. The above system of
equations describes the same locus as
$$T_0^2\cdot (T_0T_1)^{m-1}= (T_0T_1) \cdot (T_0T_1)^{m-1}=0\ \ .
\tag3.19$$
In this case we therefore see that $Z(j)$ is the union of a divisor
with support in the special fibre, and an embedded component at
$pt_{\Delta_0}$.
\hfill\break
Now consider the case when $[\Lambda_0]$ and $[\Lambda_1]$ have the
same distance to the fixed point set. There are two alternatives.
\par\noindent
{\bf First case:} ${\Cal B}^j$ is the midpoint of $\Delta_0$ (hence
$\alpha$ is odd). In this case $j(\Lambda_0)$ is a scalar multiple
of $\Lambda_1$ and $j(\Lambda_1)$ a scalar multiple of $\Lambda_0$,
$$\align
j(\Lambda_0) &
=p^m\cdot \Lambda_1\tag3.20\\
j(\Lambda_1) &
=p^{m+1}\Lambda_0\ \ .\\
\endalign$$
In this case we conclude in terms of the equation (3.17) that
$$p\mid \bar a\ \ ,\ \ p\mid\bar b\ \ ,\ \ {\roman{ord}}\, {\roman{det}}
\pmatrix \bar a&\bar b\\ \bar c&-\bar a\\\endpmatrix=1\ \ .\tag3.21$$
Therefore $\bar b=p\cdot \bar b_0$, and $\bar b_0$ and $\bar c$ are
both units. Hence the system of equations (3.8) is given by
$$\align
p^m\cdot T_0(\bar b_0 T_0-2\bar a-\bar c T_1) &
=0\tag3.22
\\
p^m\cdot T_1(\bar b_0T_0-2\bar a -\bar c T_1) &
=0\ \ .\\
\endalign$$
Up to a unit in the local ring at $pt_{\Delta}$ the second factor
is equal to $\bar b_0T_0-\bar cT_1$, hence $Z(j)$ is in a
neighbourhood of $pt_{\Delta}$ defined by the equations
$$T_0^{m+1}\cdot T_1^m (\bar b_0T_0- \bar c T_1)= T_0^mT_1^{m+1}(\bar b_0
T_0-\bar c T_1)=0\ \ .\tag3.23$$ We therefore see that $Z(j)$
locally at $pt_{\Delta}$ is the union of the divisor $Z(j)^h$, with
$$Z(j)^h= {\roman{Spf}}\, W[T_0, T_1]^{\wedge} /
(\bar b_0T_0- \bar c T_1, T_0T_1-p)
\simeq {\roman{Spf}}\, W[X] /(X^2-p)\tag3.24$$
$(p\ne 2)$, an embedded component at $pt_{\Delta}$, and a divisor
with support in the special fibre provided that $m>0$, i.e.,
$\alpha>1$. If $\alpha=1$, then $m=0$ and there is no divisor with
support in the special fibre present.
\par\noindent
{\bf Second case:} $\Delta_0$ lies in the fixed apartment ${\Cal
B}^j$ (hence $\alpha$ is even and $\chi(\varepsilon)=1$). In this
case
$$\align
j(\Lambda_0) &
=p^m\Lambda_0\tag3.25\\
j(\Lambda_1) &
=p^m\Lambda_1\ \ .\\\endalign$$
In this case we conclude in terms of the equation (3.17) that
$$p\mid \bar b\ \ ,\ \ {\roman{det}}
\pmatrix \bar a&\bar b\\ \bar c&-\bar a\\\endpmatrix\in
\Z_p^{\times}\ \ \tag3.26$$
We conclude that $\bar a$ is a unit. The system of equations (3.8)
may be written as
$$\align
p^m\cdot T_0(\bar b_0 T_0-2\bar a-\bar cT_1) &
=0\tag3.27\\
p^m\cdot T_1(\bar b_0T_0-2\bar a-\bar cT_1) &
=0\ \ ,\\
\endalign$$
where $p\cdot \bar b_0=\bar b$. The second factor is a unit in the
local ring at $pt_{\Delta_0}$ and therefore $Z(j)$ is locally
around $pt_{\Delta_0}$ described by
$$T_0^{m+1}\cdot T_1^m=T_0^m\cdot T_1^{m+1}=0\ \ .\tag3.28$$
If $m\ge 1$, i.e.\ $\alpha\ge 2$, then $Z(j)$ is the union of a
divisor with support in the special fibre and an embedded
component. If $m=\alpha=0$, then $Z(j)$ consists of the reduced
origin only ($T_0=T_1=0$).
\qed\enddemo
\proclaim{Corollary 3.4}
Suppose that $q(j)=\varepsilon\cdot p^{\alpha}$ with $\alpha$ even
and $\chi(\varepsilon)=1$. Then the support of $Z(j)$ is contained
in the special fibre.
\qed\endproclaim
\medskip\noindent
{\bf Remark 3.5.} Consider the special case $\alpha=0$,
$\chi(\varepsilon)=1$ of the previous corollary, i.e.\ (ii) of
Proposition 3.3. After correcting $j$ by a unit scalar, we may
assume that $j^2={\roman{id}}$. Let $(X,\iota)\in Z(j)(k)$. Then
$j$ induces a $O_B$-stable decomposition into its $\pm
1$-eigenspaces
$$(X,\iota)=(Y_1, \iota_1)\times (Y_2, \iota_2)\ \ .$$
But then $(Y_i, \iota_i)$ $(i=1,2)$ are both isomorphic to the
$p$-divisible group of a supersingular elliptic curve with its
$O_B$-action. It is well-known that $(Y_i,\iota_i)$ has no
non-trivial deformations, hence neither has $(X,\iota)$. This
argument gives an alternative proof of assertion (ii) of
Proposition 3.3.
\hfill\break
Later, when we make the connection with the global situation, we
will see a more convincing reason for the assertion of Corollary
3.4.

We next turn to the intersection of two special cycles $Z(j)$ and
$Z(j')$. We use the notation of (2.16), in particular the
$\Z_p$-module $\j$ is supposed to be nondegenerate of rank 2. We
again write $Z(\j)=Z(j)\cap Z(j')$.
\proclaim{Proposition 3.6}
$Z(\j)$ has support in the special fibre.
\endproclaim
\demo{Proof}
We momentarily change notations and as in section 2 denote by
$j,j'$ generators of $\j$ which diagonalize the bilinear form, with
matrix (2.18). We assume $Z(\j)\ne\emptyset$. We go through the
table (2.26), leaving aside the $\emptyset$-entries. Thanks to
Corollary 3.4 there are only two non-trivial possibilities.
\par\noindent
{\bf Case 1:} {\it $\alpha$ and $\beta$ even, $\chi(\varepsilon_1)=
\chi(\varepsilon_2)=-1$.} In this case the only part of $Z(j)$
resp.\ $Z(j')$ with support not in the special fibre is the
horizontal divisor $Z(j)^{\roman{ord}, h}$ resp.\
$Z(j')^{\roman{ord},h}$ which meets the special fibre in two
ordinary special points $x_1, x_2\in {\Bbb P}_{[\Lambda]}$ resp.\
$x'_1, x'_2\in {\Bbb P}_{[\Lambda]}$, cf.\ Propositions 3.2 and
3.3. Here $[\Lambda]= [\Lambda(j)]=[\Lambda(j')]$ is the unique
vertex fixed by $j$ and $j'$. We may write
$$j=p^{\alpha/2}\cdot\bar j\ \ ,\ \ j'=p^{\beta/2}\cdot \bar j'\ \ ,
\tag3.29$$
where ${\roman{red}}_{\Lambda}(\bar j)$ and
${\roman{red}}_{\Lambda}(\bar j')$ are invertible traceless
endomorphisms of $\Lambda/p\Lambda$ and whose eigen lines are $x_1$
and $x_2$ resp.\ $x'_1$ and $x'_2$, cf.\ proof of Proposition 3.2.
Since these two endomorphisms anticommute, each one interchanges
the eigen lines of the other (comp.\ proof of Proposition 2.9,
(ii)), hence they have no common eigen line. We conclude that in
this case
$$Z(j)^h\cap Z(j')^h =\emptyset\ \ .\tag3.30$$
\par\noindent
{\bf Case 2:} {\it $\alpha$ and $\beta$ odd.} In this case $j$ and
$j'$ both fix the midpoint of an edge $\Delta$ and the only part of
$Z(j)$ resp.\ $Z(j)$ with support not in the special fibre is the
horizontal divisor $Z(j)^h$ resp.\ $Z(j')^h$ passing through
$pt_{\Delta}$, cf.\ Proposition 3.3. We write, as in (3.12)
$$j=p^{\frac{\alpha-1}{2}}\cdot\bar j,\ \ \text{resp.}\ \
j'= p^{\frac{\beta-1}{2}}\cdot\bar j'\tag3.31$$ where
$$\bar j= \pmatrix \bar a&\bar b\\ \bar c&-\bar a\endpmatrix\ ,
\ \ \text{resp.}\ \ \bar j'=\pmatrix \bar a'&\bar b'\\ \bar c'&-\bar a'
\\\endpmatrix\ \ ,$$
with $p$ dividing $\bar a, \bar a', \bar b, \bar b'$ and where, if
$\bar b=p\cdot b_0$ resp.\ $\bar b'=p\cdot \bar b'_0$, the elements
$\bar b_0, \bar b'_0, \bar c, \bar c'$ are all units, cf.\ (3.21).
Then $Z(j)^h$ resp.\ $Z(j')^h$ is defined inside ${\roman{Spf}}\, W
[T_0, T_1]^{\wedge} / T_0T_1-p$ by the linear equation
$$\bar b_0 T_0-\bar cT_1=0\ \ \text{resp.}\ \
\bar b'_0T_0- \bar c'T_1 =0\ \ ,\tag3.32$$
cf.\ (3.24). We claim that these equations describe distinct
subschemes which will establish the proposition. Suppose to the
contrary that
$$\bar b_0\bar c' =\bar b'_0 \bar c\ \ .\tag3.33$$
However, $j$ and $j'$ anticommute, i.e.
$$\bar a \bar a' +p\cdot \bar b_0 \bar c' =-(\bar a\bar a' +
p\cdot \bar b'_0\bar c)\ \ .\tag3.34$$ Taking these identities
together we therefore obtain
$$\bar a\bar a' +p\bar b_0\bar c'=0\ \ .$$
But the first summand is divisible by $p^2$, whereas the second
summand is only divisible by $p$ which is the desired
contradiction.
\enddemo
\par\noindent
We note the following consequence of the proof.
\proclaim{Corollary 3.7}
Let $j$ and $j'$ be generators of $\j$ which diagonalize the
bilinear form, with matrix (2.18). For the horizontal divisors
$Z(j)^h$ and $Z(j')^h$ we have $Z(j)^h\cap Z(j')^h=\emptyset$
unless $\alpha$ and $\beta$ are both odd. In the latter case
$Z(j)^h$ and $Z(j')^h$ intersect transversally at a unique
superspecial point.
\endproclaim
Here, as in the rest of the paper, we formally set
$Z(j)^h=\emptyset$ if $q(j)=\varepsilon\cdot p^{\alpha}$ with
$\alpha$ even and $\chi(\varepsilon)=1$.

\subheading{\Sec 4. Intersection calculus of special cycles}
\par\noindent
Our next aim will be to determine the intersection numbers of
special cycles. Before doing this we will have to explain briefly
the kind of intersection theory we will want to use. The
definitions and facts we need are essentially all well-known but we
could not find a reference for them. The exposition in the
literature closest to our needs is Deligne's \cite{\deligne}.

Let $(S, \eta, s)$ be the spectrum of a discrete valuation ring and
let $X$ be an $S$-scheme $f:X\to S$ which is regular and locally of
finite type and flat over $S$ with all fibres of pure dimension 1.
Let $Y=f^{-1}(s)$. For any coherent ${\Cal O}_X$- module ${\Cal F}$
with support proper over $S$ and contained in $Y$ we may define its
{\it Euler-Poincar\'e characteristic}\/
$$\chi({\Cal F})= {\roman{lg}}\, f_*{\Cal F}-{\roman{lg}}\, R^1f_*{\Cal F}
\ \ .\tag4.1$$
Then $\chi$ is additive in short exact sequences, and if ${\Cal F}$
is a skyscraper sheaf concentrated in $y\in Y$, then
$$\chi({\Cal F}) ={\roman{lg}}_{{\Cal O}_y}({\Cal F}_y)\cdot [k(y):k(s)]\ \ .
\tag4.2$$
If $K$ is a complex of ${\Cal O}_X$-modules with finitely many
cohomology sheaves which are coherent and of the above type we set
$$\chi(K) =\sum(-1)^i\chi({\Cal H}^i(K))\ \ .\tag4.3$$
\proclaim{Lemma 4.1}
Let ${\Cal F}$ be a skyscraper sheaf concentrated at $y\in Y$ and
let ${\Cal G}$ be any coherent ${\Cal O}_X$-module with
${\roman{supp}}\, {\Cal G}$ a proper closed subscheme of $X$
locally at $y$. Then
$$\chi({\Cal F}\otimes^{\Bbb L} {\Cal G})=\chi({\Cal F}\otimes{\Cal G})
-\chi({\Cal T}or_1 ({\Cal F}, {\Cal G}))=0\ \ .\tag4.4$$
\endproclaim
\par\noindent
[Here ${\Cal F}\otimes^{\Bbb L} {\Cal G}$ denotes the derived
tensor product of ${\Cal F}$ and ${\Cal G}$, a complex of
quasicoherent ${\Cal O}_X$-modules well defined in the derived
category, with cohomology sheaves in degree $-i$ equal to ${\Cal
T}or_i^{{\Cal O}_X}({\Cal F}, {\Cal G})$. Since $X$ is supposed to
be regular of dimension 2, the Tor-terms for $i\geq 2$ vanish.]

\demo{Proof}
The sheaves ${\Cal T}or_i({\Cal F}, {\Cal G})$ are skyscraper
sheaves concentrated in $y$, hence
$$\chi({\Cal F}\otimes^{\Bbb L}{\Cal G})= \sum (-1)^i {\roman{lg}}\,
{\roman{Tor}}_i^{{\Cal O}_y} ({\Cal F}_y, {\Cal G}_y)\ \ .$$ Now
${\Cal G}_y$ has a resolution of length 2 by free ${\Cal
O}_y$-modules of finite rank,
$$0\to {\Cal O}_y^{n_1}\to {\Cal O}_y^{n_2}\to {\Cal G}_y\to 0
\ \ .\tag4.5$$
The hypothesis on ${\Cal G}$ implies that $n_1=n_2$. Tensoring
(4.5) with ${\Cal F}_y$ we obtain
$$\chi({\Cal F}_y\otimes_{{\Cal O}_y}^{\Bbb L}{\Cal G}_y) =
n_1\cdot {\roman{lg}}\, {\Cal F}_y- n_2{\roman{lg}}\, {\Cal F}_y
=0.\ \ \qed$$
\enddemo

Let now $Z$ and $Z'$ be closed subschemes of $X$ such that $(Z\cap
Z')_{\roman{red}}$ is contained in $Y$ and is proper over $S$. We
then define their {\it intersection number}\/ by
$$\align
(Z,Z') &
=\chi({\Cal O}_Z\otimes^{\Bbb L}{\Cal O}_{Z'})\tag4.6\\
&
=\chi({\Cal O}_Z\otimes {\Cal O}_{Z'})- \chi({\Cal T}or_1
({\Cal O}_{Z}, {\Cal O}_{Z'}))\ \ .\\
\endalign$$
If $Z$ and $Z'$ are divisors which meet in a finite number of
points contained in $Y$ then
$$(Z,Z')=\sum_{y\in Z\cap Z'} {\roman{lg}} ({\Cal O}_Z\otimes
{\Cal O}_{Z'})_y\tag4.7$$ is the intersection number in the most
naive sense. In \cite{\deligne} Deligne defines the intersection number of
two divisors by formula (4.6), if one of the two divisors is
concentrated in $Y$ (he also assumes $f$ to be proper). We will see
now that the general case essentially reduces to the case of
divisors.

Let $Z$ be a closed subscheme of $X$ and define $Z^{\roman{pure}}$
to be the closed subscheme of $Z$ defined by the ideal sheaf of
local sections with finite support. In particular, if $Z$ is
zero-dimensional, then $Z^{\roman{pure}}=\emptyset$.
\proclaim{Lemma 4.2}
$Z^{\roman{pure}}$ is Cohen-Macauley, i.e.\ has no embedded
components. If ${\roman{dim}}\, Z=1$, then $Z^{\roman{pure}}$ is a
divisor on $X$ (i.e.\ defined locally by one non--zero element).
\endproclaim
\demo{Proof}
The cases when ${\roman{dim}}\, Z=2$ or ${\roman{dim}}\, Z=0$ are
trivial. Now let ${\roman{dim}}\, Z=1$. Then the first assertion
implies the second (EGA IV, 21.7.2, 21.6.9, 21.11.1). The first
assertion is local around a point $y\in Y$. If $Z^{\roman{pure}}$
had an embedded prime ideal at $y$, it would have to be the maximal
ideal ${\frak m}_y$ of ${\Cal O}_y$. But then there would exist
$a\in {\Cal O}_y$ such that
$$\fm_y={\roman{rad}}(I:a)$$
(\cite{\AM}, th.\ 4.5). Here $I\subset {\Cal O}_y$ is the ideal of
$Z^{\roman{pure}}$. But then
$$\fm_y^n\cdot a\subset I\ \ ,\ \ \text{some}\ n\ \ ,$$
hence $a$ would define a local section of ${\Cal
O}_{Z^{\roman{pure}}}$ with support in $y$. By the definition of
$Z^{\roman{pure}}$ this implies $a\in I$, a contradiction.
\qed\enddemo
The next statement is analogous to Theorem 2.3.8, (iii) of Deligne
\cite{\deligne}.
\proclaim{Lemma 4.3}
Let $Z$ and $Z'$ be proper closed subschemes of $X$ such that their
intersection number (4.6) is defined. Then
$$(Z,Z')= (Z^{\roman{pure}}, Z^{'\roman{pure}})\ \ .\tag4.8$$
\endproclaim
\demo{Proof}
There are exact sequences
$$\align
0 &
\to {\Cal J}\to {\Cal O}_Z\to {\Cal O}_{Z^{\roman{pure}}}\to 0\\
0 &
\to {\Cal J}'\to {\Cal O}_{Z'}\to {\Cal O}_{Z^{'\roman{pure}}}\to 0\ \ ,\\
\endalign$$
where ${\Cal J}$ and ${\Cal J}'$ have finite support contained in
$Y$. Using Lemma 4.1 and the bilinearity of the tensor product and
the additivity of $\chi$, the result follows.
\qed\enddemo
One last fact we need is the following (Deligne \cite{\deligne}): Assume
that the morphism $f$ is proper and let $Z$ be a closed subscheme
with support in $Y$. Then
$$(Z,f^{-1}(s))=0\ \ .\tag4.9$$
\par\noindent
{\bf Remark 4.4.} The preceding theory applies equally well to the
case when $X$ is a formal scheme which is regular and where $f:X\to
S$ is an adic morphism into the formal spectrum of a complete
discrete valuation ring which is flat and locally of finite type
and with one-dimensional fibre over $s\in S$. In this case we have
defined the intersection number of closed formal subschemes $Z$ and
$Z'$ such that the sum of their defining ideals is open in ${\Cal
O}_X$, with the same properties as before. Furthermore, there is an
obvious compatibility between these notions: if the adic morphism
and the closed formal subschemes are formal completions along the
special fibre of a morphism of schemes and of closed subschemes
then both intersection numbers coincide.

We now return to the case of interest, namely to the formal scheme
${\Cal M}\simeq
\hat\Omega\times_{\roman{Spf}\,\Z_p}{\roman{Spf}}\, W$ over
${\roman{Spf}}\, W$. By Lemma 4.3 we need to know the divisors
$Z(j)^{\roman{pure}}$. The following proposition is an immediate
consequence of Propositions 3.2 and 3.3.
\proclaim{Proposition 4.5}
Let $q(j)=\varepsilon\cdot p^{\alpha}$, $\varepsilon\in
\Z_p^{\times}$, $\alpha\ge 0$.
\par\noindent
(i) If $\alpha$ is even and $\chi(\varepsilon)=1$, then
$$Z(j)^{\roman{pure}}=\sum_{[\Lambda]} {\roman{mult}}_{[\Lambda]} (j)
\cdot {\Bbb P}_{[\Lambda]}\tag4.10$$
(equality of divisors on ${\Cal M}$).
\par\noindent
(ii) If $\alpha$ is even and $\chi(\varepsilon)=-1$, then
$$Z(j)^{\roman{pure}}=Z(j)^h +\sum_{[\Lambda]} {\roman{mult}}_{[\Lambda]} (j)
\cdot {\Bbb P}_{[\Lambda]}\ \ ,\tag4.11$$
where $Z(j)^h$ (the horizontal part of $Z(j)$) is the disjoint sum
of two divisors each projecting isomorphically to ${\roman{Spf}}\,
W$ and meeting the special fibre in an ordinary special point of
${\Bbb P}_{[\Lambda (j)]}$. Here $[\Lambda(j)]$ denotes the unique
vertex fixed by $j$.
\par\noindent
(iii) If $\alpha$ is odd, then
$$Z(j)^{\roman{pure}}=Z(j)^h +\sum_{[\Lambda]} {\roman{mult}}_{[\Lambda]}
(j){\Bbb P}_{[\Lambda]}\tag4.12$$ where the divisor $Z(j)^h$ is the
formal spectrum of the ring of integers in a ramified quadratic
extension of $W_{\Q}$ which meets the special fibre in
$pt_{\Delta(j)}$, where $\Delta(j)$ is the edge containing the
unique fixed point of $j$.
\par\noindent
In particular, in all cases the divisor $Z(j)^{\roman{pure}}$ is
the sum with multiplicities of regular one-dimensional formal
schemes (or empty).
\endproclaim
\par\noindent
{\bf Remark 4.6.} Genestier \cite{\genestier} has formulated a moduli problem
(over ${\Cal M}$) whose solution is $Z(j)^{\roman{pure}}$. The
reason that we stick with the ``uglier'' subscheme $Z(j)$ is that
its definition can be transposed easily to other cases, \cite{\krsiegel},
\cite{\krHB}. Due to Lemma 4.3 the difference between $Z(j)$ and
$Z(j)^{\roman{pure}}$ has no consequences for the intersection
numbers.

Next we will have to determine the intersection numbers between all
the various summands of $Z(j)^{\roman{pure}}$ resp.\
$Z(j')^{\roman{pure}}$.
\proclaim{Lemma 4.7}
For any pair of vertices $[\Lambda]$, $[\Lambda']$ we have
$$({\Bbb P}_{[\Lambda]}, {\Bbb P}_{[\Lambda']})=\cases
1&\ \text{if}\ ([\Lambda], [\Lambda'])\ \text{is an edge}\\
-(p+1)
&\ \text{if}\ [\Lambda]=[\Lambda']\\ 0&\ \text{in all other
cases.}\\
\endcases$$
\endproclaim
\demo{Proof}
The first and the last entry on the RHS are obvious. Let
$\Gamma\subset PGL_2(\Q_p)$ be a cocompact discrete subgroup. Then
$\hat\Omega /\Gamma$ is a formal scheme which is {\it proper}\/ and
flat over ${\roman{Spf}}\, \Z_p$. If $\Gamma$ is sufficiently small
the projection
$$\hat\Omega\to \hat\Omega /\Gamma$$
is locally around ${\Bbb P}_{[\Lambda]}$ an isomorphism. The
special fibre $Y$ of $\hat\Omega /\Gamma$ is reduced and is of the
form
$$Y=\sum_{[\Lambda]\ {\roman{mod}}\, \Gamma} {\Bbb P}_{[\Lambda]}\ \ .$$
Using now (4.9) we obtain
$$\align
0 &
=({\Bbb P}_{[\Lambda]},Y)= ({\Bbb P}_{[\Lambda]}, (\sum_{[\Lambda']}
{\Bbb P}_{[\Lambda']}))\\ &
=({\Bbb P}_{[\Lambda]}, {\Bbb P}_{[\Lambda]}) +
\sum_{\matrix{\scriptstyle
[\Lambda']}\\ {\scriptstyle [\Lambda']\ne [\Lambda]}\\ \endmatrix}
({\Bbb P}_{[\Lambda]}, {\Bbb P}_{[\Lambda']})\ \ .\\
\endalign$$
But ${\Bbb P}_{[\Lambda]}$ meets precisely $p+1$ other irreducible
components and with multiplicity one. The result follows.
\qed\enddemo
\proclaim{Lemma 4.8}
Let $\j$ be a non-degenerate $\Z_p$-submodule of rank 2 and let
$j,j'$ be generators of $\j$ which diagonalize the bilinear form,
with matrix (2.18). For the horizontal part of the associated
divisors $Z(j)^{\roman{pure}}$ resp.\ $Z(j')^{\roman{pure}}$ we
have
$$(Z(j)^h,Z(j')^h)=\cases
1&\ \text{if}\ \alpha\ \text{and}\ \beta\ \text{odd}\\ 0 &\
\text{otherwise}\ \ .\endcases$$
\endproclaim
\demo{Proof}
This follows from Corollary 3.7.
\qed\enddemo
\proclaim{Lemma 4.9}
Let $q(j)=\varepsilon\cdot p^{\alpha}$ and let $[\Lambda]$ be a
vertex. Then
$$(Z(j)^h,{\Bbb P}_{[\Lambda]})=\cases
2&\ \text{if $\alpha$ is even and $\chi(\varepsilon)=-1$ and ${\Cal
B}^j =\{ [\Lambda]\}$}\\ 1&\ \text{if $\alpha$ is odd and
$d([\Lambda], {\Cal B}^j)=\frac{1}{2}$}\\ 0 &\ \text{in all other
cases.}\\
\endcases$$
\endproclaim
\demo{Proof}
This follows immediately from Propositions 3.2 and 3.3 and the
local equations for $Z(j)^h$ appearing in their proofs ((3.15)
resp.\ (3.24)).
\qed\enddemo

\subheading{\Sec 5. An invariance property of intersection numbers}
\par\noindent
The aim of this section is to establish the following invariance
property of the intersection numbers of special cycles.

\proclaim{Theorem 5.1}
Let $(j_1, j_2)$ resp.\ $(j'_1, j'_2)$ be two $\Z_p$-bases of the
same non-degenerate submodule $\j$ of $V$ with
$q(j_i)\in\Z_p\setminus \{ 0\}$ and $q(j'_i)\in\Z_p\setminus\{
0\}$, $i=1,2$. Then
$${\Cal O}_{Z(j_1)}\otimes^{\Bbb L} {\Cal O}_{Z(j_2)}\simeq
{\Cal O}_{Z(j'_1)} \otimes^{\Bbb L}{\Cal O}_{Z(j'_2)}$$ and hence
$$(Z(j_1), Z(j_2))= (Z(j'_1), Z(j'_2))\ \ .$$
\endproclaim

In the statement of the theorem it is obvious that the zero'th
cohomology sheaves of these complexes are isomorphic since both of
them are equal to the structure sheaf of $Z(\j)$. The statement is
also obvious locally around an isolated ordinary point of $Z(\j)$
since here the special cycles are divisors and hence there are no
higher Tor-terms if the intersection is zero-dimensional. But the
full statement is non-trivial. One basic ingredient of its proof
will be to show that the Genestier equations (3.5) and (3.8)
globalize to give a resolution of the structure sheaf of a special
cycle.

We therefore start with a single special cycle $Z(j)$. Let us
suppose that the matrix of $j$ in terms of the standard basis $e_1,
e_2$ is
$$j=\pmatrix
a&b\\ c&-a
\endpmatrix\ \ .$$
If $pt_{\Delta_0}\in Z(j)$, then $p\vert b$ and writing $b=pb_0$ we
have the global section of ${\Cal O}$ on $\hat{\Omega}_{\Delta_0}$,
$$f=b_0T_0-2a-cT_1\ \ .\tag5.1$$
Recall from (3.8) that the Genestier equations for the intersection
$\hat{\Omega}_{\Delta_0}\cap Z(j)$ are given as
$$f\cdot T_0=f\cdot T_1=0\ \ .$$
This leads us to consider the complex of free ${\Cal O}$-modules on
$\hat{\Omega}_{\Delta_0}$ concentrated in degrees $-2,-1,0$,
$$K(1,j):{\Cal O} \buildrel\scriptstyle\pmatrix \scriptstyle -T_1\\ \scriptstyle
T_0\endpmatrix \over \longrightarrow {\Cal O}\oplus {\Cal O}
\buildrel (f\cdot T_0, f\cdot T_1)\over \longrightarrow {\Cal O}\ \
.\tag5.2$$
We specify here that we consider an element in degree $-1$ as a
column vector.

By what we said above $K(1,j)$ is a free resolution of the
structure sheaf of $\hat{\Omega}_{\Delta_0}\cap Z(j)$. In order to
calculate with this complex, we write
$$f=f(j)= tr(j\cdot\tau\cdot{\roman{diag}}(-p^{-1}, 1))\tag5.3$$
with
$$\tau= \pmatrix p&-T_1\\ -T_0&1\endpmatrix = \pmatrix
T_0T_1&-T_1\\ -T_0&1\endpmatrix\ \ .\tag5.4$$ From now on we have
no further use of the matrix representation of $j$ and the symbols
$a,b,c$ will denote the entries of other matrices appearing in the
proofs below.

More generally, let $g\in G(\Q_p)^0$ such that $g\cdot
pt_{\Delta_0}\in Z(j)$. We have the homomorphism of sheaves
(covering the automorphism $g$)
$$g_*: {\Cal O}_{\hat{\Omega}}\longrightarrow {\Cal
O}_{\hat{\Omega}}\ \ .$$ The global sections $g_*(T_0)$, $g_*(T_1)$
of ${\Cal O}_{\hat{\Omega}_{g\Delta_0}}$ are just the global
sections constructed like $T_0, T_1$ but starting from the basis
$[ge_1, ge_2]$ instead of $[e_1, e_2]$. We then obtain a complex
$K(g,j)$ of ${\Cal O}$-modules on $\hat\Omega_{g\Delta_0}$ by using
in (5.2) instead of $T_0, T_1$ the global coordinates on
$\hat\Omega_{g\Delta_0}$ obtained from $T_0, T_1$ by transport via
$g_*$ and in (5.1) the coefficients of the matrix of $j$ in terms
of the basis $ge_1, ge_2$. We have an obvious isomorphism (the
identity)
$$g^*(K(g, gjg^{-1}))\longrightarrow K(1,j)\ \ ,$$
which we prefer to view as a $g_*$-linear homomorphism (between
sheaves on different spaces
$$g_*:K(1,j)\longrightarrow K(g,gjg^{-1})\ \ ,\tag5.5$$
i.e.\ $g_*(xs)=g_*(x)\cdot g_*(s)$ for a section $x$ of ${\Cal O}$
and a section $s$ of $K(1,j)$.
\proclaim{Lemma 5.2}
Assume $pt_{\Delta_0}\in Z(j)$. For every element $g$ in the
Iwahori subgroup (1.11) associated to $\Delta_0$ there is a
$g_*^{-1}$-linear isomorphism of complexes of ${\Cal O}$-modules on
$\hat\Omega_{\Delta_0}$
$$\lambda_g^0:K(1,j)\longrightarrow K(1,g^{-1}jg)\ \ .$$
These isomorphisms satisfy the transitivity condition
$$\lambda^0_{g_1 g_2}=\lambda^0_{g_2}\circ \lambda^0_{g_1}\ \ .$$
\endproclaim
\demo{Proof}
The isomorphism will be given in the following form
$$\matrix
K(1,j)&:&{\Cal O}&\buildrel\scriptstyle\pmatrix \scriptstyle -T_1\\
\scriptstyle T_0\endpmatrix\over\longrightarrow &{\Cal O}\oplus{\Cal O}
&\buildrel (fT_0, fT_1)\over\longrightarrow&{\Cal O}\\ {}\\
\llap{$\scriptstyle\lambda_g^0$}\Big\downarrow
&&
\llap{$\scriptstyle \alpha\cdot$}
\Big\downarrow&&\llap{$\scriptstyle\beta\cdot$}\Big\downarrow&&
\llap{$\scriptstyle 1\cdot$}\Big\downarrow\\
{}\\
 K(1,g^{-1}jg)&:&{\Cal O}&
\buildrel\scriptstyle\pmatrix\scriptstyle -T_1\\
\scriptstyle T_0\endpmatrix\over\longrightarrow &{\Cal
O}\oplus{\Cal O} &\buildrel (f'T_0, f'T_1)\over\longrightarrow
&{\Cal O}&,\endmatrix
\tag5.6$$
where
$$\alpha\in {\Cal O}^{\times}\ \ ,\ \ \beta={\roman{diag}}(\beta_0,
\beta_1)\in GL_2({\Cal O})\ \ .\tag5.7$$
Of course, we are using here the usual convention for describing a
semi-linear map by a matrix, e.g.\ in degree $-2$ a section $x$ of
${\Cal O}$ is mapped to $\alpha\cdot g_*^{-1}(x)$. Also we have
abbreviated $f(g^{-1}jg)$ into $f'$.
\hfill\break
Let us calculate $f'$. We have
$$\align
f' &
={\roman{tr}}(g^{-1}jg\cdot\tau\cdot{\roman{diag}}(-p^{-1},1))\tag5.8\\
&
={\roman{tr}}(j\cdot\tau'\cdot{\roman{diag}}(-p^{-1}, 1))\ \ ,
\endalign$$
with
$$\tau'=g\cdot\tau\cdot{\roman{diag}}(-p^{-1},1)\cdot
g^{-1}\cdot{\roman{diag}}(-p^{-1},1)^{-1}\ \ .\tag5.9$$ Now
$$\tau=\pmatrix T_1\\ -1\endpmatrix \cdot (T_0, -1)\ \ .$$
Let
$$g=\pmatrix a&b\\ c&d\endpmatrix =\pmatrix a&pb_0\\
c&d\endpmatrix\ \ .$$ Then
$$g\cdot \pmatrix T_1\\ -1\endpmatrix =\pmatrix aT_1&-b\\
cT_1&-d\endpmatrix = (-cT_1+d) \cdot \pmatrix g_*^{-1}(T_1)\\
-1\endpmatrix\ \ ,\tag5.10$$
where we used the expression (1.12) for $g_*^{-1}(T_1)$. We note
that
$$v=v(g, T_1)=-cT_1+d\tag5.11$$
is an automorphy factor for the action of the Iwahori subgroup on
${\Cal O}^{\times}$, i.e.\
$$v(gg', T_1)=v(g,{g'_*}^{-1}(T_1))\cdot v(g', T_1)\ \ .\tag5.12$$
Similarly we have
$$\align
 (T_0,-1)\cdot & {\roman{diag}}(-p^{-1},1)\cdot
g^{-1}\cdot{\roman{diag}}(-p^{-1},1)^{-1}\\ & =(T_0,
-1)\cdot{\roman{det}}(g)^{-1}\cdot\pmatrix d&b_0\\ pc&a\endpmatrix
\\
&
={\roman{det}}(g)^{-1}\cdot (dT_0-pc, b_0T_0-a)\\
&
={\roman{det}}(g)^{-1}\cdot (-b_0T_0+a)\cdot (g_*^{-1}(T_0), -1)\ \
.
\endalign$$
Again
$$u=u(g,T_0)=-b_0T_0+a\tag5.13$$
is an automorphy factor valued in ${\Cal O}^{\times}$. Plugging in
these expressions into (5.8) we therefore obtain
$$\tau'={\roman{det}}(g)^{-1}uv\cdot g_*^{-1}(\tau)\ \ ,\tag5.14$$
and
$$f'={\roman{det}}(g)^{-1}uv\cdot g_*^{-1}(f)\ \ .\tag5.15$$
We also note that
$$\frac{g_*^{-1}(T_0)}{T_0} =\frac{-cT_1+d}{-b_0T_0+a} =
\frac{v}{u} =\left( \frac{g_*^{-1}(T_1)}{T_1}\right)^{-1}\ \
.\tag5.16$$
We now fill in the diagram (5.6) by setting
$$\alpha= {\roman{det}}(g)\cdot u^{-1}v^{-1}\ \ ,\ \
\beta_0={\roman{det}}(g)\cdot u^{-2}\ \ ,\ \
\beta_1={\roman{det}}(g)\cdot v^{-2}\ \ .\tag5.17$$
We claim that the diagram commutes.
\par\noindent
{\bf Left square:} Via the NE passage the element 1 is mapped to
$$1\longmapsto \pmatrix -T_1\\ T_0\endpmatrix \longmapsto \pmatrix
-\beta_0\cdot g_*^{-1}(T_1)\\ \beta_1\cdot g_*^{-1}(T_0)\endpmatrix$$
Via the SW passage the element 1 is mapped to
$$1\longmapsto \alpha\longmapsto \alpha\cdot \pmatrix -T_1\\
T_0\endpmatrix\ \ .$$ The claim follows from the definitions (5.17)
and the identity (5.16).
\par\noindent
{\bf Right square:} The element $\pmatrix 1\\ 0\endpmatrix$ is
mapped via the NE passage to
$$\pmatrix 1\\ 0\endpmatrix \longmapsto f\cdot T_0\longmapsto
g_*^{-1}(f)\cdot g_*^{-1}(T_0)$$ and via the SW passage to
$$\pmatrix 1\\ 0\endpmatrix\longmapsto \beta_0 \cdot \pmatrix 1\\
0\endpmatrix\longmapsto \beta_0\cdot f'\cdot T_0\ \ .$$ The last
element is equal to by (5.15)
$$\beta_0\cdot {\roman{det}} (g)^{-1}uv\cdot g_*^{-1}(f)\cdot T_0=
\frac{v}{u}\cdot T_0\cdot g_*^{-1}(f)\ \ .$$
Again the relation (5.16) gives the claim. The element $\pmatrix
0\\ 1\endpmatrix$ is treated in a similar way. We therefore have
constructed the homomorphism $\lambda_g^0$ and it is immediate to
check that $\lambda_g^0$ is in fact an isomorphism. The
transitivity assertion is an immediate consequence of the fact that
$u$ and $v$ are automorphy factors.
\qed
\enddemo

We now use the system of isomorphisms constructed above to
construct a well-defined complex on each open chart
$\hat\Omega_{\Delta}$ of $\hat\Omega$. If $pt_{\Delta}\not\in Z(j)$
this complex is zero by definition. Next suppose that
$pt_{\Delta}\in Z(j)$ and let $g_1, g_2\in G(\Q_p)^0$ such that
$g_1\Delta_0=g_2\Delta_0=\Delta$. Then $g_2=g_1\cdot h$ where $h$
lies in the Iwahori subgroup for $\Delta_0$. We define a (linear)
isomorphism of complexes on $\hat\Omega_{\Delta}$
$$\lambda_{g_2,g_1}: K(g_1,j)\longrightarrow K(g_2,j)\tag5.18$$
by putting
$$\lambda_{g_2,g_1}= g_{2*}\circ \lambda_h^0\circ g_{1*}^{-1}\ \
.\tag5.19$$
The transitivity assertion of Lemma (5.2) yields
$$\lambda_{g_3, g_2}\circ \lambda_{g_2, g_1}=\lambda_{g_3, g_1}\tag5.20$$
for any $g_3$ with $g_3\Delta_0=\Delta$. This is the precise
meaning of the well-definedness of the complex on
$\hat\Omega_{\Delta}$.

Our next aim is to glue the complexes just constructed on common
overlaps of local charts.
\proclaim{Lemma 5.3} Assume $pt_{\Delta_0}\in Z(j)$.
For any $g\in GL_2(\Z_p)$ there is a $g_*^{-1}$-linear isomorphism
of complexes of ${\Cal O}$-modules on $\hat\Omega_{[\Lambda_0]}$,
$$\mu_g^0:K(1,j)_{\vert\hat\Omega_{[\Lambda_0]}}\buildrel\sim\over
\longrightarrow
K(1,g^{-1}jg)_{\vert \hat\Omega_{[\Lambda_0]}}\ \ .$$ These
isomorphisms satisfy the transitivity condition
$$\mu^0_{g_1 g_2} = \mu^0_{g_2}\circ \mu^0_{g_1}\ \ .$$
Furthermore, if $g$ lies in the Iwahori subgroup corresponding to
$\Delta_0$, then $\mu^0_g= \lambda^0_{g\vert
\hat\Omega_{[\Lambda_0]}}$.
\endproclaim
\demo{Proof}
Taking into account the definition (1.15) of the open immersion,
the complexes in question have the following form in terms of the
natural coordinate $T$ on $\hat\Omega_{[\Lambda_0]}$.
$$\matrix
K(1,j)_{\vert\hat\Omega_{[\Lambda_0]}}: & {\Cal O} &
\buildrel\scriptstyle \pmatrix \scriptstyle -T^{-1}\\
\scriptstyle pT\endpmatrix \over\longrightarrow
& {\Cal O}\oplus {\Cal O} &
\buildrel (f\cdot pT, f\cdot T^{-1})\over\longrightarrow
& {\Cal O}\\ {}\\
 \llap{$\scriptstyle\mu_g^0$}\Big\downarrow&\llap{$\scriptstyle
 \alpha\cdot$}\Big\downarrow &&
\llap{$\scriptstyle \beta\cdot$}\Big\downarrow
&& \llap{$\scriptstyle 1\cdot$}\Big\downarrow\\ {}\\
K(1,g^{-1}jg)_{\vert\hat\Omega_{[\Lambda_0]}}: & {\Cal O} &
\buildrel\scriptstyle \pmatrix \scriptstyle -T^{-1}\\ \scriptstyle pT
\endpmatrix\over\longrightarrow
& {\Cal O}\oplus {\Cal O} & \buildrel (f'\cdot pT, f'\cdot
T^{-1})\over\longrightarrow & {\Cal O}
\endmatrix\tag5.21$$
Here $f$ and $f'$ are the restrictions of the functions appearing
in (5.6) and our aim is to define $\alpha\in {\Cal O}^{\times}$ and
$\beta= {\roman{diag}}(\beta_0, \beta_1)\in GL_2({\Cal O})$ to
produce a commutative diagram. Let us calculate $f$ and $f'$. We
have
$$f={\roman{tr}}(j\cdot \pmatrix p&-T^{-1}\\ -pT&1\endpmatrix\cdot
{\roman{diag}}(-p^{-1}, 1))\ \ ,\ \ \text{i.e.,}$$
$$f={\roman{tr}}\left( j\cdot \pmatrix -1&-T^{-1}\\ T&1\endpmatrix
\right)\ \ .\tag5.22$$
Hence
$$f'= {\roman{tr}}(j\cdot g\cdot\pmatrix -1&-T^{-1}\\ T&1\endpmatrix
g^{-1})\ \ .$$ We write
$$\pmatrix -1&-T^{-1}\\ T&1\endpmatrix =\pmatrix -1\\ T\endpmatrix
(1, T^{-1})\ \ .$$ Let
$$g=\pmatrix a&b\\ c&d\endpmatrix\in GL_2(\Z_p)\ \  .$$
Then
$$g\cdot\pmatrix -1\\ T\endpmatrix =\pmatrix bT&-a\\
dT&-c\endpmatrix = (-bT+a)\pmatrix -1\\ g_*^{-1}(T)\endpmatrix\ \
,$$ where we use (1.8). We again note that
$$u=-bT+a\tag5.23$$
is an automorphy factor for the action of $GL_2(\Z_p)$ on ${\Cal
O}^{\times}$. Similarly we have
$$\align
(1,T^{-1})\cdot g^{-1} &
={\roman{det}} (g)^{-1}\cdot (d-cT^{-1}, b+aT^{-1})\\
&
={\roman{det}}(g)^{-1}\cdot (d-cT^{-1})\cdot (1, g_*^{-1}(T)^{-1})\
\ .\\
\endalign$$
Again
$$v=d-cT^{-1}\tag5.24$$
is an automorphy factor valued in ${\Cal O}^{\times}$. We thus
obtain
$$f'={\roman{det}}(g)^{-1}\cdot uv\cdot g_*^{-1}(f)\ \ .\tag5.25$$
We also note that
$$\frac{g_*^{-1}(T)}{T} = \frac{v}{u}\ \ .\tag5.26$$
We now fill in diagram (5.21) by setting
$$\alpha= {\roman{det}}(g)\cdot u^{-1}v^{-1}\ \ ,\ \
\beta_0={\roman{det}}(g)\cdot u^{-2}\ \ ,\ \ \beta_1={\roman{det}}
(g)\cdot v^{-2}\ \ .\tag5.27$$ The commutativity of the diagram is
checked as before. The transitivity property again follows from the
fact that $u$ and $v$ are automorphy factors. Finally, if $g$ lies
in the Iwahori subgroup the matrix entries $\alpha,
\beta_0,
\beta_1$ of $\mu_g^0$ are just the restrictions of the
corresponding entries of $\lambda_g^0$ since
$$\align
& u(g,T_0)_{\vert\hat\Omega_{[\Lambda_0]}}= (-b_0T_0+a)_{\vert
\hat\Omega_{[\Lambda_0]}} =-bT+a=u\\
& v(g,T_1)_{\vert\hat\Omega_{[\Lambda_0]}} =(-cT_1 +d)_{\vert
\hat\Omega_{[\Lambda_0]}} =d-cT^{-1}=v\ \ .\quad\qed
\endalign$$
\enddemo
\proclaim{Lemma 5.4}
Assume $pt_{\Delta_0}\in Z(j)$. Let
$$w=\pmatrix 0&p\\1&0\endpmatrix\ \ ,$$
so that $wGL_2(\Z_p)w^{-1}$ is the stabilizer of the lattice
$\Lambda_1=[pe_1, e_2]$. For any $g\in wGL_2(\Z_p)w^{-1}$ there is
a $g_*^{-1}$-linear isomorphism of complexes of ${\Cal O}$-modules
on $\hat\Omega_{[\Lambda_1]}$,
$$\mu_g^0:K(1,j)_{\vert \hat\Omega_{[\Lambda_1]}}\buildrel\sim\over\longrightarrow
K(1,g^{-1}jg)_{\vert \hat\Omega_{[\Lambda_1]}}\ \ .$$ These
isomorphisms satisfy the transitivity condition and if $g$ lies in
the Iwahori subgroup corresponding to $\Delta_0$ we have
$\mu_g^0=\lambda^0_{g\vert\hat\Omega_{[\Lambda_1]}}$.
\endproclaim
\demo{Proof}
The element $w$ lies in the normalizer of the Iwahori subgroup and
hence acts on $\hat\Omega_{\Delta_0}$. The action is given formally
by the same formula as in (1.12), hence interchanges $T_0$ and
$T_1$. Furthermore the action takes the open immersion of
$\hat\Omega_{[\Lambda_1]}$ in $\hat\Omega_{\Delta_0}$ into the open
immersion of $\hat\Omega_{[\Lambda_0]}$ in $\hat\Omega_{\Delta_0}$.
Hence after identifying $\hat\Omega_{[\Lambda_1]}$ with
$\hat\Omega_{[\Lambda_0]}$ via $w$, this open immersion is given by
$$T_0\longmapsto T^{-1}\ \ ,\ \ T_1\longmapsto pT\ \ .\tag5.28$$
For the restriction of $f$ we therefore obtain
$$f={\roman{tr}}(j\cdot \pmatrix p&-pT\\ -T^{-1}&1\endpmatrix\cdot
{\roman{diag}}(-p^{-1}, 1))\ \ .\tag5.29$$ We write
$$\pmatrix p&-pT\\ -T^{-1}&1\endpmatrix\cdot
{\roman{diag}}(-p^{-1}, 1)= \pmatrix -1\\ p^{-1}T^{-1}\endpmatrix
\cdot (1,pT)\ \ .$$ For the restriction of $f'$ we obtain
$$f'= {\roman{tr}} (j\cdot g\cdot \pmatrix -1\\
p^{-1}T^{-1}\endpmatrix \cdot (1,pT)\cdot g^{-1})\ \ .$$ Let
$g=\pmatrix a&b\\ c&d\endpmatrix \in wGL_2(\Z_p)w^{-1}$. Then
$$g\cdot\pmatrix -1\\p^{-1}T^{-1}\endpmatrix = \pmatrix
-a+p^{-1}bT^{-1}\\ -c+p^{-1}dT^{-1}\endpmatrix =
(a-p^{-1}bT^{-1})\cdot \pmatrix -1\\
\frac{-c+p^{-1}dT^{-1}}{a-p^{-1}bT^{-1}}\endpmatrix\ \ .\tag5.30$$
The second component of the last vector is equal to
$$p^{-1}\cdot\left( \frac{aT-p^{-1}b}{-pcT+d}\right)^{-1} =
p^{-1}{g'_*}^{-1}(T)^{-1}\ \ .\tag5.31$$ Here
$$g'=wgw^{-1}=\pmatrix d&pc\\ p^{-1}b&a\endpmatrix\in GL_2(\Z_p)\ \
.\tag5.32$$ Similarly,
$$(1,pT)\cdot g^{-1}={\roman{det}}(g)^{-1}\cdot (d-pcT)\cdot (1,p\cdot
{g'_*}^{-1}(T))\ \ .\tag5.33$$ Putting
$$u=a-p^{-1}b\cdot T^{-1}\ \ ,\ \ v=d-pcT\tag5.34$$
we therefore obtain
$$f'= {\roman{det}}(g)^{-1}uv\cdot {g'_*}^{-1}(f)\ \ .\tag5.35$$
Of course, after reversing the identification of
$\hat\Omega_{[\Lambda_0]}$ with $\hat\Omega_{[\Lambda_1]}$, the
isomorphism ${g'_*}^{-1}$ becomes $g_*^{-1}$. We define
$$\alpha ={\roman{det}}(g)^{-1}\cdot u^{-1}v^{-1},\
\beta_0={\roman{det}}(g)\cdot u^{-2},\
\beta_1={\roman{det}}(g)\cdot v^{-2}\tag5.36$$
and the proof proceeds as before.
\qed
\enddemo

We now use the previous lemmas to glue the complexes on the open
charts $\hat\Omega_{\Delta}$ with $pt_{\Delta}\in Z(j)$. Let $g_1,
g_2\in G(\Q_p)^0$ such that $g_i\cdot pt_{\Delta_0}\in Z(j)$ for
$i=1,2$ and with $g_1\Delta_0\cap g_2\Delta_0=\{ [\Lambda]\}$. We
define an isomorphism of complexes on $\hat\Omega_{[\Lambda]}$
$$\mu_{g_2, g_1}: K(g_1, j)_{\vert\hat\Omega_{[\Lambda]}}
\longrightarrow K(g_2, j)_{\vert
\hat\Omega_{[\Lambda]}}\tag5.37$$ by putting
$$\mu_{g_2, g_1}=g_{2*}\circ \mu^0_h\circ g_{1*}^{-1}\ \
.\tag5.38$$
Here $g_2=g_1\cdot h$ where $h\in GL_2(\Z_p)$ if $[\Lambda]$ is in
the $G(\Q_p)^0$-orbit of $[\Lambda_0]$ and $h\in wGL_2(\Z_p)w^{-1}$
if $[\Lambda]$ is in the $G(\Q_p)^0$-orbit of $[\Lambda_1]$. The
isomorphism $\mu_h^0$ has been defined in Lemma (5.3) resp.\ in
Lemma (5.4). Again there is a transitivity condition. If $g_3$ is a
third element with $g_3\cdot pt_{\Delta_0}\in Z(j)$ and with
$g_3\Delta_0\cap g_i\Delta_0=\{ [\Lambda]\}$, $i=1,2$, then
$$\mu_{g_3, g_1}= \mu_{g_3, g_2}\circ \mu_{g_2, g_1}\ \ .\tag5.39$$

Having constructed a global resolution of ${\Cal O}_{Z(j)}$ on the
complement of the isolated ordinary points of $Z(j)$ we now
consider tensor products of two of them. We let $Z(j_1), Z(j_2)$
and $Z(j'_1), Z(j'_2)$ be two pairs of special cycles. We assume
that
$$(j_1, j_2)\cdot\gamma =(j'_1, j'_2)\tag5.40$$
for
$$\gamma = \pmatrix a&b\\ c&d\endpmatrix \in M_2(\Z_p)\ \
.\tag5.41$$
Let us first assume that $pt_{\Delta_0}\in Z(j_1)\cap Z(j_2)$ and
$pt_{\Delta_0}\in Z(j'_1)\cap Z(j'_2)$. We will denote by
$K(1,\underline j)$ the tensor product $K(1,j_1)\otimes_{\Cal O}
K(1,j_2)$ with its bigrading, and similarly for $K(1, \underline
j')$. We denote the differentials in degree $-1$ of $K(1,j_i)$ by
$\partial_i$ and those of $K(1,j'_i)$ by $\partial'_i$. All
differentials in degree $-2$ are identical and will be denoted by
$\partial$.
\proclaim{Lemma 5.5}
There is a homomorphism of complexes
$$\varphi: K(1,\underline j')\longrightarrow K(1,\underline j)$$
which is given degree by degree as follows.
\par\noindent
{\bf in degree$0$:}
$$\varphi^0={\roman{id}}_{\Cal O}:{\Cal
O}\longrightarrow {\Cal O}$$
\par\noindent
{\bf in degree $-1 =(-1, 0)+(0, -1)$:}
$$\matrix
\varphi^{-1}:
& ({\Cal O}\oplus {\Cal O})\oplus ({\Cal O}\oplus {\Cal O}) &
\longrightarrow
& ({\Cal O}\oplus {\Cal O}) \oplus ({\Cal O}\oplus {\Cal O})\\ &
(x,y) &
\longmapsto
& (x,y)\cdot {}^t\gamma
\endmatrix$$
(here, as before, $x$ and $y$ are considered as column vectors).
\par\noindent
{\bf in bi-degree $(-2, 0)+(0,-2)$:}
$$\matrix
\varphi^{-2}:
& {\Cal O}\oplus {\Cal O} &
\longrightarrow
& {\Cal O}\oplus {\Cal O} \\ & (x,y) &
\longmapsto
& (x,y)\cdot {}^t\gamma
\endmatrix$$
\par\noindent
{\bf in bi-degree $(-1,-1)$:}
$$\matrix
\varphi^{-2}:
& ({\Cal O}\oplus {\Cal O})\otimes ({\Cal O}\oplus {\Cal O}) &
\longrightarrow
& {\Cal O}\oplus ({\Cal O}\oplus {\Cal O})\otimes ({\Cal O} \oplus
{\Cal O}) \oplus {\Cal O}
\\ & (x\otimes y) &
\longmapsto
&(ab\psi_1(x,y), ad(x\otimes y)-bc(y\otimes x),\,  cd\psi_2(x,y))
\endmatrix$$
Here $\psi_1(x,y)\in {\Cal O}$ is the unique element of bidegree
$(-2,0)$ whose image under the differential $\partial$ of
$K(1,j_1)$ is
$$\partial\psi_1(x,y)=\partial_1(x)y-\partial_1(y)x\ \
.\tag5.42$$
The element $\psi_2(x,y)$ of bidegree $(0,-2)$ is similarly
defined.
\par\noindent
{\bf in degree $-3=(-2, -1)+(-1,-2)$:}
$$\matrix
\varphi^{-3}:
& ({\Cal O}\oplus {\Cal O}) \oplus ({\Cal O}\oplus {\Cal O}) &
\longrightarrow
& ({\Cal O}\oplus {\Cal O}) \oplus ({\Cal O} \oplus {\Cal O})\ \
.\\ &(x,y) &
\longmapsto
& (adx+bcy,ady+bcx)
\endmatrix$$
\par\noindent
{\bf in degree $-4=(-2,-2)$:}
$$
\varphi^{-4}=(ad+bc){\roman{id}}: {\Cal O}\longrightarrow {\Cal O}\ \ .$$
\endproclaim
\demo{Proof}
We have to check that $\varphi$ is indeed a homomorphism of
complexes. We picture the complexes written horizontally and the
homomorphism $\varphi$ vertically from the north to the south.
\par\noindent
{\bf in degree $-1$:} Via the NE-passage we have
$$(x,y)\longmapsto \partial'_1(x)+ \partial'_2(y)\longmapsto
\partial'_1(x)+\partial'_2(y)\ \ .$$
However, and this is the crucial observation for the proof, the
differential in degree $-1$ of $K(1,j)$ depends linearly on $j$.
Hence
$$\partial'_1= a\partial_1+c\partial_2\ \ ,\ \
\partial'_2=b\partial_1+d\partial_2\ \ .\tag5.43$$
Therefore the last expression is equal to
$$\partial_1(ax+by)+\partial_2(cx+dy)=\partial^{-1}(\varphi^{-1}(x,y))\
\ ,$$
where $\partial^{-1}$ is the differential in degree $-1$ of
$K(1,\underline j)$.
\par\noindent
{\bf in degree $-2$:} The argument in degree $(-2,0)+(0,-2)$ is
identical. Let us consider an element $x\otimes y$ of bi-degree
$(-1,
-1)$. It is mapped via the NE route to
$$x\otimes y\mapsto \partial'_1(x)y-\partial'_2(y) x\mapsto
(-a\partial'_2(y) x+b\partial'_1(x)y,
-c\partial'_2(y)x+d\partial'_1(x)y)\ \ .$$
Inserting the expression (5.43) for $\partial'_1$ and $\partial'_2$
and collecting terms, this last element is equal to
$$(ab\partial(\psi_1(x,y)) +bc\partial_2(x)y -ad\partial_2(y)x, cd
\partial(\psi_2(x,y))+ ad\partial_1(x)y- bc\partial_1(y)x)\ .$$
Via the SW route the element goes to
$$\align
& x\otimes y\mapsto (ab\psi_1(x,y), ad(x\otimes y)- bc(y\otimes x),
cd\psi_2(x,y))\mapsto\\ &
\mapsto (ab\partial(\psi_1(x,y))-ad\partial_2(y)x
+bc\partial_2(x)y,\\ &
\quad\quad
cd\psi_2(x,y)+ad \partial_1(x)y- bc\partial_1(y)x),
\endalign$$
which proves the claim in this case.
\par\noindent
{\bf in degree $-3=(-2,-1)+(-1,-2)$:} Let $(x,0)$ be of bi-degree
$(-2,-1)$. Then via the NE route this element goes to
$$\align
& (x,0)\mapsto (\partial'_2(x), \partial(1)\otimes x, 0)\mapsto\\ &
\mapsto (a\partial'_2(x)+ ab\psi_1(\partial(1), x),
ad(\partial(1)\otimes x)-bc(x\otimes \partial(1)),
c\partial'_2(x)+cd\psi_2(\partial(1),x))\ .
\endalign$$
The entry of bi-degree $(-2,0)$ is equal to
$$ab(\partial_1 (x)+\psi_1(\partial(1), x))+ad\partial_2(x)= ad
\partial_2(x)\ \ .$$
Here we used that the first summand vanishes since
$$\partial(\partial_1(x)+ \psi_1(\partial(1), x))=
\partial_1(x)\cdot\partial (1)+ \partial^2(1)x
-\partial(1)\partial_1(x) = 0\ .\tag5.44$$
The same reasoning shows that the entry of bi-degree $(0,-2)$ is
equal to $bc\partial_1(x)$. On the other hand, the image of $(x,0)$
via the SW route is equal to
$$(x,0)\mapsto (adx, bdx)\mapsto (ad\partial_2(x),
ad(\partial(1)\otimes x)- bc(x\otimes
\partial(1)),bc\partial_1(x)),$$ which shows the claim for elements
of bi-degree $(-2, -1)$. The case of elements of bi-degree
$(-1,-2)$ is analogous.
\par\noindent
{\bf in degree $-4$:} Here the generator $1\in{\Cal O}$ is sent via
the NE route to
$$1\mapsto (\partial(1), \partial(1))\mapsto
(ad\partial(1)+bc\partial(1), ad\partial(1)+ bc\partial(1))\ .$$
Via the SW passage the same element is sent to
$$1\mapsto (ad+bc)\cdot 1\mapsto (ad+bc)(\partial (1),
\partial(1))\ \ ,$$
which proves the claim in this case.
\qed
\enddemo
We note that the formulas for $\varphi$ do not involve the special
endomorphisms $j_1,j_2,j'_1, j'_2$ and neither do the formulas for
the glueing maps $\lambda_g^0, \mu_g^0$. This observation is
crucial for the proof of the next lemma.
\proclaim{Lemma 5.6}
a) Let $g$ be an element of the Iwahori subgroup corresponding to
$\Delta_0$. The following diagram is commutative.
$$\matrix
K(1,\underline j') &
\buildrel\varphi\over\longrightarrow
& K(1,\underline j)\\ {}\\
\llap{$\scriptstyle \lambda_g^{'0}\otimes \lambda_g^{'0}$}\big\downarrow
&&
\big\downarrow\rlap{$\scriptstyle \lambda_g^0\otimes \lambda_g^0$}\\
{}\\ K(1,g^{-1}\underline j'g) &
\buildrel\varphi\over\longrightarrow
& K(1,g^{-1}\underline jg)
\endmatrix$$
b) Let $g\in GL_2(\Z_p)$. The following diagram of complexes on
$\hat\Omega_{[\Lambda_0]}$ is commutative.
$$\matrix
K(1,\underline j')_{\vert \hat\Omega_{[\Lambda_0]}} &
\buildrel\varphi\over\longrightarrow
& K(1,\underline j)_{\vert\hat\Omega_{[\Lambda_0]}}\\ {}\\
\llap{$\scriptstyle \mu_g^{'0}\otimes \mu_g^{'0}$}\big\downarrow
&&
\big\downarrow\rlap{$\scriptstyle  \mu_g^0\otimes \mu_g^0$}\\
{}\\  K(1,g^{-1}\underline j'g)_{\vert \hat\Omega_{[\Lambda_0]}} &
\buildrel\varphi\over\longrightarrow
& K(1,g^{-1}\underline j g)_{\vert \hat\Omega_{[\Lambda_0]}}
\endmatrix$$
The analogous statement holds for $g\in wGL_2(\Z_p)w^{-1}$ and for
the restrictions to $\hat\Omega_{[\Lambda_1]}$.
\endproclaim
\demo{Proof} a) We again check this degree by degree. In degree 0
all homomorphisms are the identity, hence the assertion is obvious.
\par\noindent
{\bf degree $-1= (-1,0)+(0, -1)$:} Via the NE passage an element
$(x,y)$ is sent to
$$(x,y)\mapsto (x,y)\cdot {}^t\gamma \mapsto \beta\cdot (x',
y')\cdot {}^t\gamma\ \ .$$ Here for brevity we have set
$x'=g_*^{-1}(x)$. The result via the SW passage is obviously the
same.
\par\noindent
{\bf degree $(-2,0)+(0,-2)$:} Same argument as before.
\par\noindent
{\bf bi-degree $(-1,-1)$:} Via the NE passage an element $x\otimes
y$ is sent to
$$\align
& x\otimes y\mapsto (ab\psi_1(x,y)), ad(x\otimes y)-bc(y\otimes x),
cd\psi_2(x,y))\\ &
\mapsto (\alpha\cdot ab\psi_1(x,y)',
ad(\beta\otimes\beta)(x'\otimes y')-bc(\beta\otimes\beta)(y'\otimes
x'), \alpha\cdot cd\psi_2(x,y)').
\endalign$$
Via the SW passage the element is sent to
$$x\otimes y\mapsto (\beta x', \beta y')\mapsto (ab\psi_1(\beta
x', \beta y'), ad(\beta x'\otimes \beta y')
-bc (\beta y'\otimes \beta x'), cd \psi_2(\beta x', \beta y')).$$
Obviously the terms of degree $(-1,-1)$ coincide. Let us consider
the terms of degree $(-2,0)$, leaving aside the common factor $ab$.
But
$$\align
\partial(\alpha\cdot \psi_1(x,y)')
&
=\beta\cdot\partial(\psi_1(x,y))'\\
&
=\beta\cdot (\partial'_1(x)y-\partial'_2(y)x)'\\
&
=\partial'_1(x)'\cdot\beta y'-\partial'_2(y)'\beta x'\\
&
=\partial_1(\beta x') \beta y'-\partial_2(\beta y')\beta x\\
&
=\partial (\psi_1(\beta x', \beta y')),
\endalign$$
where we used twice that $\lambda_g^0$ is a homomorphism of
complexes. It follows that the terms of degree $(-2,0)$ are
identical. The same argument applies to the components of degree
$(0, -2)$.
\par\noindent
{\bf degree $(-2,-1)+(-1,-2)$:} Via the NE passage the element
$(x,0)$ of degree $(-2,-1)$ goes to
$$(x,0)\mapsto (ad\, x, bc\, x)\mapsto (ad\alpha\cdot\beta x',
bc\alpha\cdot\beta x)\ \ .$$ Via the SW route this element is sent
to
$$(x,0)\mapsto (\alpha\cdot\beta x',0)\mapsto (ad\alpha\cdot\beta x',
bc\alpha\cdot\beta x')\ \ ,$$ which proves the claim in this case.
Elements $(0,y)$ of degree $(-1, -2)$ are analogous.
\par\noindent
{\bf degree $(-2,-2)$:} The element 1 is sent via the NE passage to
$$1\mapsto (ad+bc)\cdot 1\mapsto (ad+bc)\cdot\alpha$$
and via the SW passage to
$$1\mapsto\alpha\mapsto (ad+bc)\alpha\ \ .$$
This concludes the proof of a). The proof of b) is formally
identical.
\qed
\enddemo

We now drop the assumption that $pt_{\Delta_0}\in Z(\underline
j)\cap Z(\underline j')$. Let $\Delta$ be an edge with
$pt_{\Delta}\in Z(j)\cap Z(j')$. Let $g\in G(\Q_p)^0$ with $\Delta
=g\Delta_0$. Then, just as we used the isomorphisms
$\lambda_h^0$ (for $h$ lying in the Iwahori subgroup) to construct
well-defined complexes on the open chart $\hat\Omega_{\Delta}$, we
use part a) of the previous lemma to define homomorphisms of
complexes on $\hat\Omega_{\Delta}$,
$$\varphi_g: K(g,\underline j')\longrightarrow K(g, \underline j)\tag5.45$$
by putting
$$\varphi_g=g_*\circ\varphi\circ g_*^{-1}\ \ .\tag5.46$$
If $g_1\Delta_0=g_2\Delta_0=\Delta$, then by part a) of the lemma
we have
$$\varphi_{g_2}\circ \lambda'_{g_2, g_1}=\lambda_{g_2,g_1}\circ
\varphi_{g_1}: K(g, \underline j')\longrightarrow K(g_2,
\underline j)\ \ ,\tag5.47$$
which means that the homomorphisms on the local chart
$\hat\Omega_{\Delta}$ are well-defined. By part b) of the previous
lemma the homomorphisms just constructed glue on the overlaps of
the open charts in question.
\proclaim{Lemma 5.7}
Suppose with the above notations that $\gamma\in GL_2(\Z_p)$ and
that $ad+bc$ is a unit, i.e., equivalently, that
$(ad)^2-(bc)^2\in\Z_p^{\times}$. Then the homomorphism of complexes
constructed above is an isomorphism.
\endproclaim
\demo{Proof}
It suffices to check this for $\varphi: K(1,\underline j')\to
K(1,\underline j)$, i.e., to show that the determinant of $\varphi$
is a unit in each degree.
\par\noindent
{\bf in degree $0$:} ${\roman{det}}(\varphi^0)=1$
\par\noindent
{\bf in degree $-1$:} ${\roman{det}}(\varphi^{-1})=
{\roman{det}}(\gamma)$
\par\noindent
{\bf in degree $-2$:}
${\roman{det}}(\varphi^2)={\roman{det}}(\gamma)^2\cdot((ad)^2-(bc)^2)$
\par\noindent
{\bf in degree $-3$:} ${\roman{det}}(\varphi^{-3})=
{\roman{det}}(\gamma)^2\cdot(ad+bc)^2$
\par\noindent
{\bf in degree $-4$:} ${\roman{det}}(\varphi^{-4})= ad+bc$.
\qed
\enddemo
\par\noindent
\demo{Proof of Theorem 5.1} By hypothesis we have the relation
(5.40), where $\gamma\in GL_2(\Z_p)$. If $ad+bc\in\Z_p^{\times}$
the result follows from the previous lemma. Otherwise we write
$\gamma
=\gamma_1\ldots\gamma_r$ as a product of elements which satisfy the
hypothesis of this lemma and obtain the desired isomorphism as the
composition of the isomorphisms of the previous lemma
$$K(1,\underline j')=K(1,\underline j\cdot\gamma)\buildrel
\varphi_r\over\longrightarrow K(1,\underline j\cdot \gamma_1\ldots
\gamma_{r-1})\longrightarrow\ldots \buildrel
\varphi_1\over\longrightarrow K(1,\underline j)\quad\qed$$
\enddemo
\par\noindent
{\bf Example 5.8.} We give an example to illustrate that in the
situation of Theorem 5.1 the relative positions of the cycles
$Z(j_1)$ and $Z(j_2)$ resp.\ of $Z(j'_1)$ and $Z(j'_2)$ can differ
radically. Assume that $j_1, j_2$ is a diagonal basis with
$$q(j_1)=\varepsilon_1\ \ \text{(thus}\ \alpha=0),\ \ \hbox{and}\ q(j_2)=
\varepsilon_2\cdot p^{\beta}\ \ ,\tag5.48$$
with $\beta >0$ even and $\chi(\varepsilon_1)=-1$. Then $Z(j_1)$ is
simply a horizontal divisor isomorphic to the disjoint sum of two
copies of ${\roman{Spf}}\, W$ meeting the special fibre in two
ordinary special points $x_1, x_2$ of ${\Bbb P}^1_{[\Lambda]}$,
where $[\Lambda]$ is the unique vertex fixed by $j_1$. This vertex
is also fixed by $j_2$ and the component ${\Bbb P}^1_{[\Lambda]}$
occurs with multiplicity $\beta/2$ in $Z(j_2)$. In particular
$(Z(\j))_{\roman{red}}=\{ x_1, x_2\}$ and
$$(Z(j_1), Z(j_2))= \beta\ \ .\tag5.49$$
As a new basis of $\j$ let us take
$$j_1,\ \hbox{and}\ j'_2 =j_2+\lambda j_1\ \ ,\ \ \lambda\in\Z_p^{\times}\ \ .$$
Then
$$q(j'_2)= \lambda^2\varepsilon_1+\varepsilon_2\cdot p^{\beta}=
\lambda^2\varepsilon_1(1+\lambda^{-2}\varepsilon_1^{-1}
\varepsilon_2\cdot p^{\beta})\in
\Z_p^{\times}-\Z_p^{\times,2}\ \ .$$
Hence $q(j'_2)= \varepsilon'_2$ with $\chi(\varepsilon'_2)=-1$ and
thus $Z(j'_2)$ is, just as $Z(j_1)$, a horizontal divisor
isomorphic to the disjoint sum of two copies of ${\roman{Spf}}\, W$
and meeting the special fibre in $x_1$ and $x_2$.

Let us calculate the intersection number $(Z(j_1), Z(j'_2))$ by
using the local equations. Since $Z(j_1)$ and $Z(j'_2)$ are coprime
divisors, we only have to determine the lengths of the local rings
at $x_1$ and $x_2$. Let us choose a basis for $\Lambda$ such that
$$j_1=\pmatrix 0&1\\ \varepsilon_1&0\endpmatrix\ \ ,\ \ j_2=
\pmatrix a&b\\ -\varepsilon_1b&-a\endpmatrix\ \ .$$
Then $j'_2= \pmatrix a&b+\lambda\\
\varepsilon_1(\lambda-b)&-a\endpmatrix$ and the local equations for
$Z(j_1)$ and $Z(j'_2)$ are respectively
$$T^2-\varepsilon_1=0\ \ ,\ \
(b+\lambda)T^2-2aT-\varepsilon_1(\lambda-b)=0\ \ .$$ The affine
ring of $Z(j_1)\cap Z(j'_2)$ is therefore
$$W[T, (T^p-T)^{-1}]^{\wedge}/ (2aT-2\varepsilon_1 b,
T^2-\varepsilon_1)$$ which is isomorphic to
$$W/(a-\eta b)\oplus W/(a+\eta b)$$
where $\eta$ is a square root of $\varepsilon_1$ in $W$. Therefore
$$\align
(Z(j_1), Z(j'_2)) &
= {\roman{ord}}_p(a-\eta b)+
{\roman{ord}}_p(a+\eta b)\\ &
={\roman{ord}}_p(a^2-\varepsilon_1b^2)=
{\roman{ord}}_pq(j_2)=\beta\ \ ,
\endalign$$
in accordance with (5.49). Of course, in this case the assertion of
Theorem 5.1 is trivial since
$${\Cal O}_{Z(j_1)} \otimes^{\Bbb L}{\Cal O}_{Z(j_2)} = {\Cal
O}_{Z(\j)}= {\Cal O}_{Z(j_1)}\otimes^{\Bbb L} {\Cal O}_{Z(j'_2)}$$
is concentrated in degree zero.

\subheading{\Sec6. Computation of intersection numbers}

In this section, we combine the information obtained so far and
give an explicit expression for the intersection number
$(Z(j),Z(j'))$.

\proclaim{Theorem 6.1} Let $j$ and $j'$ be special endomorphisms
such that their $\Z_p$-span $\j = \Z_p j + \Z_p j'$ is of rank $2$
and nondegenerate for the quadratic form. Let
$$T= \pmatrix q(j)&\frac12(j,j')\\\frac12(j',j)&q(j')\endpmatrix,$$
and suppose that $T$ is $GL_2(\Z_p)$-equivalent to
$\text{diag}(\e_1 p^\a,\e_2 p^\b)$, where $0\le \a\le \b$. Then
$$\align
(Z(j),Z(j')) &= e_p(T)=\\
\nass
{}& =\a+\b+1 - \cases p^{\a/2}+2\, \frac{\scr p^{\a/2}-1}{\scr p-1}
&\text{if $\a$ is even and $\chi(\e_1)=-1$,}\\
\nass
(\b-\a+1)p^{\a/2} + 2\,\frac{p^{\a/2}-1}{p-1}&
\text{if $\a$ is even and $\chi(\e_1)=1$,}\\
\nass
 2\,\frac{p^{(\a+1)/2}-1}{p-1}&\text{if $\a$ is odd.}
\endcases
\endalign
$$
\endproclaim

By Theorem~5.1, we may assume that the elements $j$ and $j'$ diagonalize
the quadratic form, i.e., that
$$T =\pmatrix \e_1 p^\a&{}\\{}&\e_2 p^\b\endpmatrix,$$
with $0\le \a\le \b$.

Recall that
$$\align
(Z(j),Z(j')) &= (Z(j)^{\text{pure}},Z(j')^{\text{pure}})\\
\nass
{}&=(Z(j)^h,Z(j')^h) + (Z(j)^h,Z(j')^v) + (Z(j)^v,Z(j')^h) +(Z(j)^v,Z(j')^v),
\endalign
$$
where
$$Z(j)^v = \sum_{[\Lambda]} \mult_{[\Lambda]}(j)\cdot \PP_{[\Lambda]},$$
as in Proposition~4.5. Here we assume that $q(j) = \e_1 p^\a$,
$q(j') = \e_2 p^\b$ and that $jj'=-j'j$. Recall that
$$\mult_{[\Lambda]}(j) = \max\{\ \a/2 - d([\Lambda],\Cal B^j),\ 0 \}.$$

We organize the calculation according to the cases in table (2.26),
assuming from now on that $\a\le \b$.

First consider the quantity $(Z(j)^h,Z(j')^h) + (Z(j)^h,Z(j')^v) +
(Z(j)^v,Z(j')^h)$. Recall, for example, that in the case $\a$ is
even and $\chi(\e_1)=-1$, $Z(j)^h$ consists of two copies of
$\text{\rm Spf}(W)$ meeting $\PP_{[\Lambda(j)]}$, while, if $\a$ is
odd, then $Z(j)^h$ consists of one copy of $\text{Spf}(W')$ meeting
$\text{pt}_{\Delta(j)}$, where $\Delta(j)$ is the edge containing
the fixed point of $j$. Otherwise $Z(j)^h$ is empty. Taking into
account the multiplicites of the vertical components and Lemmas~4.8
and~4.9, we obtain the following table of values of
$(Z(j)^h,Z(j')^h) + (Z(j)^h,Z(j')^v) + (Z(j)^v,Z(j')^h)$:
$$
\vcenter{
\begintab{|c|c|c|c|}
\hline
$j\setminus j'$&$\beta$ even&$\beta$ even&$\beta$ odd
\cr
&$\chi(\varepsilon_2)=-1$ & $\chi(\varepsilon_2)=1$ &
\cr
\hline
$\alpha$ even &&&
\cr
$\chi(\varepsilon_1)=-1$&$\a+\b$&$\b$&$\emptyset$
\cr
\hline
$\alpha$ even &&&
\cr
$\chi(\varepsilon_1)=1$ &$\a$&$0$&$\a$
\cr
\hline
&&&
\cr
&&&
\cr
&&&
\cr
 $\alpha$ odd & $\emptyset$ &$\b$&$\a+\b-1$
\cr
\hline
\endtab
}
\tag6.1
$$

It then remains to calculate $(Z(j)^v,Z(j')^v)$, using Lemma~4.7.
The first step is the following simple observation.
\proclaim{Lemma 6.2} For any vertex $[\Lambda]\in \Cal B$,
let $r = d([\Lambda],\Cal B^{j'})$ be the distance to the fixed
point set of $j'$. Note that $[\Lambda]$ lies on the boundary of
$\Cal T(j')$ precisely when $r=\b/2$. Then
$$(\PP_{[\Lambda]},Z(j')^v) = \cases
1-p &\text{ when $1\le r \le \b/2-1$,}\\
\nass
\chi(\e_2)-p &\text{ when $r=0$ and $\b$ is even,}\\
\nass
-p&\text{ when $r=1/2$ and $\b$ is odd,}\\
\nass
1&\text{ when $r=\b/2$,}\\
\nass
0&\text{ otherwise.}
\endcases$$
\endproclaim
\demo{Proof}
First suppose that $1\le r\le \b/2-1$, so that there is a unique edge
leading from $[\Lambda]$ along the geodesic from $[\Lambda]$ to $\Cal B^{j'}$,
and there are $p$ edges at $[\Lambda]$ leading away from $\Cal B^{j'}$.
The intersection of $\PP_{[\Lambda]}$ with $Z(j')^v$ is thus
$$ (\b/2-r+1) - (p+1)(\b/2-r) + p (\b/2-r-1) = 1-p,$$
in this case. The other cases are similar. For example, if $r=1/2$ and $\b$ is odd,
then $[\Lambda]$ is one endpoint of the edge containing the unique fixed point of
$j'$. The multiplicity in $Z(j')^v$
of the component corresponding to each these two endpoints is $(\b-1)/2$.
Thus, the intersection of $\PP_{[\Lambda]}$ with $Z(j')^v$ is
$$(\b-1)/2 -(p+1)(\b-1)/2 + p (\b-3)/2 = -p.$$
Finally, note that if $r=\b/2$, then the multiplicity of
$\PP_{[\Lambda]}$ in $Z(j')$ is zero, while the multiplicity of the
component corresponding to the unique neighboring vertex closer to
$\Cal  B^{j'}$ is $1$.
\qed\enddemo

\demo{Proof of Theorem~6.1} In all calculations, we will count, with multiplicity in
$Z(j)$, the number of
vertices $[\Lambda]$ in $\Cal T(j)$ whose associated component has a fixed
intersection number with $Z(j')$.  Recall that we always assume that
$\a\le \b$.

We begin with simplest case: $\a$ and $\b$ even with $\chi(\e_1)=\chi(\e_2)=-1$.
In this case, $\Cal T(j)\cap \Cal T(j')=\Cal T(j)$ is just the $\a/2$ ball around the vertex
$[\Lambda(j)]=[\Lambda(j')]$ fixed by both $j$ and $j'$. We have
$$\align
(Z(j)^v,Z(j')^v) &= -(p+1)\a/2 + (1-p)\sum_{r=1}^{\a/2-1}
(\a/2-r)(p+1)p^{r-1}\tag6.2\\
\nass
{}&= -(p+1)\frac{p^{\a/2}-1}{p-1}.\endalign$$

In the case $\a$ and $\b$ even with $\chi(\e_1)=-1$ and $\chi(\e_2)=1$,
$\Cal T(j)\cap \Cal T(j')=\Cal T(j)$ is again just the $\a/2$ ball around the vertex
$[\Lambda(j)]$. Thus,
$$\align
(Z(j)^v,Z(j')^v) &= -(p-1)\a/2 + (1-p)\sum_{r=1}^{\a/2-1}
(\a/2-r)(p+1)p^{r-1}\tag6.3\\
\nass
{}&=\a  -(p+1)\frac{p^{\a/2}-1}{p-1},\endalign$$
where the change from $\chi(\e_2)=-1$ to $\chi(\e_2)=1$ causes the
change in the first term.

The case $\a$ and $\b$ even with $\chi(\e_1)=1$ and $\chi(\e_2)=-1$
is more complicated. Here, the geodesic joining the fixed vertex
$[\Lambda(j')]$ to any given vertex $[\Lambda]$ runs a distance
$\ell$ along the fixed apartment $\Cal B^j$ and then a distance $r$
away from the apartment. The vertices $[\Lambda]$ in $\Cal T(j)$
with $\ell=0$ are all joined to $[\Lambda(j')]$ by a geodesic
emanating along an edge outside of $\Cal B^j$. The contribution of
such vertices is:
$$-(p+1)\a/2 + (1-p) \sum_{r=1}^{\a/2-1} (\a/2-r)(p-1) p^{r-1}
= 1-\a - p^{\a/2}.\tag6.4$$
If $1\le\ell\le (\b-\a)/2$, then vertices $[\Lambda]$ with $r=\a/2$
lie in $\Cal T(j')$. Hence the contribution of vertices with $1\le
\ell\le(\beta -\alpha)/2$ is
$$2(1-p) \sum_{\ell=1}^{(\b-\a)/2}
\bigg( \a/2 + \sum_{r=1}^{\a/2-1}(\a/2-r)(p-1)p^{r-1}\bigg)
=(\a-\b)(p^{\a/2}-1).\tag6.5$$
Next, if $(\b-\a)/2<\ell<\b/2$, then the vertices $[\Lambda]$ with
$r<\b/2$ lie strictly inside $\Cal T(j')$, and contribute
$$2(1-p) \sum_{\ell=(\b-\a)/2+1}^{\b/2-1}
\bigg(\a/2 + \sum_{r=1}^{\b/2-\ell-1}(\a/2-r)(p-1)p^{r-1}\bigg)
=2\a -4\,\frac{p^{\a/2}-1}{p-1}.\tag6.6$$
Finally, the vertices on the boundary of $\Cal T(j')$ contribute
$$\a+2 \sum_{\ell=(\b-\a)/2+1}^{\b/2-1}(\a/2 -(\b/2-\ell))(p-1)p^{\b/2-\ell-1}
=2\,\frac{p^{\a/2}-1}{p-1}.\tag6.7$$
Here the initial $\a$ is the contribution from the two points of
intersection of the apartment $\Cal B^j$  with the boundary of $\Cal T(j')$,
i.e., the two points with $\ell=\b/2$.

The set $\Cal T(j)\cap \Cal T(j')$ in the case $\a$ and $\b$ even
with $\chi(\e_1)=\chi(\e_2)=1$ is identical with the corresponding
set in the previous case! The only change in the formulas occurs in
the first term (6.4), which is now
$$-(p-1)\a/2 + (1-p) \sum_{r=1}^{\a/2-1} (\a/2-r)(p-1) p^{r-1}
=1 - p^{\a/2}\ , \tag6.8$$
due to the change in the contribution of the central vertex. In effect, the
total in this case is $\a$ plus that of the previous case.

In the case $\a$ is even, $\chi(\e_1)=1$ and $\b$ is odd, $\Cal B^{j'}$ is the
midpoint of an edge in the apartment $\Cal B^j$.  Each vertex $[\Lambda]$ is joined
to this fixed midpoint by a unique geodesic, which runs along the apartment
$\Cal B^j$ a distance $\ell+\frac12$ and then a distance $r$ outside the apartment.
The contribution of the vertices with $\ell=0$ is
$$-p \a + 2(1-p)\sum_{r=1}^{\a/2-1} (\a/2-r)(p-1)p^{r-1}
= -\a -2(p^{\a/2}-1),\tag6.9$$
where the initial $-p\a$ is the contribution of the two endpoints of the
edge containing the fixed vertex. When $1\le \ell\le (\b-\a-1)/2$,
we may travel a distance $r=\a/2-1$ and remain strictly inside $\Cal T(j')$.
Thus such vertices contribute
$$2(1-p)\sum_{\ell=1}^{(\b-\a-1)/2} \bigg(\a/2+\sum_{r=1}^{\a/2-1}
(\a/2-r)(p-1)p^{r-1}\bigg)
=(\a-\b+1)(p^{\a/2}-1).\tag6.10$$
Similarly, the vertices strictly inside $\Cal T(j')$ but with
$\ell >(\b-\a-1)/2$ contribute
$$2(1-p) \sum_{\ell=(\b-\a+1)/2}^{(\b-3)/2}
\bigg(\a/2+\sum_{r=1}^{(\b-3)/2-\ell}
(\a/2-r)(p-1)p^{r-1}\bigg)
= 2\a-4\,\frac{p^{\a/2}-1}{p-1}.\tag6.11$$
Finally, the vertices on the boundary of $\Cal T(j')$ contribute
$$\a + 2 \sum_{\ell=(\b-\a+1)/2}^{(\b-3)/2} (\a/2 -((\b-1)/2-\ell))
(p-1)p^{(\b-1)/2-\ell-1}
=2\,\frac{p^{\a/2}-1}{p-1}.\tag6.12$$
Note that here the value $r= (\b-1)/2-\ell$ puts $[\Lambda]$ on the boundary.

The next case is $\a$ odd and $\b$ even with $\chi(\e_2)=1$. Here the whole
ball $\Cal T(j)$ lies entirely inside the tube $\Cal T(j')$ around the
fixed apartment $\Cal B^{j'}$. The total contribution is
$$2(1-p)\sum_{r=0}^{(\a-1)/2} ((\a-1)/2-r) p^r
=\alpha +1-2\,\frac{p^{(\a+1)/2}-1}{p-1},\tag6.13$$
where $d([\Lambda],\Cal B^{j'}) = r+\frac12$.

The last case is $\a$ and $\b$ odd. Here $j$ and $j'$ fix the same midpoint,
and we get
$$-p(\a-1) + 2(1-p)\sum_{r=1}^{(\a-1)/2} ((\a-1)/2-r) p^r
=2 -2\,\frac{p^{(\a+1)/2}-1}{p-1}.\tag6.14$$
Here the term $-p(\a-1)$ comes from the pair of vertices of the edge containing
the fixed vertex. Note that they have multiplicity $(\a-1)/2$ in $Z(j)$.

In each case, we sum the various contributions
to obtain the following table of values of $(Z(j)^v,Z(j')^v)$:

$$
\vcenter{
\begintab{|c|c|c|c|}
\hline
$j\setminus j'$&$\beta$ even&$\beta$ even&$\beta$ odd
\cr
&$\chi(\varepsilon_2)=-1$ & $\chi(\varepsilon_2)=1$ &
\cr
\hline
$\alpha$ even &&&
\cr
$\chi(\varepsilon_1)=-1$& {\bf(i)} &{\bf(ii)}&$\emptyset$
\cr
\hline
$\alpha$ even &&&
\cr
$\chi(\varepsilon_1)=1$ &{\bf(iii)}&{\bf(iv)}&{\bf(v)}
\cr
\hline
&&&
\cr
&&&
\cr
&&&
\cr
 $\alpha$ odd & $\emptyset$ &{\bf(vi)}&{\bf(vii)}
\cr
\hline
\endtab
}
\tag6.15
$$
Here are the values:
$$\align
\text{\bf(i)} &=(6.2)= -(p+1)\, \frac{p^{\a/2}-1}{p-1}\\
\nass
\text{\bf(ii)} &=(6.3)= \a -(p+1)\, \frac{p^{\a/2}-1}{p-1}\\
\nass
\text{\bf(iii)} &=(6.4)+(6.5)+(6.6)+(6.7)= \b+1 - (\b-\a+1)p^{\a/2} - 2\, \frac{p^{\a/2}-1}{p-1}\\
\nass
\text{\bf(iv)} &={\roman{(iii)}}+\alpha= \a+\b+1 - (\b-\a+1)p^{\a/2} - 2\, \frac{p^{\a/2}-1}{p-1}\\
\nass
\text{\bf(v)} &=(6.9)+(6.10)+(6.11)+(6.12)= \b+1 - (\b-\a+1)p^{\a/2} - 2\, \frac{p^{\a/2}-1}{p-1}\\
\nass
\text{\bf(vi)} &=(6.13)= \a+1 -  2\, \frac{p^{(\a+1)/2}-1}{p-1}\\
\nass
\text{\bf(vii)} &=(6.14)= 2 -2\, \frac{p^{(\a+1)/2}-1}{p-1}\\
\endalign$$
Adding these to the expressions in (6.1), we obtain the result
claimed in the Theorem.
\qed\enddemo

\subheading{\Sec7. Intersection numbers and representation densities}

In this section, we will express the intersection number
$(Z(j),Z(j')) = e_p(T)$ given in Theorem 6.1 in terms of
representation densities and their derivatives. The analogous
result in the case of a prime of good reduction is contained in
sections 8 and 14 (Proposition 14.6) of \cite{\annals}. The result
in the present case is somewhat more complicated, and this
reflects, perhaps, the effect of bad reduction.

We begin by reviewing some facts about the representation densities of
quadratic forms, based on section 8 of \cite{\annals} and \cite{\yang}.
For simplicity, we assume that $p\ne2$.

Recall that, for $S\in \Sym_m(\Z_p)$ and $T\in \Sym_n(\Z_p)$
with $\det(S)\ne 0$ and $\det(T)\ne 0$, the classical
representation density is
$$\a_p(S,T) = \lim_{t\rightarrow\infty} p^{-t n(2m-n-1)/2}\,
|\{\, x\in M_{m,n}(\Z/p^t\Z)\ ;\  S[x] - T \in p^t \Sym_m(\Z_p)\,
\}|.\tag7.1$$ Let
$$S_r = \pmatrix S&{}&{}\\{}&1_r&{}\\{}&{}&-1_r\endpmatrix\tag7.2$$
be the orthogonal sum of $S$ and a split space of dimension $2r$. Then there is a
rational function $A_{S,T}(X)$ of $X$ such that
$$\a_p(S_r,T) = A_{S,T}(p^{-r}).\tag7.3$$
The {\it derivative of the representation density} is then
$$\align
\a_p'(S,T) &= \frac{\partial}{\partial X}\bigg(A_{S,T}(X)\bigg)\bigg|_{X=1}\tag7.4\\
\nass
{}&=\frac{\partial}{\partial X}\bigg(\a_p(S_r,T)\bigg)\bigg|_{r=0}.
\endalign
$$

We will be concerned with the case $m=3$ and $n=2$. Let
$$S=-\pmatrix 1&{}&{}\\{}&1&{}\\{}&{}&-1\endpmatrix,\tag7.5$$
and
$$S'=-\pmatrix \eta&{}&{}\\{}&p&{}\\{}&{}&-\eta p\endpmatrix,\tag7.6$$
where $\eta\in \Z_p^\times$ with $\chi(\eta)=(\eta,p)_p=-1$. Here $(a,b)_p$ is the
quadratic Hilbert symbol for $\Q_p$. Note that the form $S$ (resp. $S'$)
is the matrix, with respect to a suitable basis,
of the restriction of the norm form of $M_2(\Q_p)$ (resp.
the  division quaternion algebra $B$ ) to the integral trace zero elements
$V(\Z_p)$ (resp. $V'(\Z_p)$). In particular,
$$1=\det(S) = \det(S')\in \Q_p^\times/\Q_p^{\times,2},\tag7.7$$
and $S$ and $S'$ have opposite Hasse invariants
$$\e_p(S) = -\e_p(S') =  1.\tag7.8$$
Note that this corresponds to the choice $\kappa=-1$ in section 8 of \cite{\annals}.
\medskip\noindent
{\bf Remark:} In earlier sections we have taken the quadratic form
on $V(\Q_p)$ given by $x^2=q(x)\cdot 1_2$, so that $q(x)= -\nu(x)$,
wher $\nu$ is the restriction of the reduced norm to $V$. In
\cite{\annals}, the quadratic form $Q(x)=\nu(x)$ was used. Thus, to
make a consistent link with results of \cite{\annals} we use in
this section and in section 9 below $Q(x)$ rather than $q(x)$. This
change in convention introduces slightly awkward signs at several
points.

Also let
$$S''=-\pmatrix 1&{}&{}\\{}&p&{}\\{}&{}&-p\endpmatrix,\tag7.9$$
so that $S''$ is the matrix for the restriction of the norm form to the trace $0$
elements in the lattice
$$\{\pmatrix a&b\\c&d\endpmatrix\in M_2(\Z_p);\  c \equiv 0 (p)\}.\tag7.10$$

The following Proposition summarizes results of Kitaoka, \cite{\kitaoka},
\cite{\annals}, Corollary 8.4 and 8.5,
in the case of $S$ and results of Myers, \cite{\myers} and T. Yang,
\cite{\yang} in the case of $S'$ and $S''.$

\proclaim{Proposition 7.1} Let $T\in \Sym_2(\Z_p)$ and suppose that $T$ is
$GL_2(\Z_p)$-equivalent to the matrix $\diag(\e_1 p^\a,\e_2 p^\b)$
with $\alpha\le\beta$. Let
$$\mu_p(T) = \cases \chi(-\e_1\e_2)&\text{ if $\a$ and $\b$ are odd}\\
\chi(-\e_2)&\text{ if $\a$ is odd and $\b$ is even}\\
\chi(-\e_1)&\text{ if $\a$ is even and $\b$ is odd}\\
1&\text{ if $\a$ and $\b$ are both even.}
\endcases
$$
Then
$$\a_p(S,T)\ne 0 \iff \mu_p(T) =1,\tag \it i$$
and
$$\a_p(S',T)\ne 0 \iff \mu_p(T)=-1.$$
(ii) (Kitaoka) If $\mu_p(T)=1$, then
$$\a_p(S,T) = (1-p^{-2})\cdot
\cases
p^{\a/2} + 2 \frac{ p^{\a/2}-1}{p-1} &\text{ if $\a$ and $\b$
are even}\\
{}&\text{and $\chi(-\e_1)= -1$}\\
\nass
(\b-\a+1) p^{\a/2} + 2 \frac{p^{\a/2}-1}{p-1}&\text{ if $\a$ is even}\\
{}&\text{and
$\chi(-\e_1)=1$,}\\
\nass
2\frac{p^{(\a+1)/2}-1}{p-1} &\text{ if $\a$ is odd.}
\endcases
$$
(iii) If $\mu_p(T)=-1$, then
$$\a_p'(S,T) = -(1-p^{-2})\cdot\cases \sum_{j=0}^{(\a-1)/2}
(\a+\b-4j)p^j &\text{ if $\a$ is odd.}\\
\nass
\sum_{j=0}^{\a/2-1} (\a+\b-4j)p^j
+\frac12(\b&\!\!\!\!-\a+1)p^{\a/2}\\
{}&\text{if $\a$ is even and $\b$ is odd.}
\endcases
$$
(iv) (Myers) If $\mu_p(T)=-1$, then
$$\a_p(S',T) = 2(p+1).$$
(iv) (Myers, Yang) If $\mu_p(T)=1$, then
$$\align
\a_p'(S',T) &=
-(p+1)(\a+\b+2)\\
\nass
{}&\qquad
+2p\cdot\cases
p^{\a/2} + 2 \frac{ p^{\a/2}-1}{p-1} &\text{ if $\a$ and $\b$
are even}\\
{}&\text{and $\chi(-\e_1)= -1$}\\
\nass
(\b-\a+1) p^{\a/2} + 2 \frac{p^{\a/2}-1}{p-1}&\text{ if $\a$ is even}\\
{}&\text{and
$\chi(-\e_1)=1$,}\\
\nass
2\frac{p^{(\a+1)/2}-1}{p-1} &\text{ if $\a$ is odd.}
\endcases
\endalign$$
(v) (Yang)
If $\mu_p(T)=1$, then
$$1 = \frac{p^2}{p^2-1}\a_p(S,T)-\frac{1}{2(p-1)}\a_p(S'',T).$$
\endproclaim

The last relation follows immediately from Corollary~8.4 in \cite{\yang}. For
an explanation of the dichotomy of (i), see Proposition~1.3 of \cite{\annals}.

The following striking relation is then evident, \cite{\yang}, Theorem~8.1:
\proclaim{Corollary 7.2} If $\mu_p(T)=1$, then
$$\a_p'(S',T) = -(p+1)(\a+\b+2) + \frac{2 p^3}{p^2-1}\, \a_p(S,T).$$
\endproclaim

{\it Remark:} Kitaoka gave an explicit formula
for the representation density $\a_p(S,T)$ when
$T$ is any binary form, $n=2$, and $S$ is unimodular,
$m$ is arbitrary \cite{\kitaoka}. He also handled the case where $n=3$
and $m=4$ \cite{\kitaokaII}. In his thesis, \cite{\myers}, B. Myers gave a formula for $\a_p(S,T)$
in the case of a binary form $T$ where, when diagonalized, the entries of
$S$ have $\text{\rm ord}_p \le 1$. This formula made it possible to compute the
derivative $\a_p'(S',T)$ above.
The relation of Corollary~7.2 was first observed in the thesis of B. Myers,
where the term $-(p+1)(\a+\b+2)$ was erroneously
given as $-(p+1)(\a+3\b+2)$ in certain cases. Using a new method,
Tonghai Yang, \cite{\yang}, found an explicit formula for the
representation density $\a_p(S,T)$ for $T$ a binary form and for $S$
an {\it arbitrary} form. His result (which includes the case
$p=2$) thus allow the computation of any $\a_p'(S,T)$ for a binary $T$.
Recently, Katsurada, \cite{\katsurada}, gave an explicit (very complicated)
formula for $\a_p(S,T)$ when $T$ is arbitrary and $S$ is a sum of hyperbolic planes.
This should allow the computation of the derivatives $\a_p'(S,T)$
in cases $m=2\ell$ and $n=2\ell-1$ in the interesting case in which
$\a_p(S,T)=0$, cf. Proposition~1.3 and Theorem~6.1 of \cite{\annals}.

Comparing the expressions of Proposition~7.1 and Corollary~7.2 with the formula for
the intersection number given in  Theorem~6.1, we obtain the following result.
\proclaim{Theorem 7.3} As in Theorem~6.1, suppose that $j$ and $j'$ are special
endomorphisms with matrix of inner products $T$, defined using the
quadratic form $Q$.\footnote{\ Thus $T$ in this Theorem is the
negative of the matrix attached to $j$ and $j'$ in previous
sections, e.g., in Theorem~6.1.} Then, $\mu_p(T)=1$, and
$$(Z(j),Z(j')) =e_p(-T)
= -\frac{1}{p+1}\, \a_p'(S',T) + \frac{2 p^2}{p+1}\, \a_p(S,T) +\frac{1}{2(p-1)} \a_p(S'',T).$$
\endproclaim

Note that this identity holds for all $T$ and that the coefficients are independent of $T$, so that
there is not as much flexibility as it might at first appear. Also, the occurrence of $e_p(-T)$
is due to our change in the sign of the quadratic form in this section.

As explained in \cite{\annals}, section 8 and Appendix,
the representation densities and their derivatives are closely related to
certain generalized Whittaker functions for the metaplectic cover of
$Sp_2(\Q_p)$. In our present case, let $\ph_p$ be the characteristic function
of $V(\Z_p)$, let $\ph'_p$ be the characteristic function of $V'(\Z_p)$,
and let $\ph''_p$ be the characteristic function of the sublattice of $V(\Z_p)$
where $c$ is divisible by $p$.
Also let $\P_p(s)$, $\P'_p(s)$, and $\P''_p(s)$ be the standard sections of the induced
representation $I_{2,p}(s)$ whose values at $s=0$ are $\l_p(\ph_p)$,
$\l'_p(\ph'_p)$, and $\l_p(\ph''_p)$ respectively. Then, as in Proposition~A.6 of
\cite{\annals}, the values of the associated generalized Whittaker functions at integers $r\ge0$
are
$$\align
W_{T,p}(e,r,\P_p) &= \gamma_p(V_p)\cdot \a_p(S_r,T),\tag7.11\\
\nass
W_{T,p}(e,r,\P'_p)& = p^{-2}\gamma_p(V'_p)\cdot \a_p(S'_r,T),\\
\noalign{and}
\nass
W_{T,p}(e,r,\P''_p) &= p^{-2}\gamma_p(V_p)\cdot \a_p(S''_r,T).\tag7.12
\endalign$$
Note that, by Proposition~A.4 of \cite{\annals},
$$\gamma_p(V_p) = \gamma_p(V'_p) = \gamma_p(-1,\psi_p)^{3} =1,\tag7.13$$
for our standard $\psi_p$ and $p\ne2$.
In our case, since $\mu_p(T)=1$, the quantities of interest to us are then
$$\align
W_{T,p}'(e,0,\P'_p) &= -p^{-2}\,\log(p)\, \a_p'(S',T),\tag7.14\\
\nass
W_{T,p}(e,0,\P_p) &=  \a_p(S,T),\\
\noalign{and}
\nass
W_{T,p}(e,0,\P''_p) & = p^{-2} \a_p(S'',T).\tag7.15
\endalign$$
The relation of Theorem~7.3 can be expressed in terms of the derivative of
a Whittaker function.
\proclaim{Corollary 7.4} Let $A(s)$ and $B(s)$ be rational functions of $p^{-s}$ satisfying
$$A(0)=B(0)=0,$$
and
$$A'(0) = 2 \log(p),\qquad\text{ and }\qquad B'(0) = \frac12\frac{p+1}{p-1}\,\log(p).$$
Define a {\rm nonstandard} section of $I_{2,p}(s)$ by
$$\tilde{\P}_p(s) = \P'_p(s) + A(s)\,\P_p(s) + B(s)\,\P''_p(s).$$
Then
$$\frac{p+1}{p^2}\,\log(p)\,e_p(-T) = W'_{T,p}(e,0,\tilde\P_p).$$
\endproclaim

Note that, since the standard sections $\P_p(s)$, $\P'_p(s)$ and $\P_p''(s)$ are
linearly independent at $s=0$, the conditions on $A(s)$ and $B(s)$ uniquely determine the first
two terms of the Laurent expansion of $\tilde\P_p(s)$ at $s=0$.

For example, one possibility is
$$A(s) = 1-p^{-2s}\qquad\text{ and }\qquad B(s) = \frac14\, \frac{\displaystyle 1+p^{-1-s}}{\displaystyle 1-p^{-1-s}}\,(1-p^{-2s}).\tag 7.16$$

\subheading{\Sec8. Intersection numbers on Shimura curves}

In this section we use the results about intersections of cycles in
Drinfeld space obtained above to
compute the contribution to the height pairing of certain 0-cycles
on Shimura curves of the intersections occuring in the fiber at $p$
for a prime $p$ which ramifies in the quaternion algebra.
Our result extends that of \cite{\annals} where the
case of good reduction was considered. The notation will be slightly
different from that of the earlier sections.

We begin by briefly reviewing the global moduli problem and the
definition of cycles from section 14 of \cite{\annals}.

Let $B$ be an indefinite quaternion algebra over $\Q$ and let $D(B)$ be the
product of the primes which ramify in $B$. Fix a maximal order $O_B$ in $B$.
Let
$$V=\{ x\in B;\  {\roman{tr}}^0(x)=0\},\tag8.1$$
with quadratic form defined by $x^2 = q(x)1_B$, i.e.,
$q(x)=-\nu(x)$ where $\nu$ is the reduced norm. Note that we are
using the negative of the form used in \cite{\annals} and so the
signature of this quadratic form is $(2,1)$. The action of
$H=B^\times$ on $V$ by conjugation, $x\mapsto h x h^{-1}$, induces
isomorphisms $H\simeq GSpin(V)$ and $H/Z\simeq SO(V)$, where $Z$ is
the center of $H$.

Fix a prime $p$ with $p\mid D(B)$,
and a compact open subgroup $K\subset H(\A_f)$ of the form
$K= K_pK^p$ with $K^p\subset H(\A_f^p)$ and $K_p\subset H(\Q_p)$.
We assume that the image $K^p/Z^p$ of $K^p$ in $SO(V)(\A_f^p)$ is torsion
free and that $K_p=(O_B\tt \Z_p)^\times$.

Let $\Cal A_K$ be the functor on $Sch_{\Z_{(p)}}$ which assigns to
a scheme $S$ over $\Z_{(p)}$ the set of isomorphism classes of triples
$(A,\iota,\bar{\eta})$ where
\roster
\item"{(i)}" $A$ is an abelian scheme over $S$, up to prime to $p$ isogeny,
\item"{(ii)}"$\iota:O_B\lra \End_S(A)$ is an embedding satisfying the
special condition \cite{\zink}, p.~17, and
\item"{(iii)}" $\bar{\eta}$
is an equivalence class, modulo right multiplication by $K^p$,
of $O_B$-equivariant isomorphisms $\eta:\hat{V}^p(A) \isoarrow B(\A_f^p)$,
where
$$\hat{V}^p(A) = \prod_{\ell\ne p} T_\ell(A)\tt\Q,$$
is the rational Tate module of $A$.
\endroster
We refer to \cite{\kottwitz} for the precise meaning of this last
condition; in particular, if $S=\Spec\, k$ is the spectrum of a
field, then the equivalence class $\bar{\eta}$ is stable under
$\text{\rm Gal}(\bar{k}/k)$.

As is well known, this functor is represented by a projective
scheme over $\Z_{(p)}$ which we also denote by $\Cal A_K$. The
generic fiber $A_K = \Cal A_K\times_{\Spec\,\Z_{(p)}}
\Spec\,\Q$ is the canonical model of the Shimura curve determined by $B$ and $K$.

Next, we recall the definition of special cycles. For $t\in \Q$, let
$$V_t =\{x\in V;\  q(x)=t\}\tag8.2$$
be the corresponding hyperboloid. For a fixed negative integer $t\in \Z_{(p)}$
and a $K^p$-stable compact open subset $\o\subset V(\A_f^p)$, we
consider the functor $\Cal C(t,\o)$ which assigns to any
$\Z_{(p)}$-scheme $S$ the set of isomorphism classes of collections
$(A,\iota,\bar{\eta},x)$, where $(A,\iota,\bar{\eta})$ is as
before, and where
\roster
\item"{(iv)}" $x\in End_S(A,\iota)$ is an endomorphism with ${\roman{tr}}^0(x)=0$
and such that, for any $\eta\in\bar{\eta}$, the endomorphism $\eta_*(x)$
of $B(\A_f^p)$ is given by right multiplication by an element of $\o$.
\endroster

This functor is representable.
There is a natural morphism $\Cal C(t,\o)\rightarrow \Cal A_K$ given by
omitting the endomorphism $x$. This morphism is finite and unramified.
On the generic fiber, this defines a 0-cycle
$$C(t,\o):=\Cal C(t,\o)\times_{\Spec\, \Z_{(p)}}
\Spec\, \Q \lra A_K,\tag8.3$$
described in more detail in section 10 of \cite{\annals}.

Let $\hak$ be the formal completion of $\Cal A_K$ along its special fiber, and
let
$$\ha=\hakw = \hak \times_{\Spf\,\Z_p}\Spf\, W\tag8.4$$
be the base change of $\hak$ to $W=W(\bar{\Bbb F}_p)$. The
Drinfeld-Cherednik p-adic uniformization of $\ha$ is given as
follows. Fix a base point $\xi_0=(A_0,\iota_0,\bar{\eta}_0)\in
\ha(\bar{\Bbb F}_p)$, and let $B'= \End^0(A_0,\iota_0)$. Then $B'$
is a definite quaternion algebra over $\Q$ whose invariants agree
with those of $B$ at all primes $\ell\ne p$, $\infty$. Let $H' =
B^{\prime,\times}$. Fix $\eta_0\in\bar{\eta_0}$. This choice
induces an identification
$$B'(\A_f^p)=B(\A_f^p)^{\text{\rm op}}\tag8.5$$
determined by the condition $\eta_0(b'v)=\eta_0(v)b$, and hence also
an identification $H'(\A_f^p)=H(\A_f^p)$, with order of multiplication reversed.
Let $\Bbb X = A_0(p)$ be the p-divisible group of $A_0$.
Then $\iota_0$ induces an action of $O_B\tt\Z_p$ on $\Bbb X$ and, by the condition
(ii) above, $\Bbb X$ is a s.f. $O_B\tt\Z_p$ module over $k=\bar{\Bbb F}_p$.
As in section 1, we may also assume that $0$ and $1$ are critical indices of $\Bbb X$
and fix an identification $\End^0(\Bbb X ,\iota_0)=M_2(\Q_p)$.
This gives an identification $H'(\Q_p)=GL_2(\Q_p)$.

Recall from section 1 the category $\text{Nilp}$ of $W$-schemes $S$
such that $p$ is locally nilpotent in $\Cal O_S$ and the notation
$\bar{S}=S\times_{\Spec\, W}\Spec\, \bar{\Bbb F}_p$ for $S\in
\text{\rm Nilp}$. Let $\tha$ be the functor on $\text{\rm Nilp}$
which associates to $S$ the isomorphism classes of collections
$(A,\iota,\bar{\eta},\psi)$ where $(A,\iota,\bar{\eta})$ is as
before, and where
$$\psi: A_0\times_{\Spec\, \bar\Bbb F_p} \bar{S} \lra A\times_S \bar{S}\tag8.6$$
is an $O_B$-equivariant $p$-primary quasi-isogeny.

Let $\Cal M^\bull$ be the functor on $\text{Nilp}$ which associates
to $S\in \text{Nilp}$ the set of isomorphism classes of pairs
$(X,\rho)$ where $X$ is a $p$-divisible group over $S$ and
$$\rho:\Bbb X\times_{\Spec\, \bar{\Bbb F}_p}\bar{S}\lra X\times_S\bar S
\tag8.7$$
is a quasi-isogeny. Then $\Cal M^\bull$ is representable by a formal
scheme and breaks up as a disjoint sum
$$\Cal M^\bull = \coprod_{i\in \Z} \Cal M^i,\tag8.8$$
where $\Cal M^i$ is the locus where the height of $\rho$ is equal
to $p^i$. We fix an element $\Pi\in B'(\Q_p)^\times$ such that
${\roman{tr}}^0(\Pi)=0$, and $\ord_p(\Pi^2)=1$. Using $\Pi$, we may
identify $\Cal M^i$ with $\Cal M=\Cal M^0$ (the Drinfeld space
considered in section 1) via
$$\Cal M\isoarrow \Cal M^i, \qquad (X,\rho)\mapsto (X,\rho\circ \Pi^i).
\tag8.9$$

There is a natural morphism
$$\tha \lra \Cal M^\bull\times H(\A_f^p)/K_p.\tag8.10$$
This morphism is defined as follows. Given
$(A,\iota,\bar{\eta},\psi)\in \tha(S)$, then passing to the
corresponding $p$-divisible groups, we obtain a quasi-isogeny of
$p$-divisible groups
$$\rho(\psi): \Bbb X\times_{\Spec\, {\bar{\Bbb F}_p}}\bar{S}
\lra A(p)\times_S\bar{S},\tag8.11 $$
and $(A(p),\rho(\psi))$ is a point of $\Cal M^\bull(S)$.
Also, for a choice of $\eta\in\bar\eta$ there is a unique element $g\in H(\A_f^p)$
for which the diagram
$$\matrix
\hat{V}^p(A_0)&\overset{\psi_*}\to{\lra}& \hat{V}^p(A)\\
\nass
\eta_0\downarrow &{}&\eta\downarrow\\
\nass
B(\A_f^p)&\overset{R_g}\to{\lra}&B(\A_f^p)\endmatrix\tag8.12$$
commutes, and the coset $gK^p$ is uniquely determined by
$\bar{\eta}$. Note, in (8.12), that $\hat{V}^p(A) =
\hat{V}^p(A\times_S\bar{S})$. An element $\gamma\in H'(\Q)$ acts on
$\tha$ by changing the quasi-isogeny $\psi$ to $\psi\circ
\gamma^{-1}$. Since, for this change in $\psi$, the pair
$((A(p),\rho(\psi)),gK^p)$ changes to
$((A(p),\rho(\psi)\circ\gamma^{-1},\gamma^{-1}gK^p)$, the map
(8.10) is $H'(\Q)$-equivariant.

The theorem of Drinfeld-Cherednik asserts that passage to the
quotient induces an isomorphism of formal schemes over $W$,
\cite{\drinfeld},\cite{\rapoportzink},\cite{\boutotcarayol}:
$$\matrix
\tha&\lra& \Cal M^\bull\times H(\A_f^p)/K^p\\
\nass
\downarrow&{}&\downarrow\\
\nass
\ha^{\phantom{{}\thicksim{}}}&\isoarrow&H'(\Q)\back \bigg(\Cal M^\bull \times H(\A_f^p)/K^p\bigg).
\endmatrix\tag8.13
$$

Now let
$$V'=\{x\in B';\  tr^0(x)=0\},\tag8.14$$
with its quadratic form $q$, which we regard as the space of {\it
special endomorphisms} of $A_0$ in $B'=\End^0(A_0,\iota_0)$.
Associating to any $x\in V'(\Q)$ the corresponding endomorphism of
$\Bbb X=A_0(p)$, we obtain a special endomorphism $j=j(x)$ of $\Bbb
X$, and in fact this induces an identification of $V'(\Q_p)$ with
the space of special endomorphisms of $\Bbb X$. Note that the
notation has shifted; $V'(\Q_p)$ was denoted by $V$ in the earlier
sections. For $t\in \Q$, let
$$V'_t=\{x\in V';\  q(x)=t\}.\tag8.15$$

We may now relate the formal cycle $\hc := \widehat{\Cal
C(t,\o)}\times_{\Spf\, \Z_p}\Spf\, W$ to the cycles in $\Cal M$
considered earlier. Let $\thc = \hc\times_{\ha}\tha$, so that a
point of $\thc$ is a collection $(A,\iota,\bar{\eta},x,\psi)$. From
this data, we obtain an element $\psi^*(x)\in V'_t(\Q)$. By (8.10),
we obtain a map
$$\matrix
\thc &\lra &V'_t(\Q)\times \Cal M^\bull\times H(\A_f^p)/K^p\\
\nass
\downarrow&{}&\downarrow\\
\nass
\tha&\lra& \Cal M^\bull\times H(\A_f^p)/K^p,\endmatrix
\tag8.16$$
and, upon passage to the quotient, an injection
$$\matrix
\hc &\hookrightarrow& H'(\Q)\back\bigg(V'_t(\Q)\times \Cal M^\bull\times H(\A_f^p)/K^p\bigg)\\
\nass
\downarrow&{}&\downarrow\\
\nass
\ha^{\phantom{{}\thicksim{}}}&\isoarrow&H'(\Q)\back \bigg(\Cal M^\bull \times H(\A_f^p)/K^p\bigg).
\endmatrix\tag8.17$$
The image of this map is determined
by the following incidence relations: the point $(x,(X,\rho),gK^p)$
lies in the image if and only if:
\smallskip
\roster
\item"{(i)}" $g^{-1}x g\in\o,$ and
\smallskip
\item"{(ii)}" for $j=j(x)$, $(X,\rho)\in  Z^\bull(j)$.
\endroster
Here $Z^\bull(j)$ is the formal subscheme of $\Cal M^\bull$
defined by the obvious analogue of Definition~2.1.

{\bf Remark 8.1.} Let $Z^i(j) = Z^\bull(j)\cap \Cal M^i$. Then,
under the identification (8.9), we have
$$Z(\Pi^i j\Pi^{-i})=Z^0(\Pi^i j\Pi^{-i}) \isoarrow Z^i(j).\tag8.18$$
In particular, since $\Pi^2$ is central, we have
$$Z(j) \simeq Z^{2i}(j)\qquad \text{and}\qquad Z(j^\vee)
\simeq Z^{2i+1}(j),\tag8.19$$
where $j^\vee=\Pi j \Pi^{-1}$.

Using the fact that $H'(\Q)$ acts transitively on $V'_t(\Q)$ when $t\ne 0$,
we obtain an isomorphism
$$\hc \isoarrow H'(\Q)_x\back\bigg(Z^\bull(j)\times I(x,\o)\bigg),\tag8.20$$
where $x$ is a fixed element of $V'_t(\Q)$, $j=j(x)$, and
$$I(x,\o)=\{ gK^p\in H(\A_f^p)/K^p; g^{-1}x g \in \o\}.\tag8.21$$
Relation (8.20) can be viewed as giving a {\it p-adic uniformization of the
special cycle}, analogous to the uniformization (8.13) of the whole space.

{\bf Remark 8.2.} Assume that the generic fiber $C(t,\o) = \Cal
C(t,\o)\times_{\Spec\,\Z_{(p)}}\Spec\, \Q$ is nonempty. Then $t$
is represented by $V(\Q)$. This follows from the fact that for any
$\C$-valued point $(A,\iota,\bar{\eta})$ of $\Cal A_K$ we have an
injection
$$\End^0(A,\iota)\hookrightarrow B.$$
Since $t$ is represented by $V(\Q)$, it is {\it a fortiori} represented
by $V(\Q_p)$. For the corresponding special endomorphism $j$ in
(8.20), this implies, since $B(\Q_p)$ is a quaternion division
algebra, that $\Q_p(j)$ is a field, i.e., does not split. This is
the global reason, referred to in Remark~3.5, for the observation
in Corollary~3.4, that $Z(j)^h =\emptyset$ if $q(j)=\e p^\a$ with
$\a$ even and $\chi(\e)=1$.

Next, we consider the intersection of
a pair of cycles $\Cal C_1=\Cal C(t_1,\o_1)$ and $\Cal C_2=\Cal C(t_2,\o_2)$, following the
procedure of section 3 of \cite{\krsiegel}. We change notation and
now denote by $\Cal C$ the fiber product
$$\Cal C = \Cal C(t_1,\o_1)\times_{\Cal A_K}\Cal C(t_2,\o_2).\tag8.22$$

Let $S$ be a connected scheme. For a point
$\xi=(A_\xi,\iota,\bar{\eta})\in\Cal A_K(S)$, let
$V_\xi\subset\End^0_S(A_\xi,\iota)$ be the $\Q$-vector space of
special endomorphisms (endomorphisms of trace $0$). This space has
a $\Q$-valued quadratic form defined by $x^2= q_\xi(x)\cdot id_A$.
For a point $\xi=(A_\xi,\iota,\bar{\eta},x_1,x_2)$ of $\Cal C(S)$
the pair  $x_1$, $x_2$ of elements of $V_\xi$ determine a symmetric
matrix (the {\it fundamental matrix} associated to $\xi$)
$$T_\xi =\frac12\pmatrix (x_1,x_1)&(x_1,x_2)\\(x_2,x_1)&(x_2,x_2)\endpmatrix
=\pmatrix t_1&\frac12(x_1,x_2)\\\frac12(x_2,x_1)&t_2\endpmatrix\in \Sym_2(\Q),\tag8.23$$
where $(x,y) = q_\xi(x+y)-q_\xi(x)-q_\xi(y)$ is the bilinear form associated to
$q_\xi$. Note that $\det(T_\xi) = t_1t_2- \frac14 (x_1,x_2)^2$.
A basic fact is that $T_\xi$ is negative semi-definite.
As observed in \cite{\krsiegel}, the function $\xi\mapsto T_\xi$ is
constant on each connected component of $\Cal C$ and there is a decomposition
$$\Cal C =\coprod_T\Cal C_T,\tag8.24$$
where, for $T\in\Sym_2(\Q)$, $\Cal C_T$ is the union of the
components of $\Cal C$ where $T_\xi=T$. Note that the only $T$'s
which actually contribute lie in $\Sym_2(\Z_{(p)})$, are negative
semi-definite, and have diagonal terms $t_1$ and $t_2$ as on the
right side of (8.23). Since the signature of $V(\Q)$ is $(2,1)$,
and using an argument similar to that in Remark~8.2, we obtain:
\proclaim{Lemma 8.3} Suppose that $\det(T)\ne 0$ and hence that $T<0$. Then the image in
$\Cal A$ of the underlying point set of $\Cal C_T$ lies in the
special fiber. Moreover, $\Cal C_T$ is proper over
$\Spec\,\Z_{(p)}$.
\endproclaim
The last statement follows from Corollary~2.13.

If $t_1t_2\notin \Q^{\times,2}$, so that $\det(T)\ne 0$ for all $T$
appearing in the decomposition (8.24), then the discussion of
section 4 provides a well defined intersection number $(\Cal
C_1,\Cal C_2)$, cf. Remark~4.4. As a matter of fact, in that
section, we only considered intersection numbers $(Z,Z')$ for
closed subschemes of a regular two dimensional scheme $X$. The
extension to the case of finite unramified morphisms $Z\rightarrow
X$ and $Z'\rightarrow X$ is immediate. In our case, we obtain
$$\align
(\Cal C_1,\Cal C_2) &= \chi(\Cal C,\Cal O_{\Cal C_1}\tt^{\Bbb L}
\Cal O_{\Cal C_2})\tag8.25\\
\nass
{}&= \sum_T \chi(\Cal C_T,\Cal O_{\Cal C_1}\tt^{\Bbb L} \Cal O_{\Cal C_2}).
\endalign
$$
If $t_1t_2$ is a square, then singular $T$'s can occur in (8.24),
and we define
$$(\Cal C_1,\Cal C_2)^{\text{\rm ns}} :=
\sum_{T, \det(T)\ne0} \chi(\Cal C_T,\Cal O_{\Cal C_1}\tt^{\Bbb L} \Cal O_{\Cal C_2}).
\tag8.26$$

Our next goal is to compute the quantity
$\chi(\Cal C_T,\Cal O_{\Cal C_1}\tt^{\Bbb L} \Cal O_{\Cal C_2})$
for a given nonsingular $T\in Sym_2(\Q)$.
First we pass to the formal schemes
$\hc_1$, $\hc_2$ and
$$\hc = \hc_1\times_{\ha} \hc_2,\tag8.27$$
over $W$.
Here
$$\hc = \coprod_T \hat{\Cal C}_T,\tag8.28$$
where $\hc$ (resp.$\hat{\Cal C}_T$) is the base change to $W$ of
the formal completion of $\Cal C$ (resp. $\Cal C_T$) along its
special fiber. Passing to formal completions and making a formally
\'etale base change leaves the intersection number unchanged, cf.
Remark~4.4. Hence we obtain the following statement.

\proclaim{Lemma 8.4}
Assume that $\det(T)\ne 0$. Then
$$\chi(\Cal C_T,\Cal O_{\Cal C_1}\tt^{\Bbb L} \Cal O_{\Cal C_2})
=\chi(\hat{\Cal C}_T,\Cal O_{\hat\Cal C_1}\tt^{\Bbb L} \Cal O_{\hat\Cal C_2}).$$
\endproclaim

Let ${\hat{\Cal C}_T}^{\sim} = \hat{\Cal C}_T\times_{\ha}\tha$, so
that a point $\xi\in{\hat{\Cal C}_T}^{\sim}(S)$ is a collection
$(A_\xi,\iota,\bar{\eta},x_1,x_2,\psi)$, where $\psi$ is a
quasi-isogeny, as in (8.6). The special endomorphisms $x_1$ and
$x_2$ of $A_\xi$ determine an ordered pair
$\und{x}=[\psi^*(x_1),\psi^*(x_2)]$ of special endomorphisms of
$A_0$, and these, in turn, determine an ordered pair
$\uuj=[j_1,j_2]$ of special endomorphisms of $\Bbb X$ and their
$\Z_p$ span $\j=\Z_p j_1+\Z_p j_2$. Let
$$ V^{\prime,2}_T =\{ x\in V^{\prime,2};\  \frac12(x,x)=T\}.\tag8.29$$
Then $\und{x}$ lies in $V'(\Q)^2_T$. Thus, we obtain an
inclusion
$${\hat{\Cal C}_T}^{\sim} \hookrightarrow V'(\Q)^2_T \times \Cal M^\bull \times H(\A_f^p)/K^p,
\tag 8.30 $$
analogous to (8.16). Again, the point $(\und{x},(X,\rho),gK^p)$
lies in the image of this map if and only if:
\smallskip
\roster
\item"{(i)}" $g^{-1}\und{x} g\in\o_1\times \o_2,$ and
\smallskip
\item"{(ii)}" $(X,\rho)\in  Z^\bull(\j )$, where $Z^\bull(\j )=Z^\bull(j_1)\cap Z^\bull(j_2)$.
\endroster
\smallskip
Since $\det(T)\ne0$, the group $H'(\Q)$ acts transitively on
$V'(\Q)^2_T$, and the stabilizer of a fixed $\und{x}$ is the center $Z'(\Q)$,
which acts trivially on $V'$.
Letting
$$I(\und{x},\o_1\times\o_2)
=\{ gK^p\in H(\A_f^p);\  g^{-1}\und{x} g\in \o_1\times\o_2\}\tag8.31$$
and passing to the quotient,
we obtain an isomorphism
$$\hat{\Cal C}_T \isoarrow Z'(\Q)\back\bigg( Z^\bull(\j )\times I(\und{x},\o_1
\times \o_2)\bigg).\tag8.32$$
Again, this can be viewed as a p-adic uniformization of the component $\hat{\Cal C}_T$ of
the intersection.

We then have isomorphic fiber product diagrams of formal schemes over $W$:
\smallskip
$$
\matrix {\hat{\Cal C}_T}&\overset{\scr\pr_1}\to{\lra}&\hc_1\\
\nass
{\scr\pr_2}\downarrow&{}&\downarrow\\
\nass
\hc_2&\lra&\ha
\endmatrix \tag8.33
$$
\smallskip
and
$$
\matrix
Z'(\Q)\back\bigg(Z^\bull(\j )\times  I(\und{x},\o_1\times \o_2)\bigg)&
\overset{\scr\pr_1}\to{\lra}&
H'(\Q)_{x_1}\back \bigg(Z^\bull(j_1)\times I(x_1,\o_1)\bigg)\\
\nass
{\scr\pr_2}\downarrow&{}&\downarrow\\
\nass
H'(\Q)_{x_2}\back \bigg(Z^\bull(j_2)\times I(x_2,\o_2)\bigg)&\lra&
H'(\Q)\back\bigg( \Cal M^\bull\times H(\A_f^p)/K^p\bigg).
\endmatrix
\tag8.34
$$
\smallskip
Observe that the projection maps $\pr_i$ factor as:
$$\matrix
Z'(\Q)\back \bigg(Z^\bull(\j )\times I(\und{x},\o_1\times \o_2)\bigg)
&\overset{\scr\tilde{\pr}_i}\to{\lra}&
Z'(\Q)\back \bigg(Z^\bull(j_i)\times I(x_i,\o_i)\bigg)\\
\nass
{}&\pr_i\searrow\quad&\downarrow\\
\nass
{}&{}& H'(\Q)_{x_i}\back\bigg(Z^\bull(j_i)\times I(x_i,\o_i)\bigg).
\endmatrix\tag8.35$$

We may assume that the element $\Pi$ chosen above lies in
$B'(\Q)^\times =H'(\Q)$. Let $Z'(\Q)^0$ be the set of elements
$z\in Z'(\Q)$ such that $\ord_p(\det(z))=0$, and note that
$Z'(\Q)^0$ acts trivially on $\Cal M^\bull$. Then
$Z'(\Q)=<\Pi^2>\times Z'(\Q)^0$, and, using the isomorphisms
(8.20), we obtain an identification:
$$Z'(\Q)\back \bigg(Z^\bull(j_i)\times I(x_i,\o_i)\bigg)
\simeq \bigg(Z(j_i)\times Z'(\Q)^0\back I(x_i,\o_i)\bigg) \coprod
\bigg(Z(j_i^\vee)\times Z'(\Q)^0\back I(x_i^\vee,\o_i)\bigg).\tag8.36$$
There is an analogous decomposition
$$\multline
Z'(\Q)\back \bigg(Z^\bull(\j)\times I(\und{x},\o_1\times
\o_2)\bigg)\\
\nass
\simeq
\bigg(Z(\j)\times Z'(\Q)^0\back I(\und{x},\o_1\times \o_2)\bigg)
\coprod
\bigg(Z(\j^\vee)\times Z'(\Q)^0\back I(\und{x}^\vee,\o_1\times \o_2)\bigg)
\endmultline\tag8.37$$
for the right side of (8.32) and the decompositions (8.36) and (8.37) are
compatible with the projections $\tilde{\pr}_i$.

Therefore the restriction of $\pr_i^*\Cal O_{\hc_i}$ to the first
(resp. second) component in the decomposition (8.37) is
$$\tilde{\pr}_i^*\Cal O_{Z(j_i)\times Z'(\Q)^0\back I(x_i,\o_i)}.\tag8.38$$
(resp.
$$\tilde{\pr}_i^*\Cal O_{Z(j_i^\vee)\times Z'(\Q)^0\back
I(x_i^\vee,\o_i)}\quad).\tag8.39$$ Thus we have
$$\align
\chi({\hat{\Cal C}_T},\Cal O_{\hc_1}\tt^{\Bbb L}\Cal O_{\hc_2})
&=\chi(Z(\j ), \Cal O_{Z(j_1)}\tt^{\Bbb L}\Cal O_{Z(j_2)})
\cdot |Z'(\Q)^0\back I(\und{x},\o_1\times\o_2)|\tag8.40\\
\nass
{}&\quad+\chi(Z(\j^\vee ), \Cal O_{Z(j_1^\vee)}\tt^{\Bbb L}\Cal O_{Z(j_2^\vee)})
\cdot |Z'(\Q)^0\back I(\und{x}^\vee,\o_1\times\o_2)|.
\endalign$$

Note that $I(x^\vee,\o)=\Pi\cdot I(x,\o)$ and that passing from
$\und{x}$ to $\und{x}^\vee$ leaves the matrix of inner products of the
components unchanged. Since the intersection number
$$(Z(j_1),Z(j_2))= \chi(Z(\j ), \Cal O_{Z(j_1)}\tt^{\Bbb L}\Cal O_{Z(j_2)})$$
depends only on the matrix $T$ of inner products, we have proved the following:
\proclaim{Theorem 8.5} Assume that $\det(T)\ne 0$. Then
$$
\chi({\hat{\Cal C}_T},\Cal O_{\hc_1}\tt^{\Bbb L}\Cal O_{\hc_2})
=2\cdot\chi\big(Z(\j ), \Cal O_{Z(j_1)}\tt^{\Bbb L}\Cal O_{Z(j_2)}\big)
\cdot \big|Z'(\Q)^0\back I(\und{x},\o_1\times\o_2)\big|.$$
\endproclaim

This result is analogous to Theorem~14.11 of \cite{\annals} and
Theorem~7.2 of \cite{\krsiegel} in that it expresses the
intersection number as a product of a multiplicity and a counting
function. The analysis in our present case is much more elaborate,
however, since the ``multiplicity" is the intersection number
$\chi\big(Z(\j ), \Cal O_{Z(j_1)}\tt^{\Bbb L}\Cal O_{Z(j_2)}\big)$,
which is global on the Drinfeld space!

To obtain the final formula for the intersection number, we express the counting
function as an orbital integral.
Let $\ph_i^p\in  S(V(\A_f^p))$ be the characteristic function of
$\o_i$, so that $\ph_1^p\tt\ph_2^p\in S(V(\A_f^p)^2)$. Then, for
$\und{x}\in V(\A_f^p)^2_T$, the cardinality $|Z'(\Q)^0\back
I(\und{x},\o_1\times\o_2)|$ is given by the orbital integral
$$\multline
|Z'(\Q)^0\back I(\und{x},\o_1\times\o_2)| =
\vol(K^p)^{-1}\,O_T(\ph_1^p\tt\ph_2^p)\\
= \vol(K^p)^{-1}\int_{Z(\A_f^p)\back H(\A_f^p)}
\ph_1^p(g^{-1}x_1 g)\ph_2^p(g^{-1}x_2 g)\, dg.
\endmultline\tag8.41$$
Combining this fact with the explicit formula for $e_p(T) =
(Z(j_1),Z(j_2))$ of Theorem~6.1, which, we recall, depends only on
the $GL_2(\Z_p)$-equivalence class of $T$, we obtain the following
explicit expression for the intersection number.
\proclaim{Theorem 8.6} If $t_1t_2$ is not a square in $\Q^\times$,
then the cycles $\Cal C(t_1,\o_1)$ and $\Cal C(t_2,\o_2)$ meet only in the
special fiber of $\Cal A_K$, and their intersection number is given by
$$(\Cal C(t_1,\o_1),\Cal C(t_2,\o_2)) = 2\cdot
\sum_{T}
e_p(T)\cdot \vol(K^p)^{-1}\,O_T(\ph_1^p\tt\ph_2^p),$$
Here the sum runs over
$$T=\pmatrix t_1&n\\n&t_2\endpmatrix\in \Sym_2(\Z_{(p)}),$$
$O_T(\ph_1^p\tt\ph_2^p)$ is the orbital integral (8.41), and, if
$T$ is $GL_2(\Z_p)$-equivalent to $\text{\rm diag}(\e_1 p^\a,\e_2
p^\b)$, with $0\le \a\le \b$, then
$$e_p(T) = \a+\b+1
-\cases p^{\a/2}+2 \frac{ p^{\a/2}-1}{p-1}
&\text{if $\a$ is even and $\chi(\e_1)=-1$,}\\
\nass
(\b-\a+1)p^{\a/2} +2 \frac{ p^{\a/2}-1}{p-1}
&\text{if $\a$ is even and $\chi(\e_1)=1$,}\\
\nass
2 \frac{ p^{(\a+1)/2}-1}{p-1}
&\text{if $\a$ is odd,}
\endcases$$
as in Theorem~6.1.\hfill\break If $t_1t_2$ is a square, then the
quantity $(\Cal C(t_1,\o_1),\Cal C(t_2,\o_2))^{\text{\rm ns}}$ is
defined by (8.26), and
$$(\Cal C(t_1,\o_1),\Cal C(t_2,\o_2))^{\text{\rm ns}} =
2\cdot\sum_{ T,\
\det(T)\ne0}
e_p(T)\cdot \vol(K^p)^{-1}\,O_T(\ph_1^p\tt\ph_2^p).$$
\endproclaim

\subheading{\Sec9. Intersection numbers and Fourier coefficients}

In this section, we combine the results of sections 7 and 8 to prove
a relation between the intersection numbers of special cycles and the
Fourier coefficients of the derivative of a certain Siegel Eisenstein
series. This relation extends Theorem~14.11 of \cite{\annals} to the
primes of bad reduction $p\mid D(B)$ of the Shimura curve $\Cal A_K$.
We refer to section 7 of \cite{\annals}
for more details concerning the definition of the incoherent Siegel
Eisenstein series and its derivative.

We continue to use the notation of section~8, with one exception.
To make our results consistent with those of \cite{\annals}, we
will change the sign of the quadratic form on $V$, i.e., we now
take the quadratic form $Q(x) = \nu(x)$. Thus, a special cycle will
be associated to data $(t,\o)$ where $t\in \Q^\times_{>0}$ with
$\ord_p(t)\ge 0$ and where $\o$ is a $K^p$-stable compact open
subset of $V(\A_f^p)$. We also assume that $\o$ is locally
centrally symmetric.

For a pair of cycles $\Cal C_1 =\Cal C(t_1,\o_1)$ and $\Cal C_2 = \Cal
C(t_2,\o_2)$,
with $t_1t_2$ not a square in $\Q^\times$, let
$$<\Cal C_1,\Cal C_2>_p  := \vol(K)\,\log(p)\,\big(\Cal C_1,\Cal
C_2\big)\tag9.1$$
be the $p$ part of their height pairing.
Here $K=K_pK^p$ is the
compact open subgroup of section~8, and the intersection number
$\big(\Cal C_1,\Cal C_2\big)$ is as in (8.25).
If
$t_1t_2$ is a square, let
$$<\Cal C_1,\Cal C_2>^{\text{ns}}_p  := \vol(K)\,\log(p)\,\big(\Cal C_1,
\Cal C_2\big)^{\text{ns}}\tag9.2$$
be the `nonsingular part' of the height pairing, defined using (8.26).
Of course, (9.2)
simply reduces to (9.1) when $t_1t_2$ is not a square.
Note that the intersection number $(\Cal C_1,\Cal C_2)$
is taken on the quotient $\Cal A_K$ and hence
depends on the choice of $K$. On the
other hand, due to the factor $\vol(K)$, the quantity $<\Cal C_1,\Cal C_2>_p$
is independent of $K$, but depends on the choice of Haar measure on
$H(\A_f)$ used to calculate $\vol(K)$.
Here we fix the Tamagawa measure $dh$ on $H(\A)$ and a factorization
of this measure as $dh = d_\infty h\cdot d_fh$
as in \cite{\annals}, p.~573. The Haar measure $d_fh$ on $H(\A_f)$ is used to
compute $\vol(K)$. In addition, measures $d_\ell h$ are fixed for all $\ell$.

Next we introduce the relevant Eisenstein series. Let
$\ph_1^p$ and $\ph_2^p\in S(V(\A_f^p))^{K^p}$ be the characteristic functions of
the sets $\o_1$ and $\o_2$.
Let $\P^p_f(s)$ be the standard section of the induced representation
$I_{2,\A_f^p}(s)$
with $\P^p_f(0)=\l_f^p(\ph_1^p\tt\ph_2^p)$. Let $\tilde{\P}_p(s)$ be the {\it
nonstandard}
section of $I_{2,p}(s)$ defined in Corollary~7.4, and let
$\P_\infty^{\frac32}(s)$
be the standard section of weight $\frac32$ as in (7.14) of \cite{\annals}. Then
$$\P(s) = \P_\infty^{\frac32}(s)\tt\tilde{\P}_p(s)\tt\P_f^p(s)\tag9.3$$
is an incoherent section of the global induced representation of the metaplectic
group $G''_\A$ of genus $2$. Let $E(g'',s,\P)$ be the associated incoherent
Eisenstein series.

Let $G'_\A$ be the metaplectic cover of $Sp_1(\A)$ and recall that there is a
homomorphism $\iota:G'_\A\times G'_\A\rightarrow G''_\A$. Restricting this to the
real points,
let $g''=\iota(g_1',g_2')$, where $g'_1$ and $g'_2\in G'_\R$.
For $g'\in G'_\R$, let $W_t^{\frac32}(g')$ be the holomorphic Whittaker function
of Proposition~7.3 of \cite{\annals}.

\proclaim{Theorem 9.1} For $g'_1$ and $g'_2\in G'_\R$,
and with the notation just introduced,
$$2\pi^2\,W_{t_1}^{\frac32}(g'_1)\,W_{t_2}^{\frac32}(g'_2)\, <\Cal C_1,\Cal
C_2>_p^{\text{ns}} \
=\  \sum_{T} E'_T(\iota(g_1',g_2'),0,\P),$$
where the sum is on positive definite $T\in \Sym_2(\Z_{(p)})$ with
$$T=\pmatrix t_1&*\\*&t_2\endpmatrix$$
and with $\mu_p(T)=1$. Here the invariant $\mu_p(T)$
is defined in Proposition~7.1.
\endproclaim
\demo{Remark} The condition $\mu_p(T)=1$ implies that $T$ is
{\it not} represented by the quadratic form $Q$ on $V(\Z_p)$.
\enddemo

Note that this result is consistent with Theorem~14.11, for $p\nmid D(B)$, and
Theorem~12.6, for $p=\infty$, of \cite{\annals}. Thus the result of section 15 of
that
paper can be extended correspondingly.

\demo{Proof}
For $T\in Sym_2(\Q)$, with $\det(T)\ne0$, and for $\text{Re}(s)$ sufficiently
large,
the $T$\snug-th Fourier coefficient
of the incoherent Eisenstein series $E(g'',s,\P)$
has a product formula
$$E_T(g'',s,\P) = W_{T,\infty}(g'',s,\P^{\frac32}_\infty)\cdot
W_{T,p}(e,s,\tilde{\P}_p)\cdot
\prod_{\ell\ne p}W_{T,\ell}(e,s,\P_\ell).\tag9.4$$
We will write
$$W_{T}(s,\P_f^p)=\prod_{\ell\ne p}W_{T,\ell}(e,s,\P_\ell),\tag9.5$$
and we recall that this function, initially defined for
$\text{Re}(s)$ sufficiently large, has an entire analytic
continuation, \cite{\annals}, p.~562. We will only be interested in
those $T$ for which $W_{T,p}(e,0,\tilde{\P}_p)=0$, so that for the
derivative
$$E'_T(g'',0,\P) =  W_{T,\infty}(g'',0,\P^{\frac32}_\infty)
\cdot W'_{T,p}(e,0,\tilde{\P}_p)\cdot W_{T}(0,\P_f^p).\tag9.6$$

We would like to prove the identity
$$\align &E'_T(\iota(g_1',g_2'),0,\P)\tag9.7\\
\nass
{} &= C\cdot W_{t_1}^{\frac32}(g'_1)\,W_{t_2}^{\frac32}(g'_2)\,
\vol(K)\,\log(p)\, (\Cal C_1,\Cal C_2)_{-T}\\
\nass
{}&= C\cdot W_{t_1}^{\frac32}(g'_1)\,W_{t_2}^{\frac32}(g'_2)\,
\vol(K)\,\log(p)\, 2 e_p(-T)\cdot \vol(K^p)^{-1} O_T(\ph_1^p\tt\ph_2^p),
\endalign
$$
with constant $C=2\pi^2$. By (7.35) of \cite{\annals}, we have
$$W_{T,\infty}(g'',0,\P^{\frac32}_\infty) = W_{t_1}^{\frac32}(g'_1)\,
W_{t_2}^{\frac32}(g'_2),\tag9.8$$ while, by Corollary~7.4 above,
$$W'_{T,p}(e,0,\tilde{\P}_p) = \frac{p+1}{p^2} \log(p)\, e_p(-T).\tag9.9$$
Thus, substituting these in (9.7), we find that the desired identity is:
$$\frac{p+1}{p^2}\cdot
W_{T}(0,\P_f^p) = 2 C \cdot \vol(K_p)\cdot O_T(\ph_f^p)\ \
,\tag9.10$$ where $\varphi_f^p=\varphi_1^p\otimes \varphi_2^p$.
\par\noindent
An easy calculation shows that $\vol(K_p)= (p+1)/p^2$, so that it
remains to show that
$$W_{T}(0,\P_f^p)=2 C\cdot O_T(\ph_f^p).\tag9.11$$
This last identity can be derived from the Siegel-Weil formula as
follows.

Let $B'$ be the definite quaternion algebra defined in section 8
above and let $V'$ be the space of trace $0$ elements in $B'$ with
quadratic form $Q(x) = \nu(x)$ given by the restriction of the
reduced norm. The identification (8.5) gives an identification
$$V'(\A_f^p) =V(\A_f^p).\tag9.12$$
Let $\ph'=\tt_v\ph'_v\in S(V'(\A)^2)$ be the factorizable, locally
even,  global Schwartz function defined by
$$\ph'_v(x) = \cases e^{-\pi \text{\rm tr}(Q(x))}&\text{ if $ v=
\infty$,}\\
\nass
\ph_{1,\ell}\tt\ph_{2,\ell}(x)&\text{ if $v=\ell\ne p$,}\\
\nass
\ph'_p(x)&\text{ if $v=p$,}
\endcases\tag9.13
$$
for some (for the moment arbitrary) $\ph'_p\in S(V'(\Q_p)^2)$.
If $T\in\Sym_2(\Q)$ with $\det(T)\ne 0$ is represented by $V'$, then,
by (7.28) of \cite{\annals},
the $T$-th Fourier coefficient of the theta integral of $\ph'$ is given by
$$\align
I_T(g'',\ph') &= \frac12 \int_{Z'(\A)H'(\Q)\back H'(\A)}
\sum_{x\in V'(\Q)^2\atop Q(x) = T} \o_\psi(g'')\ph'(h^{-1}x) \, dh\\
\nass
{}&=\frac12 \int_{Z'(\A)\back H'(\A)} \o_\psi(g'')\ph'(h^{-1}x_0)
 \, dh\tag9.14\\
\nass
{}&=\frac12\cdot O_T(\o_\psi(g'')\ph'_\infty)\cdot O_T(\ph'_p)\cdot O_T(\ph_f^p).
\endalign
$$
Here $x_0\in V({\Bbb Q})^2$ is an arbitrary base point, and the
orbital integrals are formed as in (8.41). We take $\varphi'_p$
even so that $\varphi'$ is locally even.

Note that the individual factors here
depend on the choice of the Haar measures, which are fixed as above.
By the Siegel-Weil formula for the coherent Eisenstein series associated to
$\l'(\ph')$,
$$E_T(g'',0,\P')=2\cdot I_T(g'',\ph'),\tag9.15$$
i.e.,
$$W_{T,\infty}(g'',0,\P^{\frac32}_\infty)\cdot W_{T,p}(e,0,\P'_p)\cdot
W_{T}(0,\P_f^p)
=O_T(\o_\psi(g'')\ph'_\infty)\cdot O_T(\ph'_p)\cdot O_T(\ph_f^p). \tag9.16$$
By (7.33)--(7.35) of \cite{\annals}:
$$O_T(\o_\psi(g'')\ph'_\infty) = (2\pi)^2\, W_{t_1}(g'_1)
W_{t_2}(g'_2),\tag9.17$$
so that (9.16) becomes
$$W_{T,p}(e,0,\P'_p)\cdot
W_{T}(0,\P_f^p)
=(2 \pi)^2\cdot O_T(\ph'_p)\cdot O_T(\ph_f^p).\tag9.18$$
We now choose $\ph'_p$ so that $W_{T,p}(e,0,\P'_p)\ne0$, and we
compute the ratio.
\enddemo

\proclaim{Lemma 9.2} $$\frac{O_T(\ph'_p)}{W_{T,p}(e,0,\P'_p)}=1.$$
\endproclaim

\demo{Proof} We only sketch the argument which consists of two
steps. In the first step one proves that the ratio on the left hand
side of Lemma 9.2 is independent of $T\in {\roman{Sym}}_2({\Bbb
Q}_p)$ (with det $(T)\neq 0$) and of $\varphi'_p\in S(V'({\Bbb
Q}_p)^2)$. In the second step one calculates the ratio by making a
special choice of $T$ and $\varphi'_p$. Namely, if $T_0=-1_2$ and
$\varphi_0={\roman{char}}\, V'(\Z_p)^2$, then
$$O_{T_0}(\varphi_0)= {\roman{vol}}(K'_p)\ \ .\tag9.19$$
With the measures as described on p.\ 573 of \cite{\annals} the
naive volume of $K'_p$, i.e., without the convergence factor
$\lambda_p=(1-p^{-1})$ is
$$p^{-4}(p^2-1)(p^2-p)\ \ ,$$
so that, dividing by the convergence factor the right hand side of
(9.19) is $1-p^2$.
\hfill\break
On the other hand, by Kitaoka's formula (cf.\ Proposition 8.3 of
\cite{\annals} and noting that $\gamma(V'(\Q_p))=1$)
$$W_{T_0}(e, 0, \Phi_0)= \alpha_p(S,T_0)=1-p^{-2}\ \ .\qquad\qed$$
\enddemo

In particular,
$$W_{T}(0,\P_f^p)
=(2\pi)^2\cdot O_T(\ph_f^p),$$
and the proof of Theorem~9.1 is complete.

\redefine\vol{\oldvol}

\subheading{References}
\widestnumber\key{666}

\ref\key{\AM}
\by M. F. Atiyah and I. G.  MacDonald
\book Introduction to Commutative Algebra
\yr 1969
\publ Addison-Wesley
\publaddr Reading, MA
\endref

\ref\key{\boutotcarayol} \by J.-F. Boutot and H. Carayol
\paper Uniformisation p-adique des courbes de Shimura: les th\'eor\`emes
de Cerednik et de Drinfeld
\inbook in: Courbes modulaires et courbes de Shimura
\publ Ast\'erisque
{\bf 196--197}
\yr 1991
\pages 45--158
\endref

\ref\key{\deligne}
\by P. Deligne
\paper Intersections sur les surfaces r\'eguli\`eres, expos\'e X
\inbook in P.\ Deligne and N.\ Katz, SGA 7, II: Groupes de monodromie en
g\'eom\'etrie alg\'ebrique, SLN 340, Springer 1973, 1-38
\endref

\ref\key{\delignetwo}
\by P.\ Deligne
\paper La classe de cohomologie associ\'ee \`a un cycle par A.\
Grothendieck
\inbook in P.\ Deligne: SGA $4{1\over 2}$, SLN 569, Springer
1973, 129-153
\endref

\ref\key{\drinfeld}
\by V. G. Drinfeld
\paper Coverings of p-adic symmetric regions
\jour Funct. Anal. Appl.
\vol 10
\yr 1977
\pages 29--40
\endref

\ref\key{\genestierbook}
\by  A. Genestier
\book Espaces Sym\'etriques de Drinfeld
\publ Ast\'erisque {\bf 234}
\yr 1996
\endref

\ref\key{\genestier}
\by  A. Genestier
\paper Letter to M. Rapoport
\jour August 12, 1996
\endref

\ref\key{\gross}
\by B. H. Gross
\paper On canonical and quasi-canonical liftings
\jour Inventiones math.
\vol 84
\yr 1986
\pages 321--326
\endref

\ref\key{\grosskeating}
\by B. H. Gross and K. Keating
\paper On the intersection of modular correspondences
\jour Inventiones math.
\vol 112
\yr 1993
\pages 225--245
\endref

\ref\key{\katsurada}
\by H. Katsurada
\paper An explicit formula for the Fourier coefficients of Siegel-Eisenstein series
\jour preprint
\yr 1997
\endref

\ref\key{\kitaoka}
\by Y. Kitaoka
\paper A note on local densities of quadratic forms
\jour Nagoya Math. J.
\vol 92
\yr 1983
\pages 145--152
\endref

\ref\key{\kitaokaII}
\by Y. Kitaoka
\paper Fourier coefficients of Eisenstein series of degree 3
\jour Proc. Japan Acad.
\vol 60
\yr 1984
\pages 259--261
\endref

\ref\key{\kottwitz}
\by R. Kottwitz
\paper Points on some Shimura varieties over finite fields
\jour JAMS
\vol 5
\yr 1992
\pages 373--444
\endref

\ref\key{\kottwitztwo}   
\by R. Kottwitz
\paper Calculation of some orbital integrals
\inbook R.P.\ Langlands and D.\ Ramakrishnan (ed.), {\it The zeta
functions of Picard modular surfaces,} Publ.\ CRM. Montreal (1992),
349-362
\endref

\ref\key{\annals}
\by S. Kudla
\paper Central derivatives of Eisenstein series and height pairings
\jour Annals of Math.
\yr 1997
\vol 146
\pages 545--646
\endref

\ref\key{\krsiegel}
\by S. Kudla and M. Rapoport
\paper Cycles on Siegel 3-folds and derivatives of Eisenstein series
\jour preprint
\yr 1997
\endref

\ref\key{\krHB}
\by S. Kudla and M. Rapoport
\paper Arithmetic Hirzebruch-Zagier cycles
\jour in preparation
\endref

\ref\key{\myers}
\by B. Myers
\paper Local representation densities of non-unimodular quadratic forms
\jour thesis, University of Maryland
\yr 1994
\endref

\ref\key{\rapoportzink}
\by M. Rapoport and T. Zink
\book Periods of p-divisible groups
\bookinfo Annals of Math. Studies {\bf 141}
\publ Princeton U. Press
\publaddr Princeton, NJ
\yr 1996
\endref

\ref\key{\yang}
\by Tonghai Yang
\paper An explicit formula for local densities of quadratic forms
\jour to appear in J. Number Theory
\endref

\ref\key{\zink}
\by T. Zink
\paper \"Uber die schlechte Reduktion einiger Shimuramannigfaltigkeiten
\jour Compositio Math.
\vol 45
\yr 1981
\pages 15--107
\endref
\vskip3 cm
\noindent
\line{Stephen S.\ Kudla\hfill Michael Rapoport}
\line{Department of Mathematics\hfill Mathematisches Institut}
\line{University of Maryland\hfill der Universit\"at zu K\"oln}
\line{College Park, MD 20742\hfill Weyertal 86-90}
\line{{}\hfill D -- 50931 K\"oln}
\line{USA\hfill Germany}

\bye

%% file: drinfigmittext.tex
\input psfig.sty
\centerline{
\psfig{figure=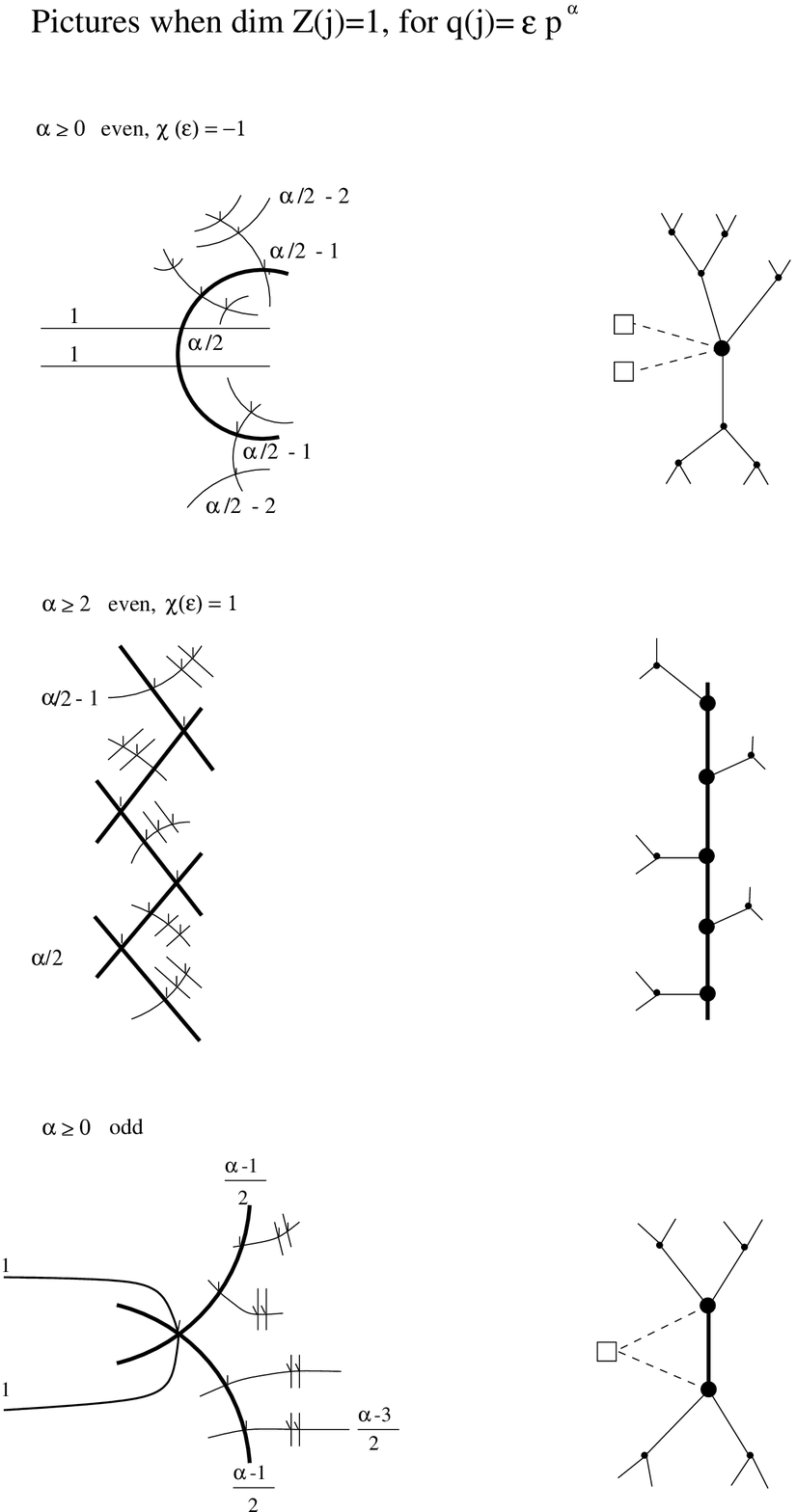,height=14.0cm}
}
\vskip1true cm\noindent
{\eightbf Explanation:} 
{\eightrm On the left the irreducible components of the
special cycles are depicted; on the right the dual graph is given.
The ``central'' vertical irreducible components (1 in the first
case, $\infty$ in the second case, 2 in the third case) are in bold
face. They occur with multiplicity $\left\lbrack
\frac{\alpha}{2}\right\rbrack$ in $Z(j)$. The other vertical
components have multiplicity $\left\lbrack
\frac{\alpha}{2}\right\rbrack -i$, where $i$ is the ``distance''
from a central component. Also the horizontal components are given;
they are indicated by a box on the right (2 in the first case, none
in the second case, 1 in the third case). The little barbs indicate
embedded components.
}

%% file: tabelle.tex
\catcode`\_=11%
\def\my_self#1{\def#1{\noexpand#1}}%
\my_self\Uebersetze_Ende
\my_self\Uebersetzung_Ende
\newtoks\Spalten_Format
\newtoks\Spalteneintrag_Format
\newtoks\Spaltenzwischenraum_Format
\newtoks\Spaltenende_Format
\newtoks\Vielfach_Format
\newif\if_vor_erster_Spalte_
\def\_vor_erster_Spalte_true{\global\let\if_vor_erster_Spalte_=\iftrue}%
\def\_vor_erster_Spalte_false{%
 \global\let\if_vor_erster_Spalte_=\iffalse
}%
\newif\if_Mathe_
\newif\if_Mathe_normal_
\newcount\ne_Zahl
\newcount\noch_ne_Zahl
\newdimen\tabrulewidth
\newdimen\hfree
\newdimen\vfree
\newskip\Init_tabskip
\def\Frei_raum#1{%
 \vbox{%
  \vskip\vfree
  \vtop{%
   \hbox{%
    \hskip\hfree
    \ignorespaces
    #1%
    \unskip
    \hskip\hfree
   }%
   \vskip\vfree
  }%
 }%
}%
\def\wider#1{%
 \relax
 \ifmmode
  \def\next{\Frei_raum{$#1$}}%
 \else
  \def\next{\Frei_raum{#1}}%
 \fi
 \next
}%
\tabrulewidth=.4pt
\hfree=2pt
\vfree=2pt
\newskip\par_width
\def\par_box#1#2{%
 \par_width=#1%
 \hbox to \par_width{%
  \vtop{%
   \normalbaselines
   \noindent
   \leftskip=0pt%
   \rightskip=-\par_width
   \advance\rightskip by\hsize
   #2%
  }%
  \hss
 }%
}%
\def\verknuepfe_Tokens_zur_ner_Liste#1#2#3{\global #3={#1#2}}%
\def\fuege_Tokens_an#1#2#3{%
 \expandafter\verknuepfe_Tokens_zur_ner_Liste
 \expandafter{\the#2}{#1}#3%
}%
\def\verknuepfe_Token_Listen#1#2#3{%
 \expandafter\fuege_Tokens_an
 \expandafter{\the#2}#1#3%
}%
\def\fuege_Tokens_ans_Spaltenformat#1{%
 \fuege_Tokens_an
 {#1}%
 \Spalten_Format
 \Spalten_Format
}%
\def\fuege_Tokens_ans_Vielfachformat#1{%
 \fuege_Tokens_an
 {#1}%
 \Vielfach_Format
 \Vielfach_Format
}%
\def\Spalteneintrag_hfil{\fuege_Tokens_ans_Spaltenformat\hfil}%
\def\Mathe_zu_dollar{%
 \if_Mathe_
  \fuege_Tokens_ans_Spaltenformat{\relax$\relax}%
 \fi
}%
\def\Spalten_Eintrag{%
 \Mathe_zu_dollar
 \verknuepfe_Token_Listen
  \Spalten_Format
  \Spalteneintrag_Format
  \Spalten_Format
 \Mathe_zu_dollar
}%
\def\Uebersetze_#1{%
 \def\Uebersetzung_##1#1##2##3\Uebersetzung_Ende{##2}%
 \Uebersetzung_
  +\Uebersetze_plus
  |\Uebersetze_vrule
  *\Uebersetze_star
  @\Uebersetze_at
  c\Uebersetze_c
  d\Uebersetze_d
  h\Uebersetze_h
  l\Uebersetze_l
  m\Uebersetze_m
  p\Uebersetze_p
  r\Uebersetze_r
  \Uebersetze_Ende\empty
  0\Uebersetze_zero
  \Uebersetzung_Ende
}%
\def\Uebersetze_Eintrag#1{%
 \if_vor_erster_Spalte_
  \_vor_erster_Spalte_false
 \else
  \verknuepfe_Token_Listen
   \Spalten_Format
   \Spaltenzwischenraum_Format
   \Spalten_Format
 \fi
 #1%
 \Uebersetze_
}%
\def\Uebersetze_l{%
 \Uebersetze_Eintrag{%
  \Spalten_Eintrag
  \Spalteneintrag_hfil
 }%
}%
\def\Uebersetze_r{%
 \Uebersetze_Eintrag{%
  \Spalteneintrag_hfil
  \Spalten_Eintrag
 }%
}%
\def\Uebersetze_c{%
 \Uebersetze_Eintrag{%
  \Spalteneintrag_hfil
  \Spalten_Eintrag
  \Spalteneintrag_hfil
 }%
}%
\def\Uebersetze_vrule{%
 \fuege_Tokens_ans_Spaltenformat\vline
 \Uebersetze_
}%
\def\vline{%
 \vrule width\tabrulewidth\relax
 \futurelet\next\v_line
}%
\def\v_line{%
 \ifx\next\vline
  \hskip\hfree\relax
 \fi
}%
\def\Uebersetze_at#1{%
 \fuege_Tokens_ans_Spaltenformat{#1}%
 \Uebersetze_
}%
\def\Uebersetze_p#1{%
 \Uebersetze_Eintrag{%
  \expandafter\Uebersetze_p_hilf
  \expandafter{\the\Spalteneintrag_Format}%
  {#1}%
 }%
}%
\def\Uebersetze_p_hilf#1#2{%
 \fuege_Tokens_ans_Spaltenformat
  {\wider{\par_box{#2}{\un_wider#1}}\hfil}%
}%
\def\un_wider{\futurelet\next\skip_wider}%
\def\skip_wider{%
 \ifx\next\wider
  \let\next\eat_wider
 \else
  \let\next\relax
 \fi
 \next
}%
\def\eat_wider#1#2{\un_wider#2}%
\def\Uebersetze_d{%
 {%
  \Uebersetze_ r0@{\hskip-2\hfree}l\Uebersetze_Ende
 }%
 \Uebersetze_
}%
\def\Uebersetze_star#1#2{%
 {%
  \ne_Zahl=#1%
  \loop
  \ifnum\ne_Zahl>0
   \advance \ne_Zahl by-1
   \Uebersetze_ #2\Uebersetze_Ende
  \repeat
 }%
 \Uebersetze_
}%
\def\Uebersetze_zero{\Uebersetze_tabskip{0pt}}%
\def\Uebersetze_plus{\Uebersetze_tabskip{0pt plus 1fil}}%
\def\Uebersetze_tabskip#1{%
 \if_vor_erster_Spalte_
  \Init_tabskip=#1\relax
 \else
  \fuege_Tokens_ans_Spaltenformat{\tabskip=#1\relax}%
 \fi\Uebersetze_
}%
\def\Uebersetze_h{\expandafter\global\_Mathe_false\Uebersetze_}%
\def\Uebersetze_m{\expandafter\global\_Mathe_true\Uebersetze_}%
\def\Beginn_der_Tabelle#1#2#3#4{%
 \_vor_erster_Spalte_true
 \global\Spalten_Format={#1}%
 \if_Mathe_normal_
  \Uebersetze_ 0m#4\Uebersetze_Ende
 \else
  \Uebersetze_ 0h#4\Uebersetze_Ende
 \fi
 \verknuepfe_Token_Listen
  \Spalten_Format
  \Spaltenende_Format
  \Spalten_Format
 \v_box\bgroup
  \offinterlineskip
  \tabskip=\Init_tabskip
  #2%
  \halign #3\bgroup
   \span\the\Spalten_Format
}%
\def\Ende_der_Tabelle{\crcr\egroup\egroup}%
\def\beginmatrix{%
 \global\let\v_box=\vcenter
 \expandafter\global\_Mathe_true
 \global\Spalteneintrag_Format={\wider{##}}%
 \global\Spaltenzwischenraum_Format={&}%
 \global\Spaltenende_Format={\cr}%
 \null
 \,%
 \_Mathe_normal_true
 \Beginn_der_Tabelle{}{\mathsurround=0pt}{}%
}%
\def\endmatrix{\Ende_der_Tabelle\,}%
\def\begintab{%
  \global\let\v_box=\vtop
  \expandafter\global\_Mathe_false
  \global\Spalteneintrag_Format={\wider{##}}%
  \global\Spaltenzwischenraum_Format={&}%
  \global\Spaltenende_Format={\cr}%
  \_Mathe_normal_false
  \Beginn_der_Tabelle{}{}{}%
}%
\let\endtab=\Ende_der_Tabelle
\def\beginfixtab#1{%
 \global\let\v_box=\vtop
 \expandafter\global\_Mathe_false
 \global\Spalteneintrag_Format={\wider{##}}%
 \global\Spaltenzwischenraum_Format={&}%
 \global\Spaltenende_Format={\cr}%
 \_Mathe_normal_false
 \Beginn_der_Tabelle{}{}{to #1}%
}%
_der_Tabelle
\def\hline{%
 \noalign{%
  \unskip
  \hrule height\tabrulewidth
  \vskip\vfree
  \vskip-\vfree
 }%
}%
\def\span_omit{\fuege_Tokens_ans_Vielfachformat{\span\omit}}%
\def\und_omit{\fuege_Tokens_ans_Vielfachformat{&\omit\unskip}}%
\def\c_line[#1-#2]{%
 \noalign{\vskip-\tabrulewidth}%
 \omit
 \ne_Zahl=#1%
 \noch_ne_Zahl=#2%
 \global\Vielfach_Format={}%
 \advance \ne_Zahl by -1%
 \advance \noch_ne_Zahl by -\ne_Zahl%
 \loop
 \ifnum \ne_Zahl>1
  \advance\ne_Zahl by-1
  \span_omit
 \repeat
 \ifnum \ne_Zahl>0
  \und_omit
 \fi
 \fuege_Tokens_ans_Vielfachformat{%
  \leaders\hrule height\tabrulewidth\hfill}%
 \ifnum \noch_ne_Zahl>1
  \advance \noch_ne_Zahl by-1
  \span_omit
 \repeat
 \the\Vielfach_Format
 \crcr
}%
\def\multicolumn#1#2#3{%
 \multispan{#1}%
 {\_vor_erster_Spalte_true
  \global\Spaltenzwischenraum_Format={}%
  \global\Spalteneintrag_Format={#3}%
  \global\Spalten_Format={}%
  \if_Mathe_normal_
   \Uebersetze_ m#2\Uebersetze_Ende
  \else
   \Uebersetze_ h#2\Uebersetze_Ende
  \fi
  \the\Spalten_Format
 }%
 \ignorespaces
}%
\catcode`\_=8